\documentstyle[12pt]{article}
\newcommand{\al}{{\alpha}}
\newcommand{\be}{{\bar \varepsilon}}
\newcommand{\tp}{{\tilde p}}
\newcommand{\bd}{{\bar d}}
\newcommand{\Epip}{{({\cal E}, \Pi, p)}}
\newcommand{\Eopip}{{({\cal E}_0, \Pi_0, p_0)}}
\newcommand{\Eopipgamma}{{({\cal E}_0^\gamma, \Pi_0^\gamma, p_0^\gamma)}}
\newcommand{\cE}{{\cal E}} 
\newcommand{\tcE}{\tilde {\cal E}}
\newcommand{\bC}{$C$}
\newcommand{\dl}{{\delta}}
\newcommand{\g}{{\cal G}}
\newcommand{\bg}{{\bar \g}}
\newcommand{\bbg}{{\bar \bg}}
\newcommand{\uhg}{U_h(\g)}
\newcommand{\uog}{U_0(\g)}

\newcommand{\om}{{\omega}}
\newcommand{\Tomi}{T_{\omega_i}}
\newcommand{\Tomj}{T_{\omega_j}}
\newcommand{\hTomi}{{{\hat T}_{\omega_i}}}
\newcommand{\hTomj}{{{\hat T}_{\omega_j}}}
\newcommand{\hTomk}{{{\hat T}_{\omega_k}}}

\newcommand{\bps}{{\bar \psi}}
\newcommand{\hbps}{{\hat {\bar \psi}}}

\newcommand{\Epipw}
{{({\cal E}^w, \Pi^w, p^w)}}
\newcommand{\Epipt}
{{({\cal E}^\tau, \Pi^\tau, p^\tau)}}

\newcommand{\Epipdash}{{({\cal E'}, \Pi', p')}}
\newcommand{\Epipdagger}{{({\cal E^\dagger}, \Pi^\dagger, p^\dagger)}}
\newcommand{\Epidagger}{{({\cal E^\dagger}, \Pi^\dagger)}}
\newcommand{\Epipsi}
{{({\cal E}^{\sigma(i)}, \Pi^{\sigma(i)}, p^{\sigma(i)})}}

\newcommand{\Epipxy}
{{({\cal E}, \Pi_{(x,y)}, p)}}
\newcommand{\Epipnisan}
{{({\cal E}, \Pi_{(2,3)}, p)}}
\newcommand{\EpipNniNsan}
{{({\cal E}, \Pi_{(N-3,N-2)}, p)}}

\newcommand{\cwhl}{{[}}
\newcommand{\cwhr}{{]}}

\newcommand{\whl}{{[\![}}
\newcommand{\whr}{{]\!]}}

\newcommand{\Da}{\Delta^{(\al)}}
\newcommand{\hRa}{{\hat R}_{\al}}
\newcommand{\Epipc}{{({\cal E}, \Pi', p)}}

\newcommand{\uhgd}{U_h({\dot {\cal G}})}

\newcommand{\uhgs}{U_h({\cal G})^\sigma}
\newcommand{\uhgas}{U_h({\cal G}^{(\al)})^\sigma}

\newcommand{\uhgsa}{U_h({\cal G}^{\sigma_\al})}

\newcommand{\uhhs}{U_h({\cal H})^\sigma}

\newcommand{\tufhbs}{{\tilde U}^{\flat}_h({\cal B}^+)^\sigma}
\newcommand{\ufhbs}{U^{\flat}_h({\cal B}^+)^\sigma}

\newcommand{\chkl}{C[[h]]\left[K_\lambda \right]}
\newcommand{\chs}{C[[h]]\langle\sigma\rangle}
\newcommand{\chh}{C[[h]]\left[{\cal H}\right]}
\newcommand{\fchh}{C((h))\left[{\cal H}\right]}
\newcommand{\fchkl}{C((h))\left[K_\lambda \right]}

\newcommand{\dotwhite}{
\begin{picture}(10,10)(0,3)
\put(05, 05){\circle{5}}
\end{picture}}

\newcommand{\dotgray}{
\begin{picture}(10,10)(0,3)
\put(05, 05){\circle{5}}
\put(03, 03){\line(1,1){4}}
\put(03, 07){\line(1,-1){4}}
\end{picture}}

\newcommand{\dotblack}{
\begin{picture}(10,10)(0,3)
\put(05, 05){\circle*{5}}
\end{picture}}

\newcommand{\dotcross}{
\begin{picture}(10,10)(0,3)
\put(03, 03){\line(1,1){4}}
\put(03, 07){\line(1,-1){4}}
\end{picture}}

\newcommand{\dotwhiteblack}{
\begin{picture}(10,10)(0,3)
\put(05, 05){\circle{5}}\put(05, 05){\circle*{2}}
\end{picture}}

\newcommand{\semiline}{
\begin{picture}(10,18)(0,6)
\put(3, 13){\line(1,1){4}} 
\put(3, 17){\line(1,-1){4}}
\put(3, 0){\line(1,1){4}} 
\put(3, 4){\line(1,-1){4}}
\put(4, 7){\vdots}
\put(6, 7){\vdots}
\end{picture}}

\newcommand{\noline}{
\begin{picture}(10,18)(0,6)
\put(05, 15){\circle{5}}
\put(05, 02){\circle{5}}
\end{picture}}

\newcommand{\twoline}{
\begin{picture}(10,18)(0,6)
\put(05, 15){\circle{5}}
\put(3, 13){\line(1,1){4}} 
\put(3, 17){\line(1,-1){4}}
\put(05, 02){\circle{5}}
\put(3, 0){\line(1,1){4}} 
\put(3, 4){\line(1,-1){4}}
\put(4, 4){\line(0,1){9}}
\put(6, 4){\line(0,1){9}}
\end{picture}}

\newcommand{\DiagramA}{
\setlength{\unitlength}{1mm}
\begin{picture}(180,35)(0,0)

\put(5,25){$A$}
\put(33, 23){\line(1,1){4}} \put(34, 30){$1$} \put(29,15){$\be_1-\be_2$} 
\put(33, 27){\line(1,-1){4}}
\put(37, 25){\line(1,0){16}}

\put(53, 23){\line(1,1){4}} \put(54, 30){$2$} \put(49,15){$\be_2-\be_3$}
\put(53, 27){\line(1,-1){4}}
\put(58, 25){\line(1,0){10}}

\put(70, 25){\dots}
\put(76, 25){\line(1,0){27}}

\put(103, 23){\line(1,1){4}} \put(100, 30){$N-1$}
\put(96,15){$\be_{N-1}-\be_{N}$}
\put(103, 27){\line(1,-1){4}}

\end{picture}}

\newcommand{\DiagramB}{
\setlength{\unitlength}{1mm}
\begin{picture}(120,35)(0,0)

\put(5,25){$B$}

\put(33, 23){\line(1,1){4}} \put(34, 30){$1$} \put(29,15){$\be_1-\be_2$}
\put(33, 27){\line(1,-1){4}}
\put(38, 25){\line(1,0){10}}

\put(50, 25){\dots}
\put(56, 25){\line(1,0){7}}

\put(63, 23){\line(1,1){4}} \put(59, 30){$N-1$}
\put(54,15){$\be_{N-2}-\be_{N-1}$}
\put(63, 27){\line(1,-1){4}}

\put(85, 25){\circle{5}}\put(85, 25){\circle*{2}}
 \put(84, 30){$N$} \put(83,15){$\be_N$} 
\put(82, 25){\line(-1,-1){4}}
\put(82, 25){\line(-1,1){4}}
\put(80, 24){\line(-1,0){13}}
\put(80, 26){\line(-1,0){13}}

\put(90, 24){.}
\end{picture}}

\newcommand{\DiagramC}{
\setlength{\unitlength}{1mm}
\begin{picture}(110,35)(0,0)

\put(5,25){$C$}

\put(33, 23){\line(1,1){4}} \put(34, 30){$1$} \put(29,15){$\be_1-\be_2$}
\put(33, 27){\line(1,-1){4}}
\put(38, 25){\line(1,0){10}}

\put(50, 25){\dots}
\put(56, 25){\line(1,0){7}}

\put(63, 23){\line(1,1){4}} \put(59, 30){$N-1$}
\put(54,15){$\be_{N-2}-\be_{N-1}$}
\put(63, 27){\line(1,-1){4}}
\put(68, 25){\line(1,1){4}}
\put(68, 25){\line(1,-1){4}}
\put(70, 24){\line(1,0){13}}
\put(70, 26){\line(1,0){13}}

\put(85, 25){\circle{5}}
 \put(84, 30){$N$} \put(82,15){$2\be_N$} 
\put(90, 24){.}

\end{picture}}

\newcommand{\DiagramD}{
\setlength{\unitlength}{1mm}
\begin{picture}(120,55)(0,0)

\put(5,25){$D$}

\put(33, 23){\line(1,1){4}} \put(34, 30){$1$} \put(29,15){$\be_1-\be_2$}
\put(33, 27){\line(1,-1){4}}
\put(38, 25){\line(1,0){10}}

\put(50, 25){\dots}
\put(56, 25){\line(1,0){7}}

\put(63, 23){\line(1,1){4}} \put(59, 30){$N-2$}
\put(54,15){$\be_{N-2}-\be_{N-1}$}
\put(63, 27){\line(1,-1){4}}

\put(68, 25){\line(1,1){14}} 
\put(68, 25){\line(1,-1){14}}
\put(83, 37){\line(1,1){4}} \put(77,48){$\be_{N-1}-\be_{N}$} \put(90,
38){$N-1$}
\put(83, 41){\line(1,-1){4}}
\put(83, 10){\line(1,1){4}} \put(77,2){$\be_{N-1}+\be_{N}$} \put(90,
11){$N$}
\put(83, 14){\line(1,-1){4}}
\put(84, 30){\vdots}
\put(86, 30){\vdots}
\put(84, 24){\vdots}
\put(86, 24){\vdots}
\put(84, 18){\vdots}
\put(86, 18){\vdots}

\put(100,24){.}

\end{picture}}

\newcommand{\DiagramAAa}{
\setlength{\unitlength}{1mm}
\begin{picture}(30,15)(0,0)

\put(5, 5){\circle{5}}
 \put(4, 10){$0$} 
\put(8, 4){\line(1,0){15}}
\put(8, 6){\line(1,0){15}}

\put(25, 5){\circle{5}} \put(24, 10){$1$}

\end{picture}}

\newcommand{\DiagramAA}{
\setlength{\unitlength}{1mm}
\begin{picture}(180,55)(0,0)

\put(5,25){$(AA)$}
\put(33, 23){\line(1,1){4}} \put(34, 30){$1$} \put(29,15){$\be_1-\be_2$} 
\put(33, 27){\line(1,-1){4}}
\put(37, 25){\line(1,0){16}}
\put(37, 25){\line(3,2){31}}

\put(53, 23){\line(1,1){4}} \put(54, 30){$2$} \put(49,15){$\be_2-\be_3$}
\put(53, 27){\line(1,-1){4}}
\put(58, 25){\line(1,0){10}}

\put(70, 25){\dots}
\put(76, 25){\line(1,0){27}}
\put(103, 25){\line(-3,2){31}}

\put(103, 23){\line(1,1){4}} \put(104, 30){$N-1$}
\put(96,15){$\be_{N-1}-\be_{N}$}
\put(103, 27){\line(1,-1){4}}

\put(68, 43){\line(1,1){4}} \put(69, 50){$0$}
\put(78,44){$\dl-\be_1+\be_N$} 
\put(68, 47){\line(1,-1){4}}

\end{picture}}

\newcommand{\DiagramBBb}{
\setlength{\unitlength}{1mm}
\begin{picture}(30,15)(0,0)

\put(5, 5){\circle{5}}\put(5, 5){\circle*{2}}
 \put(4, 10){$0$} 
\put(8, 4){\line(1,0){15}}
\put(8, 6){\line(1,0){15}}

\put(25, 5){\circle{5}}\put(25, 5){\circle*{2}}
 \put(24, 10){$1$}

\end{picture}}

\newcommand{\DiagramBB}{
\setlength{\unitlength}{1mm}
\begin{picture}(180,35)(0,0)

\put(5,25){$(BB)$}
\put(35, 25){\circle{5}}\put(35, 25){\circle*{2}}
 \put(34, 30){0} \put(30,15){$\dl-\be_1$} 
\put(38, 25){\line(1,1){4}}
\put(38, 25){\line(1,-1){4}}
\put(40, 24){\line(1,0){13}}
\put(40, 26){\line(1,0){13}}

\put(53, 23){\line(1,1){4}} \put(54, 30){$1$} \put(49,15){$\be_1-\be_2$}
\put(53, 27){\line(1,-1){4}}
\put(58, 25){\line(1,0){10}}

\put(70, 25){\dots}
\put(76, 25){\line(1,0){7}}

\put(83, 23){\line(1,1){4}} \put(79, 30){$N-1$}
\put(74,15){$\be_{N-2}-\be_{N-1}$}
\put(83, 27){\line(1,-1){4}}

\put(105, 25){\circle{5}}\put(105, 25){\circle*{2}}
 \put(104, 30){$N$} \put(103,15){$\be_N$} 
\put(102, 25){\line(-1,-1){4}}
\put(102, 25){\line(-1,1){4}}
\put(100, 24){\line(-1,0){13}}
\put(100, 26){\line(-1,0){13}}

\end{picture}}

\newcommand{\DiagramCBb}{
\setlength{\unitlength}{1mm}
\begin{picture}(25,15)(0,0)

\put(5, 5){\circle{5}} \put(4,10){$0$}
\put(25, 5){\circle{5}}\put(25, 5){\circle*{2}}
 \put(24, 10){$1$}  
\put(22, 5){\line(-1,-1){4}}
\put(22, 5){\line(-1,1){4}}
\put(20, 3.5){\line(-1,0){13}}
\put(20, 4.5){\line(-1,0){13}}
\put(20, 5.5){\line(-1,0){13}}
\put(20, 6.5){\line(-1,0){13}}

\end{picture}}

\newcommand{\DiagramCB}{
\setlength{\unitlength}{1mm}
\begin{picture}(180,35)(0,0)

\put(5,25){$(CB)$}
\put(35, 25){\circle{5}}
 \put(34, 30){0} \put(29,15){$\dl-2\be_1$} 
\put(37, 24){\line(1,0){13}}
\put(37, 26){\line(1,0){13}}
\put(52, 25){\line(-1,-1){4}}
\put(52, 25){\line(-1,1){4}}

\put(53, 23){\line(1,1){4}} \put(54, 30){$1$} \put(49,15){$\be_1-\be_2$}
\put(53, 27){\line(1,-1){4}}
\put(58, 25){\line(1,0){10}}

\put(70, 25){\dots}
\put(76, 25){\line(1,0){7}}

\put(83, 23){\line(1,1){4}} \put(79, 30){$N-1$}
\put(74,15){$\be_{N-2}-\be_{N-1}$}
\put(83, 27){\line(1,-1){4}}

\put(105, 25){\circle{5}}\put(105, 25){\circle*{2}}
 \put(104, 30){$N$} \put(103,15){$\be_N$} 
\put(102, 25){\line(-1,-1){4}}
\put(102, 25){\line(-1,1){4}}
\put(100, 24){\line(-1,0){13}}
\put(100, 26){\line(-1,0){13}}

\end{picture}}

\newcommand{\DiagramDBb}{
\setlength{\unitlength}{1mm}
\begin{picture}(40,40)(0,0)

\put(13, 39){\line(1,1){4}} \put(8,40){$0$}
\put(13, 43){\line(1,-1){4}}
\put(13, 8){\line(1,1){4}} \put(8, 9){$1$}
\put(13, 12){\line(1,-1){4}}
\put(14, 30){\vdots}
\put(16, 30){\vdots}
\put(14, 24){\vdots}
\put(16, 24){\vdots}
\put(14, 18){\vdots}
\put(16, 18){\vdots}
\put(29, 28){\line(-1,1){11}} 
\put(31, 29){\line(-1,1){13}} \put(32, 27){\line(0,1){4}}
                              \put(32, 27){\line(-1,0){4}}

\put(29, 22){\line(-1,-1){11}}
\put(31, 21){\line(-1,-1){13}} \put(32, 23){\line(0,-1){4}}
                               \put(32, 23){\line(-1,0){4}}
\put(35, 25){\circle{5}} \put(34, 30){$2$} 
\put(35, 25){\circle*{2}}

\end{picture}}

\newcommand{\DiagramDB}{
\setlength{\unitlength}{1mm}
\begin{picture}(180,35)(0,0)

\put(5,25){$(DB)$}
\put(33, 37){\line(1,1){4}} \put(27,48){$\dl-\be_{1}-\be_{2}$} \put(28,
38){$0$}
\put(33, 41){\line(1,-1){4}}
\put(33, 10){\line(1,1){4}} \put(29,2){$\be_{1}-\be_{2}$} \put(28, 11){$1$}
\put(33, 14){\line(1,-1){4}}
\put(34, 30){\vdots}
\put(36, 30){\vdots}
\put(34, 24){\vdots}
\put(36, 24){\vdots}
\put(34, 18){\vdots}
\put(36, 18){\vdots}
\put(52, 25){\line(-1,1){14}} 
\put(52, 25){\line(-1,-1){14}}

\put(53, 23){\line(1,1){4}} \put(54, 30){$2$} \put(49,15){$\be_2-\be_3$}
\put(53, 27){\line(1,-1){4}}
\put(58, 25){\line(1,0){10}}

\put(70, 25){\dots}
\put(76, 25){\line(1,0){7}}

\put(83, 23){\line(1,1){4}} \put(79, 30){$N-1$}
\put(74,15){$\be_{N-2}-\be_{N-1}$}
\put(83, 27){\line(1,-1){4}}

\put(105, 25){\circle{5}}\put(105, 25){\circle*{2}}
 \put(104, 30){$N$} \put(103,15){$\be_N$} 
\put(102, 25){\line(-1,-1){4}}
\put(102, 25){\line(-1,1){4}}
\put(100, 24){\line(-1,0){13}}
\put(100, 26){\line(-1,0){13}}

\end{picture}}

\newcommand{\DiagramDC}{
\setlength{\unitlength}{1mm}
\begin{picture}(180,35)(0,0)

\put(5,25){$(DC)$}
\put(33, 37){\line(1,1){4}} \put(27,48){$\dl-\be_{1}-\be_{2}$} \put(28,
38){$0$}
\put(33, 41){\line(1,-1){4}}
\put(33, 10){\line(1,1){4}} \put(29,2){$\be_{1}-\be_{2}$} \put(28, 11){$1$}
\put(33, 14){\line(1,-1){4}}
\put(34, 30){\vdots}
\put(36, 30){\vdots}
\put(34, 24){\vdots}
\put(36, 24){\vdots}
\put(34, 18){\vdots}
\put(36, 18){\vdots}
\put(52, 25){\line(-1,1){14}} 
\put(52, 25){\line(-1,-1){14}}

\put(53, 23){\line(1,1){4}} \put(54, 30){$2$} \put(49,15){$\be_2-\be_3$}
\put(53, 27){\line(1,-1){4}}
\put(58, 25){\line(1,0){10}}

\put(70, 25){\dots}
\put(76, 25){\line(1,0){7}}

\put(83, 23){\line(1,1){4}} \put(79, 30){$N-1$}
\put(74,15){$\be_{N-2}-\be_{N-1}$}
\put(83, 27){\line(1,-1){4}}
\put(88, 25){\line(1,1){4}}
\put(88, 25){\line(1,-1){4}}
\put(90, 24){\line(1,0){13}}
\put(90, 26){\line(1,0){13}}

\put(105, 25){\circle{5}}
 \put(104, 30){$N$} \put(102,15){$2\be_N$}

\end{picture}}

\newcommand{\DiagramCC}{
\setlength{\unitlength}{1mm}
\begin{picture}(180,30)(0,0)

\put(5,25){$(CC)$}
\put(35, 25){\circle{5}}
 \put(34, 30){0} \put(29,15){$\dl-2\be_1$} 
\put(37, 24){\line(1,0){13}}
\put(37, 26){\line(1,0){13}}
\put(52, 25){\line(-1,-1){4}}
\put(52, 25){\line(-1,1){4}}

\put(53, 23){\line(1,1){4}} \put(54, 30){$1$} \put(49,15){$\be_1-\be_2$}
\put(53, 27){\line(1,-1){4}}
\put(58, 25){\line(1,0){10}}

\put(70, 25){\dots}
\put(76, 25){\line(1,0){7}}

\put(83, 23){\line(1,1){4}} \put(79, 30){$N-1$}
\put(74,15){$\be_{N-2}-\be_{N-1}$}
\put(83, 27){\line(1,-1){4}}
\put(88, 25){\line(1,1){4}}
\put(88, 25){\line(1,-1){4}}
\put(90, 24){\line(1,0){13}}
\put(90, 26){\line(1,0){13}}

\put(105, 25){\circle{5}}
 \put(104, 30){$N$} \put(102,15){$2\be_N$}

\end{picture}}

\newcommand{\DiagramCD}{
\setlength{\unitlength}{1mm}
\begin{picture}(180,55)(0,0)

\put(5,25){$(CD)$}
\put(35, 25){\circle{5}}
 \put(34, 30){0} \put(29,15){$\dl-2\be_1$} 
\put(37, 24){\line(1,0){13}}
\put(37, 26){\line(1,0){13}}
\put(52, 25){\line(-1,-1){4}}
\put(52, 25){\line(-1,1){4}}

\put(53, 23){\line(1,1){4}} \put(54, 30){$1$} \put(49,15){$\be_1-\be_2$}
\put(53, 27){\line(1,-1){4}}
\put(58, 25){\line(1,0){10}}

\put(70, 25){\dots}
\put(76, 25){\line(1,0){7}}

\put(83, 23){\line(1,1){4}} \put(79, 30){$N-2$}
\put(74,15){$\be_{N-2}-\be_{N-1}$}
\put(83, 27){\line(1,-1){4}}

\put(88, 25){\line(1,1){14}} 
\put(88, 25){\line(1,-1){14}}
\put(103, 37){\line(1,1){4}} \put(97,48){$\be_{N-1}-\be_{N}$} \put(110,
38){$N-1$}
\put(103, 41){\line(1,-1){4}}
\put(103, 10){\line(1,1){4}} \put(97,2){$\be_{N-1}+\be_{N}$} \put(110,
11){$N$}
\put(103, 14){\line(1,-1){4}}
\put(104, 30){\vdots}
\put(106, 30){\vdots}
\put(104, 24){\vdots}
\put(106, 24){\vdots}
\put(104, 18){\vdots}
\put(106, 18){\vdots}

\end{picture}}

\newcommand{\DiagramDDd}{
\setlength{\unitlength}{1mm}
\begin{picture}(35,30)(0,0)

\put(15, 20){\circle{5}}
\put(13, 18){\line(1,1){4}} \put(8,19){$0$}
\put(13, 22){\line(1,-1){4}}
\put(17, 18){\line(1,-1){11}}

\put(14, 18){\line(0,-1){11}}
\put(16, 18){\line(0,-1){11}}
\put(17, 20){\line(1,0){11}}

\put(30, 20){\circle{5}}
\put(28, 18){\line(1,1){4}} \put(35, 19){$2$} 
\put(28, 22){\line(1,-1){4}}

\put(29, 18){\line(0,-1){11}}
\put(31, 18){\line(0,-1){11}}
\put(17, 5){\line(1,0){11}}

\put(15, 5){\circle{5}}
\put(13, 3){\line(1,1){4}} \put(8,4){$1$}
\put(13, 7){\line(1,-1){4}}
\put(17, 7){\line(1,1){11}}

\put(30, 5){\circle{5}}
\put(28, 3){\line(1,1){4}} \put(35, 4){$3$} 
\put(28, 7){\line(1,-1){4}}

\end{picture}}

\newcommand{\DiagramDD}{
\setlength{\unitlength}{1mm}
\begin{picture}(180,55)(0,0)

\put(5,25){$(DD)$}
\put(33, 37){\line(1,1){4}} \put(27,48){$\dl-\be_{1}-\be_{2}$} \put(28,
38){$0$}
\put(33, 41){\line(1,-1){4}}
\put(33, 10){\line(1,1){4}} \put(29,2){$\be_{1}-\be_{2}$} \put(28, 11){$1$}
\put(33, 14){\line(1,-1){4}}
\put(34, 30){\vdots}
\put(36, 30){\vdots}
\put(34, 24){\vdots}
\put(36, 24){\vdots}
\put(34, 18){\vdots}
\put(36, 18){\vdots}
\put(52, 25){\line(-1,1){14}} 
\put(52, 25){\line(-1,-1){14}}

\put(53, 23){\line(1,1){4}} \put(54, 30){$2$} \put(49,15){$\be_2-\be_3$}
\put(53, 27){\line(1,-1){4}}

\put(58, 25){\line(1,0){10}}

\put(70, 25){\dots}
\put(76, 25){\line(1,0){7}}

\put(83, 23){\line(1,1){4}} \put(79, 30){$N-2$}
\put(74,15){$\be_{N-2}-\be_{N-1}$}
\put(83, 27){\line(1,-1){4}}

\put(88, 25){\line(1,1){14}} 
\put(88, 25){\line(1,-1){14}}
\put(103, 37){\line(1,1){4}} \put(97,48){$\be_{N-1}-\be_{N}$} \put(110,
38){$N-1$}
\put(103, 41){\line(1,-1){4}}
\put(103, 10){\line(1,1){4}} \put(97,2){$\be_{N-1}+\be_{N}$} \put(110,
11){$N$}
\put(103, 14){\line(1,-1){4}}
\put(104, 30){\vdots}
\put(106, 30){\vdots}
\put(104, 24){\vdots}
\put(106, 24){\vdots}
\put(104, 18){\vdots}
\put(106, 18){\vdots}

\end{picture}}

\newcommand{\hoshi}{
\setlength{\unitlength}{1mm}
\begin{picture}(15,15)(0,14)
  
\put(0.5,14){$\times$}
\put(1,18){$i$}

\put(3,16){\line(1,1){5}}
\put(7.5,21){$\times$}\put(11,21.5){$j$}

\put(3,14){\line(1,-1){5}}
\put(7.5, 7){$\times$}\put(11,8.5){$k$}

\put(8, 14){$*$}

\end{picture}}

\newcommand{\DiagramCDd}{
\setlength{\unitlength}{1mm}
\begin{picture}(50,50)(0,0)

\put(5, 25){\circle{5}}
 \put(4, 30){0} 
\put(7, 24){\line(1,0){13}}
\put(7, 26){\line(1,0){13}}
\put(22, 25){\line(-1,-1){4}}
\put(22, 25){\line(-1,1){4}}

\put(23, 23){\line(1,1){4}} \put(24, 30){$1$} 
\put(23, 27){\line(1,-1){4}}
\put(25, 25){\circle{5}}
\put(27, 27){\line(1,1){11}}
\put(40, 40){\circle{5}} \put(45,39){$2$}
\put(27, 23){\line(1,-1){11}}
\put(40, 10){\circle{5}} \put(45,9){$3$}

\end{picture}}

\newcommand{\DiagramDDdd}{
\setlength{\unitlength}{1mm}
\begin{picture}(35,30)(0,0)

\put(15, 25){\circle{5}}
\put(13, 23){\line(1,1){4}} \put(7,24){$\al_0$}
\put(13, 27){\line(1,-1){4}}

\put(13.5, 23){\line(0,-1){16}}
\put(16.5, 23){\line(0,-1){16}}

\put(15, 5){\circle{5}}
\put(13, 3){\line(1,1){4}} \put(7,4){$\al_1$}
\put(13, 7){\line(1,-1){4}}

\put(25,15){\circle{5}} \put(30,14){$\al_2+\al_3$}
\put(24,18){\line(-1,1){6}}
\put(22,16){\line(-1,1){6}}
\put(22,14){\line(-1,-1){6}}
\put(24,12){\line(-1,-1){6}}

\end{picture}}

\newcommand{\twistdo}{
\setlength{\unitlength}{1mm}
\begin{picture}(30,45)(0,0)

\put(3, 23){\line(1,1){4}} \put(4, 30){$1$} 
\put(3, 27){\line(1,-1){4}}
\put(5, 25){\circle{5}}
\put(7, 27){\line(1,1){11}}
\put(9, 34){{\footnotesize {$-1$}}}
\put(20, 40){\circle{5}} \put(25,39){$2$}
\put(7, 23){\line(1,-1){11}}
\put(9, 14){{\footnotesize {$-x$}}}

\put(20, 10){\circle{5}} \put(25,9){$3$}

\end{picture}}

\newcommand{\twistda}{
\setlength{\unitlength}{1mm}
\begin{picture}(50,45)(0,0)

\put(5, 25){\circle{5}}
 \put(4, 30){0} 
\put(7, 25){\line(1,0){16}}
\put(11, 28){{\footnotesize {$x+1$}}}

\put(23, 23){\line(1,1){4}} \put(24, 30){$1$} 
\put(23, 27){\line(1,-1){4}}
\put(25, 25){\circle{5}}
\put(27, 27){\line(1,1){11}}
\put(29, 34){{\footnotesize {$-1$}}}
\put(40, 40){\circle{5}} \put(45,39){$2$}
\put(27, 23){\line(1,-1){11}}
\put(29, 14){{\footnotesize {$-x$}}}

\put(40, 10){\circle{5}} \put(45,9){$3$}

\end{picture}}

\newcommand{\twistdaa}{
\setlength{\unitlength}{1mm}
\begin{picture}(50,45)(0,0)

\put(5, 25){\circle{5}}
 \put(4, 30){0} 
\put(7, 25){\line(1,0){16}}
\put(11, 28){{\footnotesize {$x+1$}}}

\put(23, 23){\line(1,1){4}} \put(24, 30){$1$} 
\put(23, 27){\line(1,-1){4}}
\put(25, 25){\circle{5}}
\put(27, 27){\line(1,1){11}}
\put(29, 34){{\footnotesize {$-x$}}}
\put(40, 40){\circle{5}} \put(45,39){$2$}
\put(27, 23){\line(1,-1){11}}
\put(29, 14){{\footnotesize {$-1$}}}

\put(40, 10){\circle{5}} \put(45,9){$3$}

\end{picture}}

\newcommand{\twistdb}{
\setlength{\unitlength}{1mm}
\begin{picture}(35,30)(0,0)

\put(15, 20){\circle{5}}
\put(13, 18){\line(1,1){4}} \put(14, 25){$0$}
\put(13, 22){\line(1,-1){4}}
\put(17, 18){\line(1,-1){11}}

\put(15, 18){\line(0,-1){11}}
\put( 5, 12){{\tiny{$-(x+1)$}}}

\put(17, 20){\line(1,0){11}}
\put(22, 21){{\footnotesize {$x$}}}
\put(19, 16){{\scriptsize {$1$}}}
\put(24, 16){{\scriptsize {$1$}}}

\put(30, 20){\circle{5}}
\put(28, 18){\line(1,1){4}} \put(29, 25){$2$} 
\put(28, 22){\line(1,-1){4}}

\put(30, 18){\line(0,-1){11}}
\put(31, 12){{\tiny{$-(x+1)$}}}
\put(17, 5){\line(1,0){11}}
\put(22, 2){{\footnotesize {$x$}}}

\put(15, 5){\circle{5}}
\put(13, 3){\line(1,1){4}} \put(14,-3){$1$}
\put(13, 7){\line(1,-1){4}}
\put(17, 7){\line(1,1){11}}

\put(30, 5){\circle{5}}
\put(28, 3){\line(1,1){4}} \put(29,-3){$3$} 
\put(28, 7){\line(1,-1){4}}

\end{picture}}

\newcommand{\twistdbb}{
\setlength{\unitlength}{1mm}
\begin{picture}(35,30)(0,0)

\put(15, 20){\circle{5}}
\put(13, 18){\line(1,1){4}} \put(14, 25){$0$}
\put(13, 22){\line(1,-1){4}}
\put(17, 18){\line(1,-1){11}}

\put(15, 18){\line(0,-1){11}}
\put( 5, 12){{\tiny{$-(x+1)$}}}

\put(17, 20){\line(1,0){11}}
\put(22, 21){{\footnotesize {$1$}}}
\put(19, 16){{\scriptsize {$x$}}}
\put(24, 16){{\scriptsize {$x$}}}

\put(30, 20){\circle{5}}
\put(28, 18){\line(1,1){4}} \put(29, 25){$2$} 
\put(28, 22){\line(1,-1){4}}

\put(30, 18){\line(0,-1){11}}
\put(31, 12){{\tiny{$-(x+1)$}}}
\put(17, 5){\line(1,0){11}}
\put(22, 2){{\footnotesize {$1$}}}

\put(15, 5){\circle{5}}
\put(13, 3){\line(1,1){4}} \put(14,-3){$1$}
\put(13, 7){\line(1,-1){4}}
\put(17, 7){\line(1,1){11}}

\put(30, 5){\circle{5}}
\put(28, 3){\line(1,1){4}} \put(29,-3){$3$} 
\put(28, 7){\line(1,-1){4}}

\end{picture}}

\newcommand{\twistdc}{
\setlength{\unitlength}{1mm}
\begin{picture}(50,45)(0,0)

\put(23, 23){\line(1,1){4}} \put(24, 30){$2$} 
\put(23, 27){\line(1,-1){4}}
\put(25, 25){\circle{5}}
\put(23, 27){\line(-1,1){11}}
\put(16, 34){{\footnotesize {$-1$}}}
\put(10, 40){\circle{5}} \put(3,39){$0$}
\put(23, 23){\line(-1,-1){11}}
\put(16, 14){{\footnotesize {$-x$}}}

\put(10, 10){\circle{5}} \put(3,9){$1$}

\put(45, 25){\circle{5}}
 \put(44, 30){$3$} 
\put(27, 25){\line(1,0){16}}
\put(31, 28){{\footnotesize {$x+1$}}}

\end{picture}}

\newcommand{\twistdcc}{
\setlength{\unitlength}{1mm}
\begin{picture}(50,45)(0,0)

\put(23, 23){\line(1,1){4}} \put(24, 30){$2$} 
\put(23, 27){\line(1,-1){4}}
\put(25, 25){\circle{5}}
\put(23, 27){\line(-1,1){11}}
\put(16, 34){{\footnotesize {$-x$}}}
\put(10, 40){\circle{5}} \put(3,39){$0$}
\put(23, 23){\line(-1,-1){11}}
\put(16, 14){{\footnotesize {$-1$}}}

\put(10, 10){\circle{5}} \put(3,9){$1$}

\put(45, 25){\circle{5}}
 \put(44, 30){$3$} 
\put(27, 25){\line(1,0){16}}
\put(31, 28){{\footnotesize {$x+1$}}}

\end{picture}}

\newcommand{\tateyaa}{
\setlength{\unitlength}{1mm}
\begin{picture}(5,18)(0,0)

\put(0, 16){\vector(0,-1){16}}
\put(0, 0){\vector(0,1){16}}
\put(2, 7){{\scriptsize  1}}
\end{picture}}

\newcommand{\tateyab}{
\setlength{\unitlength}{1mm}
\begin{picture}(5,18)(0,0)

\put(0, 16){\vector(0,-1){16}}
\put(0, 0){\vector(0,1){16}}
\put(2, 7){{\scriptsize  2}}
\end{picture}}

\newcommand{\yokoyaa}{
\setlength{\unitlength}{1mm}
\begin{picture}(5,18)(0,0)

\put(16, 0){\vector(-1,0){16}}
\put(0, 0){\vector(1,0){16}}
\put(7, 2){{\scriptsize  0}}
\end{picture}}

\newcommand{\yokoyab}{
\setlength{\unitlength}{1mm}
\begin{picture}(5,18)(0,0)

\put(16, 0){\vector(-1,0){16}}
\put(0, 0){\vector(1,0){16}}
\put(7, 2){{\scriptsize  3}}
\end{picture}}

\newcommand{\comtwistd}{
\setlength{\unitlength}{1mm}
\begin{picture}(140,165)(0,0)

\put(0,110){\twistda}\put(80,110){\twistdaa}
\put(6,70){\twistdb}\put(86,70){\twistdbb}
\put(0,0){\twistdc}\put(80,0){\twistdcc}

\put(30,100){\tateyaa}\put(110,100){\tateyaa}
\put(30,45){\tateyab}\put(110,45){\tateyab}
\put(56,25){\yokoyaa}
\put(56,135){\yokoyab}

\put(60,77){\yokoyaa}
\put(60,87){\yokoyab}

\put(0,165){.}
\end{picture}}

\newcommand{\sfoo}{
\setlength{\unitlength}{1mm}
\begin{picture}(35,15)(0,0)

\put(1.5, 8){{\scriptsize  1}}
\put(11.5, 8){{\scriptsize  2}}
\put(21.5, 8){{\scriptsize  3}}
\put(31.5, 8){{\scriptsize  4}}

\put(2, 5){\circle{3}}
\put(3.5, 5){\line(1,0){7}}

\put(12,5){\circle{3}}
\put(13.5,4){$\Longrightarrow$}

\put(22,5){\circle{3}}
\put(23.5, 5){\line(1,0){7}}

\put(32,5){\circle{3}}
\put(30.5,4){$\times$}

\end{picture}}

\newcommand{\foa}{
\setlength{\unitlength}{1mm}
\begin{picture}(45,15)(0,0)

\put(1.5, 8){{\scriptsize  0}}
\put(11.5, 8){{\scriptsize  1}}
\put(21.5, 8){{\scriptsize  2}}
\put(31.5, 8){{\scriptsize  3}}
\put(41.5, 8){{\scriptsize  4}}

\put(2, 5){\circle{3}}
\put(3.5, 5){\line(1,0){7}}

\put(12,5){\circle{3}}
\put(13.5, 5){\line(1,0){7}}

\put(22,5){\circle{3}}
\put(23.5,4){$\Longrightarrow$}

\put(32,5){\circle{3}}
\put(33.5,5){\line(1,0){7}}

\put(42,5){\circle{3}}

\end{picture}}

\newcommand{\faa}{
\setlength{\unitlength}{1mm}
\begin{picture}(35,15)(0,0)

\put(1.5, 8){{\scriptsize  1}}
\put(11.5, 8){{\scriptsize  2}}
\put(21.5, 8){{\scriptsize  3}}
\put(31.5, 8){{\scriptsize  4}}
\put(36.5, 15){{\scriptsize  0+1+2+3}}

\put(2, 5){\circle{3}}
\put(3.5, 5){\line(1,0){7}}

\put(12,5){\circle{3}}
\put(13.5,4){$\Longrightarrow$}

\put(22,5){\circle{3}}
\put(23.5, 5){\line(1,0){7}}

\put(32,5){\circle{3}}
\put(30.5,4){$\times$}

\put(32.5,6.5){\line(1,1){5}}
\put(33,6){\line(1,1){5}}
\put(33.5,5.5){\line(1,1){5}}

\put(39,12){\circle{3}}

\end{picture}}

\newcommand{\fab}{
\setlength{\unitlength}{1mm}
\begin{picture}(35,15)(0,0)

\put(-0.5, 8){{\scriptsize  0+1}}
\put(11.5, 8){{\scriptsize  2}}
\put(21.5, 8){{\scriptsize  3}}
\put(31.5, 8){{\scriptsize  4}}
\put(36.5, 15){{\scriptsize  1+2+3}}

\put(2, 5){\circle{3}}
\put(3.5, 5){\line(1,0){7}}

\put(12,5){\circle{3}}
\put(13.5,4){$\Longrightarrow$}

\put(22,5){\circle{3}}
\put(23.5, 5){\line(1,0){7}}

\put(32,5){\circle{3}}
\put(30.5,4){$\times$}

\put(32.5,6.5){\line(1,1){5}}
\put(33,6){\line(1,1){5}}
\put(33.5,5.5){\line(1,1){5}}

\put(39,12){\circle{3}}

\end{picture}}

\newcommand{\fac}{
\setlength{\unitlength}{1mm}
\begin{picture}(35,15)(0,0)

\put(1.5, 8){{\scriptsize  0}}
\put(9.5, 8){{\scriptsize  1+2}}
\put(21.5, 8){{\scriptsize  3}}
\put(31.5, 8){{\scriptsize  4}}
\put(36.5, 15){{\scriptsize  2+3}}

\put(2, 5){\circle{3}}
\put(3.5, 5){\line(1,0){7}}

\put(12,5){\circle{3}}
\put(13.5,4){$\Longrightarrow$}

\put(22,5){\circle{3}}
\put(23.5, 5){\line(1,0){7}}

\put(32,5){\circle{3}}
\put(30.5,4){$\times$}

\put(32.5,6.5){\line(1,1){5}}
\put(33,6){\line(1,1){5}}
\put(33.5,5.5){\line(1,1){5}}

\put(39,12){\circle{3}}

\end{picture}}

\newcommand{\fad}{
\setlength{\unitlength}{1mm}
\begin{picture}(35,15)(0,0)

\put(1.5, 8){{\scriptsize  0}}
\put(11.5, 8){{\scriptsize  1}}
\put(19.5, 8){{\scriptsize  2+3}}
\put(31.5, 8){{\scriptsize  4}}
\put(38.5, 15){{\scriptsize  3}}

\put(2, 5){\circle{3}}
\put(3.5, 5){\line(1,0){7}}

\put(12,5){\circle{3}}
\put(13.5,4){$\Longrightarrow$}

\put(22,5){\circle{3}}
\put(23.5, 5){\line(1,0){7}}

\put(32,5){\circle{3}}
\put(30.5,4){$\times$}

\put(32.5,6.5){\line(1,1){5}}
\put(33,6){\line(1,1){5}}
\put(33.5,5.5){\line(1,1){5}}

\put(39,12){\circle{3}}

\end{picture}}

\newcommand{\fba}{
\setlength{\unitlength}{1mm}
\begin{picture}(35,15)(0,0)

\put(1.5, 8){{\scriptsize  1}}
\put(11.5, 8){{\scriptsize  2}}
\put(21.5, 8){{\scriptsize  3}}
\put(35.5, 8){{\scriptsize  4}}
\put(26.5, 15){{\scriptsize  0+1+2+3}}

\put(2, 5){\circle{3}}
\put(3.5, 5){\line(1,0){7}}

\put(12,5){\circle{3}}
\put(13.5, 4.5){\line(1,0){7}}
\put(13.5, 5.5){\line(1,0){7}}

\put(22,5){\circle{3}}
\put(20.5,4){$\times$}
\put(23.5, 5){\line(1,0){11}}

\put(36,5){\circle{3}}
\put(34.5,4){$\times$}

\put(22.5,6.5){\line(1,1){5}}
\put(23.5,5.5){\line(1,1){5}}

\put(29,12){\circle{3}}
\put(27.5,11){$\times$}

\put(30.5,11.5){\line(1,-1){5}}
\put(30,11){\line(1,-1){5}}
\put(29.5,10.5){\line(1,-1){5}}

\end{picture}}

\newcommand{\fbb}{
\setlength{\unitlength}{1mm}
\begin{picture}(35,15)(0,0)

\put(-0.5, 8){{\scriptsize  0+1}}
\put(11.5, 8){{\scriptsize  2}}
\put(21.5, 8){{\scriptsize  3}}
\put(35.5, 8){{\scriptsize  4}}
\put(26.5, 15){{\scriptsize  1+2+3}}

\put(2, 5){\circle{3}}
\put(3.5, 5){\line(1,0){7}}

\put(12,5){\circle{3}}
\put(13.5, 4.5){\line(1,0){7}}
\put(13.5, 5.5){\line(1,0){7}}

\put(22,5){\circle{3}}
\put(20.5,4){$\times$}
\put(23.5, 5){\line(1,0){11}}

\put(36,5){\circle{3}}
\put(34.5,4){$\times$}

\put(22.5,6.5){\line(1,1){5}}
\put(23.5,5.5){\line(1,1){5}}

\put(29,12){\circle{3}}
\put(27.5,11){$\times$}

\put(30.5,11.5){\line(1,-1){5}}
\put(30,11){\line(1,-1){5}}
\put(29.5,10.5){\line(1,-1){5}}

\end{picture}}

\newcommand{\fbc}{
\setlength{\unitlength}{1mm}
\begin{picture}(35,15)(0,0)

\put(1.5, 8){{\scriptsize  0}}
\put(9.5, 8){{\scriptsize  1+2}}
\put(21.5, 8){{\scriptsize  3}}
\put(35.5, 8){{\scriptsize  4}}
\put(26.5, 15){{\scriptsize  2+3}}

\put(2, 5){\circle{3}}
\put(3.5, 5){\line(1,0){7}}

\put(12,5){\circle{3}}
\put(13.5, 4.5){\line(1,0){7}}
\put(13.5, 5.5){\line(1,0){7}}

\put(22,5){\circle{3}}
\put(20.5,4){$\times$}
\put(23.5, 5){\line(1,0){11}}

\put(36,5){\circle{3}}
\put(34.5,4){$\times$}

\put(22.5,6.5){\line(1,1){5}}
\put(23.5,5.5){\line(1,1){5}}

\put(29,12){\circle{3}}
\put(27.5,11){$\times$}

\put(30.5,11.5){\line(1,-1){5}}
\put(30,11){\line(1,-1){5}}
\put(29.5,10.5){\line(1,-1){5}}

\end{picture}}

\newcommand{\fbd}{
\setlength{\unitlength}{1mm}
\begin{picture}(35,15)(0,0)

\put(1.5, 8){{\scriptsize  0}}
\put(11.5, 8){{\scriptsize  1}}
\put(18.5, 8){{\scriptsize  2+3}}
\put(35.5, 8){{\scriptsize  4}}
\put(28.5, 15){{\scriptsize  3}}

\put(2, 5){\circle{3}}
\put(3.5, 5){\line(1,0){7}}

\put(12,5){\circle{3}}
\put(13.5, 4.5){\line(1,0){7}}
\put(13.5, 5.5){\line(1,0){7}}

\put(22,5){\circle{3}}
\put(20.5,4){$\times$}
\put(23.5, 5){\line(1,0){11}}

\put(36,5){\circle{3}}
\put(34.5,4){$\times$}

\put(22.5,6.5){\line(1,1){5}}
\put(23.5,5.5){\line(1,1){5}}

\put(29,12){\circle{3}}
\put(27.5,11){$\times$}

\put(30.5,11.5){\line(1,-1){5}}
\put(30,11){\line(1,-1){5}}
\put(29.5,10.5){\line(1,-1){5}}

\end{picture}}

\newcommand{\fca}{
\setlength{\unitlength}{1mm}
\begin{picture}(40,15)(0,0)

\put(10.5, 12){$*$}
\put(25.5, 12){$*$}

\put(1.5, 8){{\scriptsize  1}}
\put(9.5, 8){{\scriptsize  2+3}}
\put(25.5, 8){{\scriptsize  3}}
\put(33.5, 8){{\scriptsize  0+1+2}}
\put(18.5, 15){{\scriptsize 4}}

\put(2,5){\circle{3}}
\put(3.5, 4.5){\line(1,0){7}}
\put(3.5, 5.5){\line(1,0){7}}

\put(12,5){\circle{3}}
\put(10.5,4){$\times$}
\put(13.5, 4.5){\line(1,0){11}}
\put(13.5, 5.5){\line(1,0){11}}

\put(26,5){\circle{3}}
\put(24.5,4){$\times$}

\put(27.5, 4.5){\line(1,0){7}}
\put(27.5, 5.5){\line(1,0){7}}

\put(36,5){\circle{3}}

\put(13,6){\line(1,1){5}}

\put(19,12){\circle{3}}

\put(20,11){\line(1,-1){5}}

\end{picture}}

\newcommand{\fcb}{
\setlength{\unitlength}{1mm}
\begin{picture}(40,15)(0,0)

\put(10.5, 12){$*$}
\put(25.5, 12){$*$}

\put(-0.5, 8){{\scriptsize 1+2}}
\put(11.5, 8){{\scriptsize  3}}
\put(25.5, 8){{\scriptsize 2+3}}
\put(33.5, 8){{\scriptsize  0+1}}
\put(18.5, 15){{\scriptsize  4}}

\put(2,5){\circle{3}}
\put(3.5, 4.5){\line(1,0){7}}
\put(3.5, 5.5){\line(1,0){7}}

\put(12,5){\circle{3}}
\put(10.5,4){$\times$}
\put(13.5, 4.5){\line(1,0){11}}
\put(13.5, 5.5){\line(1,0){11}}

\put(26,5){\circle{3}}
\put(24.5,4){$\times$}

\put(27.5, 4.5){\line(1,0){7}}
\put(27.5, 5.5){\line(1,0){7}}

\put(36,5){\circle{3}}

\put(13,6){\line(1,1){5}}

\put(19,12){\circle{3}}

\put(20,11){\line(1,-1){5}}

\end{picture}}

\newcommand{\fcc}{
\setlength{\unitlength}{1mm}
\begin{picture}(40,15)(0,0)

\put(10.5, 12){$*$}
\put(25.5, 12){$*$}

\put(1.5, 8){{\scriptsize  2}}
\put(11.5, 8){{\scriptsize  3}}
\put(23.5, 8){{\scriptsize  1+2+3}}
\put(35.5, 8){{\scriptsize  0}}
\put(18.5, 15){{\scriptsize  4}}

\put(2,5){\circle{3}}
\put(3.5, 4.5){\line(1,0){7}}
\put(3.5, 5.5){\line(1,0){7}}

\put(12,5){\circle{3}}
\put(10.5,4){$\times$}
\put(13.5, 4.5){\line(1,0){11}}
\put(13.5, 5.5){\line(1,0){11}}

\put(26,5){\circle{3}}
\put(24.5,4){$\times$}

\put(27.5, 4.5){\line(1,0){7}}
\put(27.5, 5.5){\line(1,0){7}}

\put(36,5){\circle{3}}

\put(13,6){\line(1,1){5}}

\put(19,12){\circle{3}}

\put(20,11){\line(1,-1){5}}

\end{picture}}

\newcommand{\fd}{
\setlength{\unitlength}{1mm}
\begin{picture}(45,15)(0,0)

\put(1.5, 8){{\scriptsize  0}}
\put(11.5, 8){{\scriptsize  1}}
\put(21.5, 8){{\scriptsize  2}}
\put(31.5, 8){{\scriptsize  3}}
\put(41.5, 8){{\scriptsize  4}}

\put(2, 5){\circle{3}}
\put(3.5, 5){\line(1,0){7}}

\put(12,5){\circle{3}}
\put(13.5, 5){\line(1,0){7}}

\put(22,5){\circle{3}}
\put(23.5, 4.5){\line(1,0){7}}
\put(23.5, 5.5){\line(1,0){7}}

\put(32,5){\circle{3}}
\put(30.5,4){$\times$}

\put(33.5,4.5){\line(1,0){7}}
\put(33.5,5){\line(1,0){7}}
\put(33.5,5.5){\line(1,0){7}}

\put(42,5){\circle{3}}

\end{picture}}

\newcommand{\fyokoyaa}{
\setlength{\unitlength}{1mm}
\begin{picture}(5,18)(0,0)

\put(8, 0){\vector(-1,0){8}}
\put(0, 0){\vector(1,0){8}}
\put(3, 2){{\scriptsize  3}}
\end{picture}}

\newcommand{\fyokoyab}{
\setlength{\unitlength}{1mm}
\begin{picture}(5,18)(0,0)

\put(8, 0){\vector(-1,0){8}}
\put(0, 0){\vector(1,0){8}}
\put(3, 2){{\scriptsize  4}}
\end{picture}}

\newcommand{\ftateyaa}{
\setlength{\unitlength}{1mm}
\begin{picture}(5,18)(0,0)

\put(0, 10){\vector(0,-1){10}}
\put(0, 0){\vector(0,1){10}}
\put(2, 4){{\scriptsize  0}}
\end{picture}}

\newcommand{\ftateyab}{
\setlength{\unitlength}{1mm}
\begin{picture}(5,18)(0,0)

\put(0, 10){\vector(0,-1){10}}
\put(0, 0){\vector(0,1){10}}
\put(2, 4){{\scriptsize  1}}
\end{picture}}

\newcommand{\ftateyac}{
\setlength{\unitlength}{1mm}
\begin{picture}(5,18)(0,0)

\put(0, 10){\vector(0,-1){10}}
\put(0, 0){\vector(0,1){10}}
\put(2, 4){{\scriptsize  2}}
\end{picture}}

\newcommand{\fyoncom}{
\setlength{\unitlength}{1mm}
\begin{picture}(150,110)(0,0)

\put(0,90){\faa}\put(42,95){\fyokoyab}\put(17,80){\ftateyaa}

\put(52,90){\fba}\put(92,95){\fyokoyaa}\put(69,80){\ftateyaa}

\put(102,90){\fca}\put(114,80){\ftateyaa}\put(124,80){\ftateyac}

\put(0,60){\fab}\put(42,65){\fyokoyab}\put(17,50){\ftateyab}

\put(52,60){\fbb}\put(92,65){\fyokoyaa}\put(69,50){\ftateyab}

\put(102,60){\fcb}\put(119,50){\ftateyab}

\put(0,30){\fac}\put(42,35){\fyokoyab}\put(17,20){\ftateyac}

\put(52,30){\fbc}\put(92,35){\fyokoyaa}\put(69,20){\ftateyac}

\put(102,30){\fcc}

\put(0,0){\fad}\put(42,5){\fyokoyab}

\put(52,0){\fbd}\put(92,5){\fyokoyaa}

\put(102,0){\fd}

\end{picture}}

\newcommand{\fxi}{
\setlength{\unitlength}{1mm}
\begin{picture}(45,15)(0,0)

\put(2,4){{\it {if}}}

\put(11.5, 8){{\scriptsize  {\it i}}}
\put(21.5, 8){{\scriptsize  {\it j}}}
\put(31.5, 8){{\scriptsize  {\it k}}}
\put(41.5, 8){{\scriptsize  {\it l}}}

\put(12,5){\circle{3}}
\put(13.5, 5){\line(1,0){7}}

\put(22,5){\circle{3}}
\put(23.5, 4.5){\line(1,0){7}}
\put(23.5, 5.5){\line(1,0){7}}

\put(32,5){\circle{3}}
\put(30.5,4){$\times$}

\put(33.5,4.5){\line(1,0){7}}
\put(33.5,5){\line(1,0){7}}
\put(33.5,5.5){\line(1,0){7}}

\put(42,5){\circle{3}}

\end{picture}}

\newcommand{\fxii}{
\setlength{\unitlength}{1mm}
\begin{picture}(45,15)(0,0)

\put(2,4){{\it {if}}}

\put(11.5, 8){{\scriptsize  {\it i}}}
\put(21.5, 8){{\scriptsize  {\it j}}}
\put(31.5, 8){{\scriptsize  {\it k}}}
\put(41.5, 8){{\scriptsize  {\it l}}}

\put(12,5){\circle{3}}
\put(13.5,4){$\Longrightarrow$}

\put(22,5){\circle{3}}
\put(23.5, 5){\line(1,0){7}}

\put(32,5){\circle{3}}
\put(30.5,4){$\times$}

\put(33.5,4.5){\line(1,0){7}}
\put(33.5,5){\line(1,0){7}}
\put(33.5,5.5){\line(1,0){7}}

\put(42,5){\circle{3}}

\end{picture}}

\newcommand{\sgoo}{
\setlength{\unitlength}{1mm}
\begin{picture}(25,12)(0,0)

\put(1.5, 8){{\scriptsize  1}}
\put(11.5, 8){{\scriptsize  2}}
\put(21.5, 8){{\scriptsize  3}}

\put(2, 5){\circle{3}}
\put(0.5, 4){$\times$}
\put(3.5, 5){\line(1,0){7}}

\put(12,5){\circle{3}}
\put(13.5,4){$<$}
\put(16,4.5){\line(1,0){4.5}}
\put(16,5){\line(1,0){4.5}}
\put(16,5.5){\line(1,0){4.5}}

\put(22,5){\circle{3}}

\end{picture}}

\newcommand{\sgoa}{
\setlength{\unitlength}{1mm}
\begin{picture}(35,12)(0,0)

\put(1.5, 8){{\scriptsize  0}}
\put(11.5, 8){{\scriptsize  1}}
\put(21.5, 8){{\scriptsize  2}}
\put(31.5, 8){{\scriptsize  3}}

\put(2, 5){\circle{3}}
\put(3.5, 4.25){\line(1,0){7}}
\put(3.5, 4.75){\line(1,0){7}}
\put(3.5, 5.25){\line(1,0){7}}
\put(3.5, 5.75){\line(1,0){7}}

\put(12, 5){\circle{3}}
\put(10.5, 4){$\times$}
\put(13.5, 5){\line(1,0){7}}

\put(22,5){\circle{3}}
\put(23.5,4){$<$}
\put(26,4.5){\line(1,0){4.5}}
\put(26,5){\line(1,0){4.5}}
\put(26,5.5){\line(1,0){4.5}}

\put(32,5){\circle{3}}

\end{picture}}

\newcommand{\sgoaa}{
\setlength{\unitlength}{1mm}
\begin{picture}(35,10)(0,0)

\put(2, 5){\circle{3}}
\put(3.5, 4.25){\line(1,0){7}}
\put(3.5, 4.75){\line(1,0){7}}
\put(3.5, 5.25){\line(1,0){7}}
\put(3.5, 5.75){\line(1,0){7}}

\put(12, 5){\circle{3}}
\put(10.5, 4){$\times$}
\put(13.5, 5){\line(1,0){7}}

\put(22,5){\circle{3}}
\put(23.5,4){$<$}
\put(26,4.5){\line(1,0){4.5}}
\put(26,5){\line(1,0){4.5}}
\put(26,5.5){\line(1,0){4.5}}

\put(32,5){\circle{3}}

\end{picture}}

\newcommand{\sgoab}{
\setlength{\unitlength}{1mm}
\begin{picture}(25,15)(0,0)

\put(2,5){\circle{3}}
\put(0.5,4){$\times$}
\put(3.5, 4.25){\line(1,0){11}}
\put(3.5, 4.75){\line(1,0){11}}
\put(3.5, 5.25){\line(1,0){11}}
\put(3.5, 5.75){\line(1,0){11}}

\put(16, 5){\circle{3}}
\put(14.5,4){$\times$}

\put(2.5,6.5){\line(1,1){5}}
\put(3,6){\line(1,1){5}}
\put(3.5,5.5){\line(1,1){5}}

\put(9,12){\circle{3}}
\put(7.5,11){$\times$}
\put(10, 11){\line(1,-1){5}}

\put(10.5, 11.5){\line(1,0){7}}
\put(10.5, 12){\line(1,0){7}}
\put(10.5, 12.5){\line(1,0){7}}

\put(19,12){\circle{3}}

\end{picture}}

\newcommand{\sgoac}{
\setlength{\unitlength}{1mm}
\begin{picture}(25,15)(0,0)

\put(2,5){\circle{3}}
\put(0.5,4){$\times$}

\put(3.5, 4.66){\line(1,0){9}}
\put(3.5, 5.34){\line(1,0){9}}
\put(8,4){$\Longrightarrow$}

\put(16, 5){\circle{3}}

\put(2.5,6.5){\line(1,1){5}}
\put(3,6){\line(1,1){5}}
\put(3.5,5.5){\line(1,1){5}}

\put(9,12){\circle{3}}
\put(7.5,11){$\times$}
\put(10, 11){\line(1,-1){5}}

\put(10.5, 11.5){\line(1,0){7}}
\put(10.5, 12){\line(1,0){7}}
\put(10.5, 12.5){\line(1,0){7}}

\put(19,12){\circle{3}}

\end{picture}}

\newcommand{\sgoad}{
\setlength{\unitlength}{1mm}
\begin{picture}(35,10)(0,0)

\put(2, 5){\circle{3}}
\put(3.5, 4.25){\line(1,0){7}}
\put(3.5, 4.75){\line(1,0){7}}
\put(3.5, 5.25){\line(1,0){7}}
\put(3.5, 5.75){\line(1,0){7}}

\put(12, 5){\circle{3}}
\put(10.5, 4){$\times$}

\put(13.5,4.5){\line(1,0){7}}
\put(13.5,5){\line(1,0){7}}
\put(13.5,5.5){\line(1,0){7}}

\put(22,5){\circle{3}}
\put(23.5, 5){\line(1,0){7}}

\put(32,5){\circle{3}}

\end{picture}}

\newcommand{\sgoae}{
\setlength{\unitlength}{1mm}
\begin{picture}(35,10)(0,0)

\put(2, 5){\circle{3}}
\put(3.5, 5){\line(1,0){7}}

\put(12, 5){\circle{3}}

\put(13.5,4.5){\line(1,0){7}}
\put(13.5,5){\line(1,0){7}}
\put(13.5,5.5){\line(1,0){7}}

\put(22,5){\circle{3}}\put(20.5, 4){$\times$}
\put(23.5, 4){$\Longrightarrow$}

\put(32,5){\circle*{3}}

\end{picture}}

\newcommand{\gsancom}{
\setlength{\unitlength}{1mm}
\begin{picture}(110,30)(0,0)
\put(0,15){\sgoaa}
\put(45,15){\sgoab} 
\put(75,15){\sgoac} 
\put(10,0){\sgoad} 
\put(50,0){\sgoae} 
\end{picture}}

\newcommand{\gxiii}{
\setlength{\unitlength}{1mm}
\begin{picture}(25,10)(0,0)

\put(-8,4){\it {if}}  
\put(1.5,8){{\scriptsize {\it i}}}
\put(11.5,8){{\scriptsize {\it j}}}
\put(21.5,8){{\scriptsize {\it k}}}

\put(2,5){\circle{3}}
\put(0.5,4){$\times$}

\put(3.5,5){\line(1,0){7}}

\put(12,5){\circle{3}}
\put(10.5,4){$\times$}

\put(13.5, 4.5){\line(1,0){7}}
\put(13.5, 5){\line(1,0){7}}
\put(13.5, 5.5){\line(1,0){7}}

\put(22,5){\circle{3}}

\end{picture}}

\newcommand{\gxix}{
\setlength{\unitlength}{1mm}
\begin{picture}(25,10)(0,0)

\put(-8,4){\it {if}}  
\put(1.5,8){{\scriptsize {\it i}}}
\put(11.5,8){{\scriptsize {\it j}}}
\put(21.5,8){{\scriptsize {\it k}}}

\put(2, 5){\circle{3}}

\put(3.5,4.5){\line(1,0){7}}
\put(3.5,5){\line(1,0){7}}
\put(3.5,5.5){\line(1,0){7}}

\put(12,5){\circle{3}}\put(10.5, 4){$\times$}
\put(13.5, 4){$\Longrightarrow$}

\put(22,5){\circle*{3}}

\end{picture}}

\newcommand{\gxx}{
\setlength{\unitlength}{1mm}
\begin{picture}(35,10)(0,0)

\put(1.5, 8){{\scriptsize  {\it l}}}
\put(11.5, 8){{\scriptsize {\it k}}}
\put(21.5, 8){{\scriptsize {\it j}}}
\put(31.5, 8){{\scriptsize {\it i}}}

\put(-8, 4){\it {if}}
\put(2, 5){\circle{3}}
\put(3.5, 4.25){\line(1,0){7}}
\put(3.5, 4.75){\line(1,0){7}}
\put(3.5, 5.25){\line(1,0){7}}
\put(3.5, 5.75){\line(1,0){7}}

\put(12, 5){\circle{3}}
\put(10.5, 4){$\times$}
\put(13.5, 5){\line(1,0){7}}

\put(22,5){\circle{3}}
\put(23.5,4){$<$}
\put(26,4.5){\line(1,0){4.5}}
\put(26,5){\line(1,0){4.5}}
\put(26,5.5){\line(1,0){4.5}}

\put(32,5){\circle{3}}

\end{picture}}

\newcommand{\gxxi}{
\setlength{\unitlength}{1mm}
\begin{picture}(35,10)(0,0)

\put(1.5, 8){{\scriptsize  {\it l}}}
\put(11.5, 8){{\scriptsize {\it k}}}
\put(21.5, 8){{\scriptsize {\it j}}}
\put(31.5, 8){{\scriptsize {\it i}}}

\put(-8, 4){\it {if}}

\put(2, 5){\circle{3}}
\put(3.5, 4.25){\line(1,0){7}}
\put(3.5, 4.75){\line(1,0){7}}
\put(3.5, 5.25){\line(1,0){7}}
\put(3.5, 5.75){\line(1,0){7}}

\put(12, 5){\circle{3}}
\put(10.5, 4){$\times$}

\put(13.5,4.5){\line(1,0){7}}
\put(13.5,5){\line(1,0){7}}
\put(13.5,5.5){\line(1,0){7}}

\put(22,5){\circle{3}}
\put(23.5, 5){\line(1,0){7}}

\put(32,5){\circle{3}}

\end{picture}}

\newcommand{\aii}{
\setlength{\unitlength}{1mm}
\begin{picture}(17,12)(0,4)

\put(2,4){{\it {if}}}

\put(11.5, 8){{\scriptsize  {\it i}}}

\put(12,5){\circle{3}}
\put(10.5,4){$\times$}

\put(16,4){{\it {,}}}

\end{picture}}

\newcommand{\aiv}{
\setlength{\unitlength}{1mm}
\begin{picture}(50,12)(0,4)

\put(2,4){{\it {if}}}

\put(11.5, 8){{\scriptsize  {\it i}}}
\put(21.5, 8){{\scriptsize  {\it j}}}
\put(31.5, 8){{\scriptsize  {\it k}}}

\put(15, 7){{\scriptsize  $-x$}}
\put(26, 7){{\scriptsize  $x$}}

\put(10.5,4){$\times$}
\put(13.5, 5){\line(1,0){7}}

\put(22,5){\circle{3}}
\put(20.5,4){$\times$}
\put(23.5, 5){\line(1,0){7}}

\put(30.5,4){$\times$}

\put(36,4){{\it {($x\ne0$),}}}

\end{picture}}

\newcommand{\av}{
\setlength{\unitlength}{1mm}
\begin{picture}(37,12)(0,4)

\put(2,4){{\it {if}}}

\put(11.5, 8){{\scriptsize  {\it i}}}
\put(21.5, 8){{\scriptsize  {\it j}}}
\put(31.5, 8){{\scriptsize  {\it k}}}

\put(12,5){\circle{3}}
\put(10.5,4){$\times$}
\put(13.5, 5){\line(1,0){7}}

\put(22,5){\circle{3}}
\put(20.5,4){$\times$}
\put(23.5,4){$\Longleftarrow$}

\put(32,5){\circle{3}}

\put(36,4){{\it {,}}}

\end{picture}}

\newcommand{\avi}{
\setlength{\unitlength}{1mm}
\begin{picture}(47,12)(0,4)

\put(2,4){{\it {if}}}

\put(11.5, 8){{\scriptsize  {\it i}}}
\put(21.5, 8){{\scriptsize  {\it j}}}
\put(31.5, 8){{\scriptsize  {\it k}}}
\put(41.5, 8){{\scriptsize  {\it l}}}

\put(10.5,4){$\times$}
\put(13.5, 5){\line(1,0){7}}

\put(22,5){\circle{3}}
\put(23.5, 5){\line(1,0){7}}

\put(32,5){\circle{3}}
\put(30.5,4){$\times$}
\put(33.5,4){$\Longleftarrow$}

\put(42,5){\circle{3}}

\put(46,4){{\it {,}}}

\end{picture}}

\newcommand{\avii}{
\setlength{\unitlength}{1mm}
\begin{picture}(133,12)(0,4)

\put(2,4){{\it {if}}}

\put(9, 4){{\scriptsize  {\it j}}}
\put(28, 4){{\scriptsize  {\it k\,\,($ab\ne 0$)}}}
\put(18, 14){{\scriptsize {\it i}}}

\put(11.5,5){.}
\put(13.5, 5){\line(1,0){11}}

\put(13.5, 9){{\scriptsize  $a$}}
\put(23.5, 9){{\scriptsize  $b$}}
\put(14.5, 1){{\scriptsize  $-a-b$}}

\put(25.5,5){.}
\put(13,6){\line(1,1){5}}

\put(18.5,12){.}
\put(20,11){\line(1,-1){5}}

\put(45, 4){\mbox{and}\,
$p(\al_i)p(\al_j)+p(\al_i)p(\al_j)
+p(\al_i)p(\al_j)\equiv1,$}

\end{picture}}

\newcommand{\aviii}{
\setlength{\unitlength}{1mm}
\begin{picture}(58,12)(0,4)

\put(11.5, 8){{\scriptsize  {\it i}}}
\put(21.5, 8){{\scriptsize  {\it j}}}
\put(32, 11){{\scriptsize  {\it k}}}
\put(32, -3){{\scriptsize  {\it l}}}
\put(35, 4){$(x\ne \pm 1,\,0),$}

\put(2,4){{\it {if}}}

\put(12,5){\circle{3}}
\put(13.5, 5){\line(1,0){7}}
\put(13.5, 7){{\scriptsize  $x+1$}}

\put(22,5){\circle{3}}
\put(20.5,4){$\times$}
\put(23,6){\line(1,1){5}}
\put(26, 7){{\scriptsize  $-1$}}

\put(29,12){\circle{3}}

\put(23,4){\line(1,-1){5}}
\put(26,2){{\scriptsize  $-x$}}

\put(29,-2){\circle{3}}

\end{picture}}

\newcommand{\aix}{
\setlength{\unitlength}{1mm}
\begin{picture}(47,12)(0,4)

\put(2,4){{\it {if}}}

\put(11.5, 8){{\scriptsize  {\it i}}}
\put(21.5, 8){{\scriptsize  {\it j}}}
\put(31.5, 8){{\scriptsize  {\it k}}}
\put(41.5, 8){{\scriptsize  {\it l}}}

\put(12,5){\circle{3}}
\put(13.5,4){$\Longrightarrow$}

\put(22,5){\circle{3}}
\put(23.5, 5){\line(1,0){7}}

\put(32,5){\circle{3}}
\put(30.5,4){$\times$}
\put(33.5, 4.5){\line(1,0){7}}
\put(33.5, 5.5){\line(1,0){7}}

\put(42,5){\circle{3}}

\put(46,4){{\it {,}}}

\end{picture}}

\newcommand{\aixa}{
\setlength{\unitlength}{1mm}
\begin{picture}(47,12)(0,4)

\put(11.5, 8){{\scriptsize  0}}
\put(21.5, 8){{\scriptsize  1}}
\put(31.5, 8){{\scriptsize  2}}
\put(41.5, 8){{\scriptsize  3}}

\put(16,7){\scriptsize 2}
\put(26,7){\scriptsize 1}
\put(35.5,7){\scriptsize $-2$}

\put(12,5){\circle{3}}
\put(13.5,4){$\Longrightarrow$}

\put(22,5){\circle{3}}
\put(23.5, 5){\line(1,0){7}}

\put(32,5){\circle{3}}
\put(30.5,4){$\times$}
\put(33.5, 4.5){\line(1,0){7}}
\put(33.5, 5.5){\line(1,0){7}}

\put(42,5){\circle{3}}

\put(46,4){{\it {.}}}

\end{picture}}

\newcommand{\ax}{
\setlength{\unitlength}{1mm}
\begin{picture}(30,15)(0,4)

\put(2,4){{\it {if}}}

\put(11.5, 8){{\scriptsize  {\it i}}}
\put(25.5, 8){{\scriptsize  {\it j}}}
\put(18.5, 15){{\scriptsize {\it k}}}

\put(12,5){\circle{3}}
\put(10.5,4){$\times$}
\put(13.5, 5){\line(1,0){11}}

\put(26,5){\circle{3}}
\put(24.5,4){$\times$}

\put(13,6){\line(1,1){5}}

\put(19,12){\circle{3}}

\put(20,11){\line(1,-1){5}}

\put(30,4){{\it {,}}}

\end{picture}}

\newcommand{\aivu}{
\setlength{\unitlength}{1mm}
\begin{picture}(33,10)(0,4)

\put(2,4){{\it {if}}}

\put(11.5, 8){{\scriptsize  {\it i}}}
\put(21.5, 8){{\scriptsize  {\it j}}}
\put(31.5, 8){{\scriptsize  {\it u}}}

\put(10.5,4){$\times$}
\put(13.5, 5){\line(1,0){7}}

\put(22,5){\circle{3}}
\put(20.5,4){$\times$}
\put(23.5, 5){\line(1,0){7}}

\put(30.5,4){$\times$}

\end{picture}}

\begin{document}
\title{On defining relations of the affine Lie superalgebras
and their quantized universal enveloping superalgebras
\thanks{}}
\author{Hiroyuki Yamane\thanks{Dept. of Math., Osaka Univ.,
Toyonaka 560 Japan \newline and DPMMS, Cambridge Univ., 16 Mill Lane,
Cambridge CB2 1SB, UK}}
\date{11th March 1996}
\maketitle

{\bf Introduction.} In this paper, we give defining relations of the affine Lie superalgebras and defining relations of a super-version of the Drinfeld[D1]-Jimbo[J] 
affine quantized enveloping algebras. As a result, we can exactly define
the affine quantized universal enveloping superalgebras with generators and relations. Moreover we give a Drinfeld's realization of $U_h({\hat {sl}}(m|n)^{(1)})$.
\par
For the Kac-Moody Lie algebra $G$, Gabber-Kac [GK] proved the Serre theorem which states that $G$ is defined with the Chevalley generators $H_i$, $E_i$, $F_i$ $(1\leq i \leq \mbox{rank}G)$ and relations
$$
[H_i,H_j]=0,\,[H_i,E_j]=(\al_i,\al_j)E_j,\,[H_i,F_j]=-(\al_i,\al_j)F_j,
$$
$$
[E_i,F_j]=\delta_{ij}H_i,
$$
$$
ad(E_i)^{1-a_{ij}}(E_j)=0,\,ad(F_i)^{1-a_{ij}}(F_j)=0 
$$ where $\{\al_i(1\leq i \leq \mbox{rank}G)\}$ are simple roots of $G$,
$(\,,\,)$ is an invariant form of $G$ and $(a_{ij})$ is the Cartan matrix of $G$. 
We call these relations Serre's relations.
\par
Kac [K2] classified the finite dimensional simple Lie superalgebras, which are
$sl(m|n)$, $osp(m|n)$, $D(2,1;x)$ $(x\ne 0,1)$, $F_4$ and $G_3$. Van de Leur [VdL] classified the Kac-Moody Lie superalgebras $\g$ of finite growth, which are the finite dimensional simple Lie superalgebras and:
\[\begin{array}{c}
\mbox{${\hat {sl}}(m|n)^{(i)}$ $(i=1,2,4)$, ${\widehat {osp}}(m|n)^{(i)}$ $(i=1,2)$,} \\ \mbox{$D(2,1;x)^{(1)}$ $(x\ne 0,1)$, $F_4^{(1)}$, $G_3^{(1)}$.}  
\end{array}\] In
 this paper we call complex infinite dimensional Kac-Moody Lie superalgebras 
of finite growth affine Lie superalgebras.
\par 
Our first result is to give a Serre theorem for the affine Lie superalgebra ${\g}$, i.e., to give defining relations between the Chevalley generators $H_i$, $E_i$, $F_i$.
We give the defining relations associated to each Cartan matrix of $\g$.
(In general, $\g$ does not have a unique Cartan matrix.) To do this, we use Weyl-group-type isomorphisms $\{L_i\}$ between $\g$. 
Let ${\cal H}$ be a Cartan subalgebra of $\g$.
We note that the Cartan matrix defined for ${\cal H}$ 
does not necessarily coincide with the one defined for $L_i({\cal H})$,
 though $\{L_i\}$ are introduced as counterparts of the inner automorphisms 
$\{ \exp(-\mbox{ad}F_i)\exp(-{\frac {2} {(\al_i,\al_i)}}\mbox{ad}E_i)
\exp(-\mbox{ad}F_i)\,\,(1\leq i \leq \mbox{rank}G)\,\}$ of the Kac-Moody Lie algebra $G$.
We introduce another Lie superalgebra $\bg$
associated to each $\g$. We define the Lie superalgebras ${\bg}$ by a universal condition that $\{L_i\}$ can be lifted to isomorphisms 
$\{{\bar L}_i\}$
 between ${\bg}$. We directly calculate defining relations of $\bg$.
In the case of the Kac-Moody Lie algebra $G$, defining relations of ${\bar G}$ 
are given by Serre's relations. However, for $\g$, we need other relations such as
$$
[[E_i,E_j],[E_j,E_k]]=0\quad \mbox{for}\, (\al_i,\al_k)=(\al_j,\al_j)=0,\,
(\al_i,\al_j)=-(\al_j,\al_k)\ne 0.
$$

There is an epimorphism $j:\bg\rightarrow\g$ satisfying 
$j\circ L_i={\bar L}_i\circ j$.
In the case of the Kac-Moody Lie algebra, most of 
the proof of the Gabber-Kac theorem [GK] was to prove $\ker j =0$. \par
However, in the case of the affine Lie superalgebras $\g$, 
$\bg\ne\g$ if and only if $\g={\hat {sl}}(m|n)^{(i)}$ $(i=1,2,4)$ and $m=n$. 
(In the case of $\g={\hat {sl}}(m|m)^{(1)}$, 
$\bg=sl(m|m)\otimes C[t,t^{-1}]\oplus Cc\oplus Cd$ and 
$\ker j = I\otimes C[t,t^{-1}]$
where $I$ is the unit matrix.) Nevertheless, we can also look for defining relations of $\g={\hat {sl}}(m|m)^{(i)}$ because 
we have concretely known ${\hat {sl}}(m|m)^{(i)}$.  
\par Our second result is to give relations of quantized universal enveloping superalgebras $\uhg$ such that, after $h\rightarrow 0$, 
the relations become the defining relations of $\bg$ 
obtained as our first result. Here $\uhg$ is an $h$-adic topological $C[[h]]$-Hopf superalgebra introduced in [Y1]. In [Y1], we
 showed
an existence of a non-degenerate symmetric bilinear form
defined on a Borel part of $\uhg$.\par
Applying the Drinfeld's quantum double construction to $\uhg$ by using
 the bilinear form, we can see that $\uhg$ is topologically free and that the universal $R$-matrix of $\uhg$ exists. \par 
Since $\uog=\uhg/h\uhg$ is a
 cocommutative Hopf 
$C$-superalgebra, applying Milnor-Moor theorem [MM]
to $\uog$, we see that $\uog$ is a universal enveloping superalgebra $U(\g_0)$ of the Lie superalgebra $\g_0={\cal P}(\uog)$ of primitive elements of $\uog$. By definition of $\g$ as the Kac-Moody Lie superalgebra, $\g$ must be a quotient of $\g_0$. On the other hand, we see that $\g_0$ is a quotient of $\bg$ by our second result. Hence, if $\g=\bg$, we see that $U_0(\g)$ coincides with $U(\g)$ and that our relations must be defining relations of $\uhg$ by the topologically freedom of $\uhg$. \par 
Finally, we calculate relations of $U_h({\hat {sl}}(m|m)^{(1)})$ which become ones generating $\ker j$ after $h\rightarrow 0$, while showing
Drinfeld's realization of $U_h({\hat {sl}}(m|n)^{(1)})$ for general $m$, $n$. 
Gathering up the relations and the ones obtained as our second result, we get defining relations of $U_h({\hat {sl}}(m|m)^{(1)})$. To do these, we 
introduce a Braid group action on $U_h({\hat {sl}}(m|n)^{(1)})$ which become
 the action on ${\hat {sl}}(m|n)^{(1)}$ 
defined by $\{L_i\}$
after $h\rightarrow 0$,
and follow Beck's argument [B]. We won't consider
$U_h({\hat {sl}}(m|m)^{(2)})$
or $U_h({\hat {sl}}(m|m)^{(4)})$. \par
Results in this paper have already been announced in [Y2]. The same results for the finite dimensional $A-G$ type simple Lie superalgebras have already been given in [Y1].

\section{Preliminary}
{\bf 1.1.} 
In \S 1, we mainly refer to [K1-2] and [VdL].\par
Let $\g$ be a $C$-vector space with
a direct sum decomposition 
$\g = \g(0)\oplus \g(1) $.
For $X\in \g$, $p(X) \in \{0,\,1\}$ means that
$X \in \g(p(X))$. We call $p(X)$ the parity of $X$.
A Lie superalgebra $\g$ is defined 
with the bilinear map $[\,\,,\,\,] : \g \times \g
 \rightarrow \g$ such that
$$
 [X,Y]=-(-1)^{p(X)p(Y)}[Y,X],
$$
$$
 [X,[Y,Z]]=[[X,Y],Z]+(-1)^{p(X)p(Y)}[Y,[X,Z]].
$$
If a bilinear form $(\,\,|\,\,):\g \times \g
 \rightarrow C$ satisfies $(X|Y)=(-1)^{p(X)p(Y)}(Y|X)$
and $([X,Y]|Z)=(X|[Y,Z])$, then we call it an invariant form. \par
A Lie superalgebra
${\hat \g}=\g\otimes_{C}C[t,t^{-1}]
\oplus Cc \oplus Cd$ is defined by
$$
[X\otimes t^{m}+a_1c+b_1d,\,Y\otimes t^{n}+a_2c+b_2d]
$$ 
$$
=[X,Y]\otimes t^{m+n} + m\delta_{m+n,0}(X|Y)c 
+ b_1nY\otimes t^{n}-b_2mX\otimes t^{m}.
$$ where ${\hat \g}(0) =\g(0)\otimes_{C}C[t,t^{-1}]
\oplus Cc \oplus Cd$ and
${\hat \g}(1) =\g(1)\otimes_{C}C[t,t^{-1}]$.\par
Let $\gamma : \g \rightarrow \g$ be an automorphism
of finite order $r$\,\,(i.e. $\gamma([X,Y])=[\gamma(X),\gamma(Y)]$). 
Put 
$$
\g^\gamma_n = \{X \in \g| 
\gamma(X)=(exp{\frac{2\pi\sqrt{-1}}{n}})X\} 
\quad (0 \le n < r).
$$
Then $\g_0^\gamma$ is a subalgebra of $\g$
and $\g_i^\gamma$ $(1\le i\le r-1)$ is the $\g_0^\gamma$-module.  
We can define a subalgebra ${\hat \g}^{(\gamma)}$ by
$$
{\hat \g}^{(\gamma)}= \bigoplus_{n=0}^{r-1}
(\bigoplus_{m\in Z} \g^\gamma_n \otimes t^{mr+n})
\oplus Cc \oplus Cd\,.
$$ Obviously ${\hat \g}^{(1)} = {\hat \g}$.

{\bf 1.2.} Here we introduce a definition of the Kac-Moody Lie
superalgebra 
in an abstract manner similar to the abstract definition of 
the Kac-Moody Lie algebra given in [K1;1.3]. Let ${\cal E}$ be a finite 
dimensional $C$ vector space with a nondegenerate symmetric bilinear
 form $(\,,\,):{\cal E}\times{\cal E}\rightarrow C$. Let
 $\Pi = \{\alpha_0,...,\alpha_1\}$ be a finite linearly independent 
subset of ${\cal E}$. We call an element $\alpha_i \in \Pi$ the simple
root.
Let $p:\Pi\rightarrow Z/2Z$ be a function. We call $p$  
the parity function. Put ${\cal H} = {\cal E}^*$. We call ${\cal H}$ the
Cartan subalgebra. We identify an element $\gamma \in {\cal E}$ with 
$H_\gamma \in {\cal H}$ satisfying $\delta(H_\gamma)=(\delta,\gamma)$
$(\delta \in {\cal E})$. For a datum $\Epip$, we define a 
Lie superalgebra 
${\widetilde \g} = {\widetilde \g}\Epip$  
with the generators $E_i$, $F_i$ $(0 \leq i \leq n)$, $H \in {\cal H}$,
the parities $p(E_i)=p(F_i)=p(\al_i)$, $p(H)=0$
 and the relations:

$${[H, H']} = 0 \quad (H, H' \in {\cal H}), \eqno{(1.2.1)}$$
$${[H, E_i]} = \al_i(H)E_i,\,\, [H, F_i]=-\al_i(H)F_i, \eqno{(1.2.2)}$$
$${[E_i,F_j]} = \delta_{ij}H_{\al_i}. \eqno{(1.2.3)}$$

Then we have the triangular decomposition
${\widetilde \g}=
{\widetilde {\cal N}}^+\oplus
{\cal H}\oplus{\widetilde {\cal N}}^-$. Here ${\widetilde {\cal N}}^+$
(resp. ${\widetilde {\cal N}}^-$) is the free superalgebra 
with generators $E_i$ (resp. $F_i$).
Define the quotient Lie superalgebra 
$\g = \g\Epip$ 
of
${\widetilde \g} = {\widetilde \g}\Epip$ 
by
\[
{\g\Epip} = 
{\widetilde \g}\Epip/{r}\,\,.
\]
where $r=r\Epip$ is the ideal which is maximal of the ideals
$r_1$ satisfying 
$r_1 \cap {\cal H} = 0$. 
We call $\g = \g\Epip$
the Kac-Moody Lie superalgebra. We 
have the triangular decomposition
\[
\g = {\cal N}^{+}\oplus{\cal H}\oplus {\cal N}^{-}
\]
where ${\cal N}^{+}$ and ${\cal N}^{-}$ are
the subalgebras of $\g$ generated by $E_i$ and $F_i$ 
respectively. Let $r_{\pm}=r \cap {\widetilde {\cal N}}^{\pm}$. Then we have
$r=r_{-}\oplus r_{+}$ and 
${\cal N}^{\pm}={\widetilde {\cal N}}^{\pm}/r_{\pm}$. 
We also have ${\widetilde {\cal N}}^{+}\cong {\widetilde {\cal N}}^{-}$ 
($E_i \leftrightarrow F_i$).
 Let
 $\g = {\cal H} \oplus 
(\bigoplus_{{\al}\in{\cal E}}\g_{\al})$ be the root space
decomposition where $\g_{\al} = \{X\in\g|[H,X]=\al(H)X\,
(H\in{\cal H})\}$. Let $r_\al = r\cap \g_{\al}$. 
We put $\Phi = \Phi\Epip = 
\{\al\in {\cal E}| \g_{\al} \ne 0\}$. Let 
 $P_{+}=
Z_{+}\al_0\oplus\cdots\oplus Z_{+}\al_n$ and 
$\Phi_{+}=\Phi\cap P_{+}$, $\Phi_{-}=-\Phi_{+}$. Then we have
$\Phi = \Phi_{+} \cup \Phi_{-}$.
Clearly,  
we have $r = \bigoplus_{{\al}\in\Phi}r_\al$. \par
 For $\beta, \al \in P_+$, 
we say $\beta\leq\al$ if $\al -\beta \in P_+$.
Let $r_{+, \leq\al}$ be the ideal of ${\widetilde {\cal N}}^+$ generated by 
$r_\beta$ with $\beta\leq\al$. Then 
$r_{+}$ $ = \cup_{\gamma \in P_{+}} r_{+,\leq\gamma}$.
By the same argument in [K1], we have:\newline
{\bf Proposition 1.2.1.} {\it For $\Epip$, let $\rho\in {\cal E}$ be 
an element such that $(\rho,\al_i)=(\al_i,\al_i)/2$.
If $\al\in P_+$ satisfies $(\al,\al)\ne 2(\rho,\al)$,
then $r_\al$ is included in the ideal of ${\widetilde {\cal N}}^+$
generated by $r_\beta$ such that $\beta\geq\al$.
In particular,
$$
 r_+ = \bigcup_{\gamma \in P_{+},
\,(\gamma,\gamma) = 2(\rho,\gamma) } r_{+,\leq\gamma}\quad .
\eqno{(1.2.4)}
$$ }
 \newline\newline
{\bf Lemma 1.2.1.} {\it For $\al_i\in\Pi$, we have
$\dim\g_{\al_i}=1$. If $p(\al_i)=1$
and $(\al_i,\al_i)\ne 0$, then 
 $\dim\g_{2\al_i}=1$.} \newline
{\bf Proof.} By $[E_i.F_i]=H_{\al_i}\ne 0$,
$\dim\g_{\al_i}\ne 0$. Hence 
$\dim\g_{\al_i}= 1$. Similarly 
we can get a proof of the latter half.
\begin{flushright}
Q.E.D.
\end{flushright}
$$
$$
{\bf 1.3.} 
Here we introduce the Dynkin diagram $\Gamma$
for a datum $\Epip$. We first prepare the three-type 
dots:\setlength{\unitlength}{1mm}
\begin{center}
\begin{picture}(60,30)(0,0)

\put(05, 15){\circle{5}\,,}
\put(25, 15){\circle{5}\,,}
\put(23, 13){\line(1,1){4}}
\put(23, 17){\line(1,-1){4}}
\put(45, 15){\circle*{5}}
\put(50, 15){.}

\end{picture}
\end{center}
 We call those the white dot, the gray dot and the black dot 
respectively. To the $i$-th simple root $\al_i$, we put the
corresponding $i$-th dot defined such that:
 
\setlength{\unitlength}{1mm}
\begin{center}
\begin{picture}(70,40)(0,0)

\put(05, 25){\circle{5}}
\put(15, 24){if $(\al_i,\al_i)\ne 0$ and $p(\al_i)=0$\,,}
\put(05, 15){\circle{5}}
\put(03, 13){\line(1,1){4}}
\put(03, 17){\line(1,-1){4}}
\put(15, 14){if $(\al_i,\al_i)=0$ and $p(\al_i)=1$\,,}
\put(05, 05){\circle*{5}}
\put(15, 04){if $(\al_i,\al_i)\ne 0$ and $p(\al_i)=1$\,.}

\end{picture}
\end{center} 
The dot $\dotcross$ stands for $\dotwhite$
or $\dotgray$.
The dot $\dotwhiteblack$ stands for $\dotwhite$
or $\dotblack$.
We write a 
$|(\al_i,\al_j)|$-line between the i-th dot and the j-th dot
or write as follows:
 
\setlength{\unitlength}{1mm}
\begin{center}
\begin{picture}(40,13)(0,0)

\put(05, 04){$\cdot$}\put(04.5, 08){$i$}
\put(07, 05){\line(1,0){26}}
\put(13, 9){$(\al_i,\al_j)$}
\put(35, 04){$\cdot$}\put(34, 08){$j$}

\end{picture}
\end{center}
Moreover we add a pile pointing to the smaller of 
$|(\al_i,\al_i)|$ and $|(\al_j,\al_j)|$. If   
$|(\al_i,\al_i)|=0$ or $|(\al_j,\al_j)|=0$, then
we sometimes omit the pile. The semilines $\semiline$
stands for $\noline$ or $\twoline$. \vspace{1.5cm}
If $(\al_i,\al_j)\ne 0$, $(\al_i,\al_k)\ne 0$ and
$\displaystyle{{\frac {|(\al_i,\al_j)|} {(\al_i,\al_j)}}\cdot 
{\frac {|(\al_i,\al_k)|} {(\al_i,\al_k)}}=1}$, then we put $*$ in a sector enclosed by an edge between $i$-th and $j$-th dots and an edge between $i$-th and $k$-th dots. Namely we describe the situation as $\hoshi$.\vspace{1.5cm} However we sometimes omit $*$.
\newline\newline
{\bf 1.4.} We have already known:\newline
{\bf Theorem 1.4.1}[K2].{\it 
The Kac-Moody Lie superalgebra $\g\Epip$
is finite dimensional as a $C$-vector space 
if and only if the datum is one of the following Dynkin diagram.
(In any diagram, there is an only one dot whose parity is odd.)}
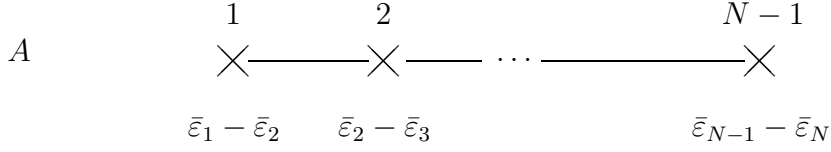
\begin{center}
Diagram 1.4.1
\end{center}

\begin{center}
\setlength{\unitlength}{1mm}
\begin{picture}(50,15)(0,0)
  
\put(-15, 3){$A_{N-1}$}
\put(2, 5){\circle{3}}
\put(3.5, 5){\line(1,0){4}}
\put(9.5, 4){$\cdots$}
\put(16.5, 5){\line(1,0){4}}
\put(22, 5){\circle{3}}\put(20.5, 4){$\times$}
\put(23.5, 5){\line(1,0){4}}
\put(29.5, 4){$\cdots$}
\put(36.5, 5){\line(1,0){4}}
\put(42, 5){\circle{3}}

\end{picture}

\setlength{\unitlength}{1mm}
\begin{picture}(50,15)(0,0)
  
\put(-10,3){$B_N$}
\put(2, 5){\circle{3}}
\put(3.5, 5){\line(1,0){4}}
\put(9.5, 4){$\cdots$}
\put(16.5, 5){\line(1,0){4}}
\put(22, 5){\circle{3}}\put(20.5, 4){$\times$}
\put(23.5, 5){\line(1,0){4}}
\put(29.5, 4){$\cdots$}
\put(36.5, 5){\line(1,0){4}}
\put(42, 5){\circle{3}}
\put(43.5, 4){$\Longrightarrow$}

\put(52,5){\circle{3}}

\end{picture}

\setlength{\unitlength}{1mm}
\begin{picture}(50,15)(0,0)
  
\put(-10, 3){$B_N$}
\put(2, 5){\circle{3}}
\put(3.5, 5){\line(1,0){14}}
\put(19.5, 4){$\cdots$}
\put(26.5, 5){\line(1,0){14}}

\put(42, 5){\circle{3}}
\put(43.5, 4){$\Longrightarrow$}

\put(52,5){\circle*{3}}

\end{picture}

\setlength{\unitlength}{1mm}
\begin{picture}(50,15)(0,0)
  
\put(-10, 3){$C_N$}
\put(2, 5){\circle{3}}\put(0.5, 4){$\times$}
\put(3.5, 5){\line(1,0){14}}
\put(19.5, 4){$\cdots$}
\put(26.5, 5){\line(1,0){14}}

\put(42, 5){\circle{3}}
\put(43.5, 4){$\Longleftarrow$}

\put(52,5){\circle{3}}

\end{picture}

\setlength{\unitlength}{1mm}
\begin{picture}(50,25)(0,0)
  
\put(-10, 13){$D_N$}
\put(2, 15){\circle{3}}
\put(3.5, 15){\line(1,0){4}}
\put(9.5, 14){$\cdots$}
\put(16.5, 15){\line(1,0){4}}
\put(22, 15){\circle{3}}\put(20.5, 14){$\times$}
\put(23.5, 15){\line(1,0){4}}
\put(29.5, 14){$\cdots$}
\put(36.5, 15){\line(1,0){4}}
\put(42, 15){\circle{3}}

\put(43,16){\line(1,1){5}}
\put(49,22){\circle{3}}

\put(43,14){\line(1,-1){5}}
\put(49, 8){\circle{3}}

\put(48, 14){$*$}

\end{picture}

\setlength{\unitlength}{1mm}
\begin{picture}(50,25)(0,0)

\put(-10, 13){$D(2;1,x)$}  
\put(22, 15){\circle{3}}\put(20.5, 14){$\times$}

\put(23,16){\line(1,1){5}}\put(22,20){\tiny{$-1$}}
\put(29,22){\circle{3}}

\put(23,14){\line(1,-1){5}}\put(22,8){\tiny{$x$}}
\put(29, 8){\circle{3}}

\put(32, 15){$(x\ne 0,\,1)$}

\end{picture}

\setlength{\unitlength}{1mm}
\begin{picture}(35,15)(0,0)

\put(-10, 3){$F_4$}
\put(2, 5){\circle{3}}
\put(3.5, 5){\line(1,0){7}}

\put(12,5){\circle{3}}
\put(13.5,4){$\Longrightarrow$}

\put(22,5){\circle{3}}
\put(23.5, 5){\line(1,0){7}}

\put(32,5){\circle{3}}
\put(30.5,4){$\times$}

\end{picture}

\setlength{\unitlength}{1mm}
\begin{picture}(25,12)(0,0)
  
\put(-10, 3){$G_3$}
\put(2, 5){\circle{3}}
\put(0.5, 4){$\times$}
\put(3.5, 5){\line(1,0){7}}

\put(12,5){\circle{3}}
\put(13.5,4){$<$}
\put(16,4.5){\line(1,0){4.5}}
\put(16,5){\line(1,0){4.5}}
\put(16,5.5){\line(1,0){4.5}}

\put(22,5){\circle{3}}

\end{picture}
\end{center}

{\bf 1.5.} In 1.5-13, we give a concrete form of 
finite dimensional $\g\Eopip$
of $\Eopip$ of $ABCD$-type, which is given by using matrices.

Let ${\tcE_0}$ is an 
${\tilde N}$-dimensional
$C$-vector space with a nondegenerate symmetric bilinear form
$(\,,\,):{\tcE_0}\times{\tcE_0}\rightarrow C$.
Let $e \in \{\pm 1\}$.
Let $\be_i$ $(1\le n \le N)$ be the orthogonal basis of ${\tcE_0}$
satisfying 
$$
(\be_i,\be_j)=e\cdot\delta_{ij}\cdot (-1)^{\tp (i)}\,\,.
$$
where $\tp (i)$ is 0 or 1. Let $\bd_i=(\be_i,\be_i)$. 
Let $gl(\tcE_0,e)=gl(\tcE_0)$ be the $C$-linear space of 
${\tilde N}\times {\tilde N}$
-matrices.
Put $E_{ij}= (\delta_{xi},\delta_{yi})_
{1\le x,\,y \le {\tilde N}} \in gl(\tcE_0)$\,\, 
$(1\le i,\,j \le {\tilde N})$. We regard $gl(\tcE_0,e)$
as a superspace with a 
parity  $p$ defined by 
$p(E_{ij})=(-1)^{(\tp (i)+\tp (j))}$.
Then $gl(\tcE_0,e)$ can be regarded as a Lie superalgebra defined by
$$
[X,Y]=XY
-(-1)^{p(X)p(Y)}YX\,.
$$ Define a $C$-linear map $str:gl(\cE_0,e)\rightarrow C$
by $str(E_{ij})=\delta_{ij}\cdot \bd_i$.
Let $sl(\cE_0,e)$ denote the subalgebra of $gl(\cE_0,e)$ of the
elements $X \in gl(\cE_0,e)$ satisfying $str(X)=0$.
Let $\g\Eopip$ be a finite dimensional Kac-Moody Lie superalgebra. 
Let $\gamma : \g\Eopip \rightarrow \g\Eopip$ be an automorphism
of finite order $r$.
In 1.6-13, We give a concrete form
of an affine $ABCD$-type superalgebra
$\g\Epip$ arising from 
${\hat \g}\Eopip^{(\gamma)}$. Let $\Eopipgamma$ be
 the datum of $\g\Eopip_0^\gamma$, i.e., $\g\Eopipgamma=\g\Eopip_0^\gamma$.
For $\gamma$, we define the datum $\Epip=\Eopip^{(\gamma)}$ 
 of affine type as follows:
\newline\newline
(i) \quad\quad ${\cal E}= \cE_0^\gamma\oplus C\dl \oplus C\Lambda_0$
\newline\newline where 
$(x,\dl)=(x,\Lambda_0)=0$\, $(x \in \cE_0^\gamma)$,
 $(\dl,\dl)=$ $(\Lambda_0,\Lambda_0)=0$
and $(\dl,\Lambda_0)=1$.
\newline\newline
(ii) For the lowest weight $\psi$
 of $\g\Eopip_1^\gamma$,
let $\al_0 = \dl + \psi$ 
and $\Pi=\{\al_0\} \cup \Pi_0^\gamma$.
\newline\newline
(iii) Let $e_\psi\in\g\Eopip_1^\gamma$ be a weight vector of $\psi$. We define the parity $p(\al_0)$ of $\al_0$ by the parity of $e_\psi$ in $\g\Eopip$.
We also define 
$p(\al_i)=p^\gamma_0(\al_i)$ $(\al_i\in\Pi_0^\gamma)$ 
$$
$$
{\bf 1.6.} Here we put 
$N={\tilde N}$,\,$\cE_0=\tcE_0$,\,
$n=N-1$,  $N \ge 2$ \, and $e=1$.
Let $\Eopip$ be the datum whose Dynkin diagram is:

\begin{center}
Diagram 1.6.1
\end{center}
\begin{center}
\DiagramA
\end{center}
Then we can identify $\g\Eopip$ with 
$gl(\cE_0,e)$
where
$H_{\be_j}=\bd_jE_{jj}$,\,
$E_i=E_{ii+1}$,\,
$F_i=\bd_iE_{i+1i}$
\,. We also note that $H_{\al_i}=$
$\bd_iE_{ii}-\bd_{i+1}E_{i+1i+1}$.\par
We also note the lowest root of $gl(\cE_0,e)$ is $\theta=\be_N-\be_1$. 
 Lowest and highest root vectors $E_{\theta}$, $F_{\theta}$ 
$(\in gl(\cE_0,e))$ satisfying
$[E_{\theta},F_{\theta}]=H_\theta(=\bd_NE_{NN}-\bd_1E_{11})$ are
given by
$$
F_{\theta} = \bd_NE_{1N}, E_{\theta} = E_{N1}.
$$
\newline
{\bf Proposition 1.6.1.}  If $\sum_{i=1}^{N} \bd_i \ne 0$,
then $sl(\cE_0,e)$ is the simple Lie superalgebra. 
If $\sum_{i=1}^{N} \bd_i = 0$, then the quotient
$sl(\cE_0,e)/{\cal I}$ is the simple Lie superalgebra
where
\begin{center}
${\cal I}$\,=\,{\bC}
$\cdot\sum_{i=1}^{N-1}((\sum_{j=1}^{i}\bd_j)\cdot H_{\al_i})$ 
\,=\,{\bC}$\cdot\sum_{i=1}^{N}E_{ii}$\,.
\end{center}
\vspace{1cm} 
{\bf Proposition 1.6.2.} {\it Assume $N \ge 3$ if $p(\al_0)=1$. Let $\Eopip$ be the datum of Diagram 1.6.1 and $\Epip=\Eopip^{(I)}$. Let
$\bbg\Epip ={\hat {sl}}(\cE_0,e) + \cE_0\,\, (\subset 
{\hat {gl}}(\cE_0,e))$. There is an epimorphism 
$j:\bbg\Epip\rightarrow\g\Epip$ defined by letting 
$j(H\otimes 1 + ac + bd)=H + aH_\dl + bH_{\Lambda_0}
\, (H\in{\cal H})$ 
and $j(E_\theta\otimes t)=E_0$,  $j(F_\theta\otimes t^{-1})=F_0$. \par
If $\sum_{i=1}^{N} \bd_i \ne 0$, then $j$ is an isomorphism.   
If $\sum_{i=1}^{N} \bd_i = 0$, then 
$$
\ker j = \oplus_{i \ne 0} {\cal I}\otimes t^i.
$$}
\newline\newline
{\bf Example 1.6.1.} The Dynkin diagram 
of $\Epip$ in Proposition 1.6.2 is:
\begin{center}
Diagram 1.6.2. $(N \ge 2)$
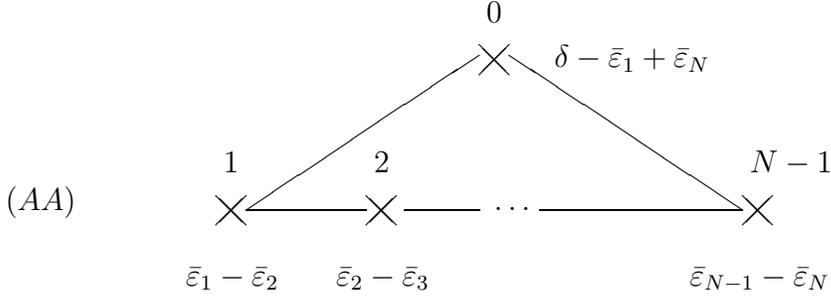
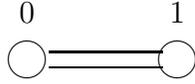
\DiagramAA
\end{center}
\hspace{2cm} \DiagramAAa
\begin{flushright} 
($\sum_{i=0}^{N-1} p(\al_i)\equiv 0$).\\
\end{flushright}

$$
$$ 
\newline\newline\newline{\bf 1.7.}\, Here we assume
\[{\tilde N}= \left\{
 \begin{array}{cl}
  2N+1 &\quad\mbox{if $\tilde N$ is odd,}\\
  2N &\quad\mbox{if $\tilde N$ is even.}\\
 \end{array}\right. \]
{\it We assume ${\tilde N} \ge 3$.}
 Let $i'={\tilde N}+1-i$. We also assume
$$
\tp (i) = \tp(i')\,\, (1\le i \le {\tilde N})
$$
and
\begin{center}
$\tp (N+1) = 0$  if ${\tilde N}$ is odd.
\end{center}
Let $g_i$\,\,$(1\le i \le {\tilde N})$ be such that
$g_i\in \{\pm 1\}$ and $g_ig_{i'}=(-1)^{\tp (i)}$.
We assume that $g_{N+1}=1$ if ${\tilde N}$ is odd.
\par  Let $\Omega$ be an automorphism of $gl(\tcE_0,e)$ of order 2
defined by:
$$
\Omega(X)_{ij}=-(-1)^{(\tp (i)\tp (j)+
\tp (j))}g_ig_jX_{j'i'}\,.
$$

Denote $sl(\tcE_0,e)_0^\Omega$ by $osp(\tcE_0,e)$.
We can identify $osp(\tcE_0,e)$ with
a finite dimensional 
Kac-Moody Lie superalgebra
$\g\Eopip$ of a datum $\Eopip$.

Here we note $n=N$. In 1.8-11 we will give:\par
(1) the datum $\Eopip$. \par 
(2) the lowest root $\theta$ of $sl(\tcE_0,e)_0^\Omega$,
a lowest root vector $E_{\theta}$ 
and a highest root vector $F_{\theta}$
such that $[E_{\theta},F_{\theta}]=H_\theta$.\par
(3) the lowest weight $\psi$ of $sl(\tcE_0,e)_1^\Omega$,
a lowest weight vector $E_{\psi}$ 
and a highest weight vector $F_{\psi}$
such that $[E_{\psi},F_{\psi}]=H_\psi$.\par
(4) $\g\Epip$ and $\bbg\Epip$ arising from
${\widehat{osp}}(\tcE_0,e)^{(I)}$ and ${\hat{sl}}(\tcE_0,e)^{(\Omega)}$.
$$
$$
{\bf 1.8.} $B$-type. If ${\tilde N}$ is odd, then 
the Dynkin diagram of $\Eopip$ is:

\begin{center}
Diagram 1.8.1.
\DiagramB
\end{center}
Here
$H_{\be_j}=\bd_j(E_{jj}-E_{j'j'})$\,\,$(1 \le j \le N)$,\,
$E_i=E_{ii+1}-(-1)^{\tp (i)\tp (i+1)}(-1)^{\tp (i+1)}
g_ig_{i+1}E_{(i+1)'i'}$
\,$(1 \le i \le N-1)$,\,
$E_N=E_{NN+1}-g_NE_{(N+1)'N'}$,\,
$F_i=e\cdot \{(-1)^{\tp (i)}E_{i+1i}
-(-1)^{\tp (i)\tp (i+1)}
g_ig_{i+1}E_{i'(i+1)'}\}$\,$(1 \le i \le N-1)$,
\,$F_N=e\cdot \{(-1)^{\tp (N)}E_{N+1N}-g_NE_{N'(N+1)'}\}$
.\par Moreover we have:\newline
If $p(\al_1)+\cdots+p(\al_{N-1})\equiv 0$, then 
$$
\theta=-\be_1-\be_2, 
F_{\theta}=e\cdot \{(-1)^{\tp (1)}E_{21'}
-(-1)^{\tp (1)\tp (2)}g_2g_{1'}E_{12'}\},
$$
$$
 E_{\theta}=E_{1'2}
-(-1)^{\tp (1)\tp (2)}(-1)^{\tp (2)}g_{1'}g_2E_{2'1},
$$
$$
\psi=-2\be_1, F_{\psi}=e\cdot 2(-1)^{\tp (1)}E_{11'},
E_{\psi}=E_{1'1}.
$$
If
 $p(\al_1)+\cdots+p(\al_{N-1})\equiv 1$, then 
$$
\theta=-2\be_1, F_{\theta}=e\cdot 2(-1)^{\tp (1)}E_{11'},
E_{\theta}=E_{1'1}.
$$
$$
\psi=-\be_1-\be_2, 
F_{\psi}=e\cdot \{(-1)^{\tp (1)}E_{21'}
+(-1)^{\tp (1)\tp (2)}g_2g_{1'}E_{12'}\},
$$
$$
 E_{\psi}=E_{1'2}
+(-1)^{\tp (1)\tp (2)}(-1)^{\tp (2)}g_{1'}g_2E_{2'1},
$$
\newline 
{\bf 1.9.} $C$-type.
If ${\tilde N}$ is even and $\tp (N)=1$, $e=-\bd_N$, then 
the Dynkin diagram of $\Eopip$ is:

\begin{center}
Diagram 1.9.1.
\DiagramC
\end{center}
Here 
$H_{\be_j}=\bd_j(E_{jj}-E_{j'j'})$\,\,$(1 \le j \le N)$,\,
$E_i=E_{ii+1}-(-1)^{(\tp (i)\tp (i+1) +\tp (i+1))}
g_ig_{i+1}E_{(i+1)'i'}$
\,$(1 \le i \le N-1)$,\,
$E_N=E_{NN'}$,\,
$F_i=e\cdot \{(-1)^{\tp (i)}E_{i+1i}
-(-1)^{\tp (i)\tp (i+1)}
g_ig_{i+1}E_{i'(i+1)'}\}$\,$(1 \le i \le N-1)$,
\,$F_N=e\cdot 2(-1)^{\tp (N)}E_{N'N}$
.\par Moreover we have:\newline
If $p(\al_1)+\cdots+p(\al_{N-1})\equiv 1$, then 
$$
\theta=-\be_1-\be_2, 
F_{\theta}=e\cdot \{(-1)^{\tp (1)}E_{21'}
-(-1)^{\tp (1)\tp (2)}g_2g_{1'}E_{12'}\},
$$
$$
 E_{\theta}=E_{1'2}
-(-1)^{\tp (1)\tp (2)}(-1)^{\tp (2)}g_{1'}g_2E_{2'1},
$$
$$
\psi=-2\be_1, F_{\psi}=e\cdot 2(-1)^{\tp (1)}E_{11'},
E_{\psi}=E_{1'1}.
$$
If
 $p(\al_1)+\cdots+p(\al_{N-1})\equiv 0$, then 
$$
\theta=-2\be_1, F_{\theta}=e\cdot 2(-1)^{\tp (1)}E_{11'},
E_{\theta}=E_{1'1},
$$
$$
\psi=-\be_1-\be_2, 
F_{\psi}=e\cdot \{(-1)^{\tp (1)}E_{21'}
+(-1)^{\tp (1)\tp (2)}g_2g_{1'}E_{12'}\},
$$
$$
 E_{\psi}=E_{1'2}
+(-1)^{\tp (1)\tp (2)}(-1)^{\tp (2)}g_{1'}g_2E_{2'1}.
$$
\newline 
{\bf 1.10.} $D$-type.
If ${\tilde N}$ is even and $\tp (N)=0$, $e=\bd_N$, then 
the Dynkin diagram of $\Eopip$ is:

\begin{center}
Diagram 1.10.1.
\DiagramD
\end{center}
Here
$H_{\be_j}=\bd_j(E_{jj}-E_{j'j'})$\,\,$(1 \le j \le N)$,\,
$E_i=E_{ii+1}-(-1)^{(\tp (i)\tp (i+1) +\tp (i+1))}
g_ig_{i+1}E_{(i+1)'i'}$
\,$(1 \le i \le N-1)$,\,
$E_N=E_{N-1N'}-g_{N-1}g_{N'}E_{N(N-1)'}$,\newline
$F_i=e\cdot \{(-1)^{\tp (i)}E_{i+1i}
-(-1)^{\tp (i)\tp (i+1)}
g_ig_{i+1}E_{i'(i+1)'}\}$\quad$(1 \le i \le N-1)$,
\,$F_N=e\cdot \{(-1)^{\tp (N-1)}E_{N'N-1}
-g_{(N-1)'}g_{N'}E_{(N-1)'N}\}$
.\par Moreover we have:\newline
If $p(\al_1)+\cdots+p(\al_{N-1})\equiv 0$, then 
$$
\theta=-\be_1-\be_2, F_{\theta}=e\cdot \{(-1)^{\tp (1)}E_{21'}
-(-1)^{\tp (1)\tp (2)}g_2g_{1'}E_{12'}\},
$$
$$
 E_{\theta}=E_{1'2}
-(-1)^{\tp (1)\tp (2)}(-1)^{\tp (2)}g_{1'}g_2E_{2'1},
$$
$$
\psi=-2\be_1, F_{\psi}=e\cdot 2(-1)^{\tp (1)}E_{11'},
E_{\psi}=E_{1'1}.
$$
If
 $p(\al_1)+\cdots+p(\al_{N-1})\equiv 1$, then 
$$
\theta=-2\be_1, F_{\theta}=e\cdot 2(-1)^{\tp (1)}E_{11'},
E_{\theta}=E_{1'1}.
$$
$$
\psi=-\be_1-\be_2, 
F_{\psi}=e\cdot \{(-1)^{\tp (1)}E_{21'}
+(-1)^{\tp (1)\tp (2)}g_2g_{1'}E_{12'}\},
$$
$$
 E_{\psi}=E_{1'2}
+(-1)^{\tp (1)\tp (2)}(-1)^{\tp (2)}g_{1'}g_2E_{2'1},
$$\newline 
{\bf 1.11.} 
{\bf Proposition 1.11.1.} $osp(\tcE_0,e)$ is the simple Lie superalgebra.
\vspace{1cm}\newline 
{\bf Proposition 1.11.2.} 
{\it Let $\Eopip$ be the datum of $osp(\tcE_0,e)$ and $\Epip=\Eopip^{(I)}$.
Put $\bbg\Epip={\widehat{osp}}(\tcE_0,e)^{(I)}$. There is an isomorphism: 
$j:\bbg\Epip\rightarrow\g\Epip$ defined by letting 
$j(H\otimes 1 + ac + bd)=H + aH_\dl + bH_{\Lambda_0}
\, (H\in{\cal H})$ 
and $j(E_\theta\otimes t)=E_0$,  $j(F_\theta\otimes t^{-1})=F_0$.}
\vspace{0.5cm}\newline
{\bf Proposition 1.11.3.} {\it If $\sum_{i=1}^{{\tilde N}} (-1)^{\tp (i)} \ne
0$,
then $sl(\tcE_0,e)_1^\Omega$ is the simple $osp(\tcE_0,e)$-module. 
If $(-1)^{\tilde N}=1$ and $\sum_{i=1}^{N} \bd_i = 0$, then the quotient
$sl(\tcE_0,e)_1^\Omega/{\cal I}$ is the simple $osp(\tcE_0,e)$-module
where}
\begin{center}
${\cal I}$\,=\,{\bC}$\cdot\sum_{i=1}^{N}
(E_{ii}+E_{i'i'})$\,.
\end{center}
$$
$$
{\bf Proposition 1.11.4.} 
{\it Assume $N \ge 3$ if $\al_N=2\be_N$ and $p(\al_1)=1$.
 Let $\Eopip$ be the datum of $sl(\tcE_0,e)$ and $\Epip=\Eopip^{(\Omega)}$.
Put $\bbg\Epip={\widehat{sl}}(\tcE_0,e)^{(\Omega)}$. There is an epimorphism: 
$j:{\hat{sl}}(\tcE_0,e)^{(2)}\rightarrow\g\Epip$ defined by letting 
$j(H\otimes 1 + ac + bd)=H + aH_\dl + bH_{\Lambda_0}
\, (H\in{\cal H})$ 
and $j(E_\psi\otimes t)=E_0$,  $j(F_\psi\otimes t^{-1})=F_0$. \par
If $\sum_{i=1}^{N} \bd_i \ne 0$, then $j$ is an isomorphism.   
If $(-1)^{\tilde N}=1$ and  $\sum_{i=1}^{N} \bd_i = 0$, then 
$$
\ker j = \oplus_{i} {\cal I}\otimes t^{2i-1}.
$$.}
\newline\newline\newline
{\bf Example 1.11.1.} The Dynkin diagrams of $\Epip$
in Proposition 1.11.2 (resp. Proposition 1.11.4) are:

\begin{center}
Diagram 1.11.1 $(N \ge 2)$
\end{center}
\vspace{0.5cm}
\begin{center}
\DiagramDB
\end{center}\vspace{1cm}
\hspace{2cm}\DiagramDBb
\begin{flushright} 
($\sum_{i=1}^Np(\al_i)\equiv 0$ (resp. 1)) \\
\end{flushright} 

\begin{center}
Diagram 1.11.2 $(N \ge 1)$
\end{center}
\begin{center}
\DiagramCB
\end{center}
\hspace{2cm}\DiagramCBb
\begin{flushright} 
($\sum_{i=1}^Np(\al_i)\equiv 1$ (resp. 0)) \\
\end{flushright} 

\begin{center}
Diagram 1.11.3 $(N \ge 3)$
\end{center}
\vspace{0.5cm}
\begin{center}
\DiagramDC
\end{center}
\begin{flushright} 
($\sum_{i=1}^{N-1}p(\al_i)\equiv 1$ (resp. 0)) \\
\end{flushright} 

\begin{center}
Diagram 1.11.4 $(N \ge 3)$
\end{center}
\vspace{0.5cm}
\begin{center}
\DiagramCC
\end{center}
\begin{flushright} 
($\sum_{i=1}^{N-1}p(\al_i)\equiv 0$ (resp. 1)) \\
\end{flushright} 

\begin{center}
Diagram 1.11.5 $(N \ge 3)$
\end{center}
\vspace{0.5cm}
\begin{center}
\DiagramDD
\end{center}
\hspace{2cm}\DiagramDDd
\begin{flushright} 
($\sum_{i=1}^{N-1}p(\al_i)\equiv 0$ (resp. 1)) \\
\end{flushright} 

\begin{center}
Diagram 1.11.6 $(N \ge 3)$
\end{center}
\vspace{0.5cm}
\begin{center}
\DiagramCD
\end{center}
\begin{flushright} 
($\sum_{i=1}^{N-1}p(\al_i)\equiv 1$ (resp. 0)) \\
\end{flushright}

{\bf 1.12.} Keep the notations in 1.10. However we denote the integer
$N$ in 1.10 by $N_1$. In 1.12, we let $N$ denotes $N_1-1$. {\it 
We assumed
${\tilde N} \ge 3$. Then $N \ge 2$.}
Let $\omega$ be an automorphism of $osp(\tcE_0,,\bd_{N_1})$ of order 2
defined by:
$$
\omega(X)_{ij}=X_{{\hat i}{\hat j}}\,
$$where 
\[{\hat i}= \left\{
 \begin{array}{cl}
  N' &\quad\mbox{if $i=N_1$,}\\
  N &\quad\mbox{if $i=N_1'$,}\\
  i &\quad\mbox{otherwise.}\\
 \end{array}\right. \]
Then we can identify 
$osp(\tcE_0,,\bd_{N_1})_0^\omega$ 
with $\g\Eopip$ of $\Eopip$ of
the Dynkin diagram:

\begin{center}
Diagram 1.12.1. 
\end{center}
\begin{center}
\DiagramB
\end{center} 
$$
$$ Here
$H_{\be_j}=$ $\bd_j(E_{jj}-E_{j'j'})$ $(1 \le j \le N-1)$,\,
$E_i=$ $E_{ii+1}-$\newline $(-1)^{(\tp (i)\tp (i+1) +\tp (i+1))}
g_ig_{i+1}E_{(i+1)'i'}$
\,$(1 \le i \le N-2)$,\,
$E_{N_1-1}=$ $E_{N_1-1N_1}+E_{N_1-1N_1'}
$ $-g_{N_1-1}g_{N_1'}(E_{N_1(N_1-1)'}+E_{N_1'(N_1-1)'})$,\newline
$F_i=e\cdot \{(-1)^{\tp (i)}E_{i+1i}
-(-1)^{\tp (i)\tp (i+1)}
g_ig_{i+1}E_{i'(i+1)'}\}$\quad$(1 \le i \le N_1-2)$,
\,$F_{N_1-1}=\frac{1}{2}\cdot e\cdot \{(-1)^{\tp
(N_1-1)}(E_{N_1N_1-1}+E_{N_1'N_1-1})
-g_{(N_1-1)'}g_{N_1'}(E_{(N_1-1)'N}+E_{(N_1-1)'N_1'})\}$.
$$
$$
Here we introduce the lowest weight $\psi$ of $osp(\tcE_0,e)_1^\omega$,
a lowest weight vector $E_{\psi}$ 
and a highest weight vector $F_{\psi}$
such that $[E_{\psi},F_{\psi}]=H_\psi$:
$$
\psi=-\be_1, 
F_{\psi}=\frac{1}{2}\cdot e\cdot \{E_{1N_1}+E_{1N_1'}
-g_1g_{N_1}(E_{N_11}+E_{N_1'1})\},
$$
$$
 E_{\psi}=E_{N1}+E_{N'1}-(-1)^{\tp(1)}g_Ng_{1'}(E_{1'N}+E_{1'N'}).
$$\newline
{\bf Proposition 1.12.2.} 
{\it $osp(\tcE_0,e)_1^\omega$
is a simple $osp(\tcE_0,e)_0^\omega$-module.}
\newline\newline
{\bf Proposition 1.12.3.} 
{\it Let $\Eopip$ be the datum of $osp(\tcE_0,e)$
 and $\Epip=\Eopip^{(\omega)}$.
Put $\bbg\Epip={\widehat{osp}}(\tcE_0,e)^{(\omega)}$.
There is an isomorphism: 
$j:\bbg\Epip$ 
$\rightarrow\g\Epip$ defined by letting 
$j(H\otimes 1 + ac + bd)=H + aH_\dl + bH_{\Lambda_0}
\, (H\in{\cal H})$ 
and $j(E_\psi\otimes t)=E_0$,  $j(F_\psi\otimes t^{-1})=F_0$.}
\newline\newline\newline
{\bf Example 1.12.1.} If $N_1 \ge 3$, the Dynkin diagrams of $\Epip$
in Proposition 1.12.1 are:
\begin{center}
Diagram 1.12.2. $(N \ge 1)$
\end{center}
\begin{center}
\DiagramBB
\end{center}
\hspace{2cm}\DiagramBBb
\begin{flushright} 
($\sum_{i=1}^Np(\al_i)\equiv 0$) \\
\end{flushright} 
$$
$$\newline 
{\bf 1.13.} Let $\tcE_0$ be the $C$-vector space in 1.5. Here we assume
that ${\tilde N}=2N+2$ $(N \ge 1)$
and shift the numbering of the basis as follows
$$
\{\be_0, \be_1,\ldots, \be_N, \be_{N+1}, \be_{N'},\ldots, \be_{1'}\}
$$ 
where $i'=2N-i+2$. We also assume:
$$
\tp(0)=1,\,\tp(N+1)=0,\,\tp(i)=\tp(i')\,(1\le i \le N),
$$ and

$$
g_i,\, g_{i'} \in \{\pm 1\}\,(1\le i \le N) \quad g_ig_{i'}=(-1)^{\tp(i)}.
$$
In 1.13, we denote this $\tcE_0$ by $\tcE_{01}$ and we again denote the
subspace with 
the basis $\{\be_1,\ldots, \be_N, \be_{N+1}, \be_{N'},\ldots, \be_{1'}\}$  
by $\tcE_{0}$.
We denote an element $X\in sl(\tcE_{01},e)$ by
\[ X=\left(
 \begin{array}{c|ccc}
  \alpha &\,& a_i &\, \\
  \hline
      &\,&       &\, \\
  b_i &\,&X_{ij} &\, \\
      &\,&       &\, \\
 \end{array}
   \right) \]
where the sizes of the matrices $(\alpha)$, $(a_i)$, $(b_i)$ and $(X_{ij})$
are
$1\times1$, $1\times(2N+1)$,
$(2N+1)\times1$ and $(2N+1)\times(2N+1)$ respectively. Let 
$\Xi$ be an automorphism of $sl(\tcE_{01},e)$ of order 4
defined by:
\[\Xi(X)=\left(
 \begin{array}{c|ccc}
  -\alpha &\,& -\sqrt{-1}g_ib_{i'} &\, \\
  \hline
      &\,&       &\, \\
  -\sqrt{-1}g_{i'}a_{i'} &\,&\Omega((X_{ij}))&\, \\
      &\,&       &\, \\
 \end{array}
   \right). \]
Then $sl(\tcE_{01},e)_0^\Xi$ consists of the matrices
\[ X=\left(
 \begin{array}{c|ccc}
  0 &\,& 0 &\, \\
  \hline
      &\,&       &\, \\
  0 &\,&X_{ij} &\, \\
      &\,&       &\, \\
 \end{array}
   \right) \]
where the $2N+1\times2N+1$-matrices
$(X_{ij})$ form $osp(\tcE_0,e)$ whose Dynkin diagram is:
\begin{center}
Diagram 1.13.1
\DiagramB
\end{center}
Here we introduce the lowest weight $\psi$ of $sl(\tcE_{01},e)^\Xi$,
a lowest weight vector $E_{\psi} \in sl(\tcE_{01},e)_1^\Xi$ 
and a highest weight vector $F_{\psi} \in sl(\tcE_{01},e)_3^\Xi$
such that $[E_{\psi},F_{\psi}]=H_\psi$:
$$
\psi=-\be_1, 
F_{\psi}=e\cdot \{(-1)^{\tp(1)}E_{01'}+g_1E_{1'0}\},
$$
$$
 E_{\psi}=E_{1'0}-g_1E_{01'}.
$$\newline
{\bf Proposition 1.13.1.}
$sl(\tcE_{01},e)_1^\Xi$
and $sl(\tcE_1,e)_3^\Xi$ are the simple
$osp(\tcE,e)=sl(\tcE_{01},e)_0^\Xi$-modules. 
\vspace{1cm}\newline
{\bf Proposition 1.13.2.} As an $osp(\tcE_0,e)$-module, 
$$
 sl(\tcE_1,e)_2^\Xi \cong sl(\tcE,e)_1^\Omega \oplus {\cal I}
$$ where ${\cal I} = C\sum_{i=1}^{2N+1}\{(-1)^{\tp(i)}E_{00}+E_{ii}\}$.
\newline\newline
{\bf Proposition 1.13.3.} 
{\it 
 Let $\Eopip$ be the datum of $sl(\tcE_{01},e)$ and $\Epip=\Eopip^{(\Omega)}$.
Put $\bbg\Epip={\widehat{sl}}(\tcE_{01},e)^{(\Xi)}$. There is an epimorphism: 
$j:{\hat{sl}}(\tcE_0,e)^{(4)}\rightarrow\g\Epip$ defined by letting 
$j(H\otimes 1 + ac + bd)=H + aH_\dl + bH_{\Lambda_0}
\, (H\in{\cal H})$ 
and $j(E_\psi\otimes t)=E_0$,  $j(F_\psi\otimes t^{-1})=F_0$. \par
If $\sum_{i=1}^{N} \bd_i \ne 0$, then $j$ is an isomorphism.   
If $(-1)^{\tilde N}=1$ and  $\sum_{i=1}^{N} \bd_i = 0$, then 
$$
\ker j = \oplus_{i} {\cal I}\otimes t^{4i-2}.
$$.}\newline\newline\newline
{\bf Example 1.13.1.} The Dynkin diagrams of $\Epip$
in Proposition 1.13.2 are:
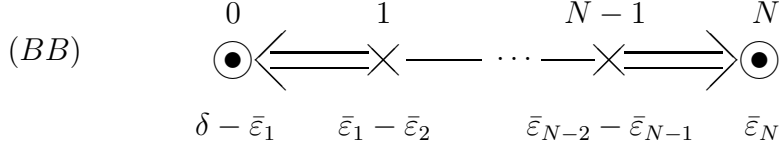
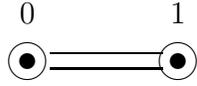
\begin{center}
Diagram 1.13.2. $(N \ge 1)$
\end{center}
\begin{center}
\DiagramBB
\end{center}
\hspace{2cm}\DiagramBBb
\begin{flushright} 
($\sum_{i=1}^Np(\al_i)\equiv 1$) \\
\end{flushright} 
{\bf 1.14.} We call the datum $\Epip$ in the following table 
the one of {\it {affine ABCD-type}}. In the following table, a subscript of a name of a Dynkin diagram shows $\displaystyle{\sum_{i=0}^{n}p(\al_i)
\,(\mbox{mod\,2})}$ of the corresponding superalgebra. 

\[\renewcommand{\arraystretch}{1.6}
\begin{array}{|c|c|c|c|} \noalign{\hrule height0.8pt}
\mbox{Name} & \mbox{Another name} & \mbox{Dynkin diagram} 
& \g=\bbg \\ 
\hline
A_{N-1}^{(1)} & \hat{sl}(\tcE)^{(1)} & (AA)_0 & 
 \displaystyle{\sum_{i=1}^{N}\bd_i\ne 0} \\
\hline
B_N^{(1)} & \widehat{osp}(\tcE_{odd})^{(1)} & (DB)_0\,(CB)_1 & 
all \\
\hline
A_{2N}^{(2)} & \hat{sl}(\tcE_{odd})^{(2)} & (DB)_1\,(CB)_0 & 
all \\
\hline
C_N^{(1)}\, D_N^{(1)} & 
\widehat{osp}(\tcE_{even})^{(1)} & (CC)_0\,(CD)_1\,(DD)_0\,(DC)_1 & 
all \\
\hline
A_{2N-1}^{(2)} & \hat{sl}(\tcE_{even})^{(2)} & (CC)_1\,(CD)_0\,(DD)_1\,(DC)_0 & 
\displaystyle{\sum_{i=1}^{N}\bd_i\ne 0} \\
\hline
D_{N+1}^{(2)} & \widehat{osp}(\tcE_{even})^{(2)} & (BB)_0
 & all \\
\hline
A_{2N+1}^{(4)} & \hat{sl}(\tcE)^{(4)} & (BB)_1
 & \displaystyle{\sum_{i=1}^{N}\bd_i\ne 0} \\
\noalign{\hrule height0.8pt}
\end{array} \]
{\bf 1.15.} We call the data whose Dynkin diagrams are Diagram 5.1.4,
Diagram 5.2.3 and Diagram 5.3.3 $D(2;1,x)^{(1)}$-type, $F_4^{(1)}$-type
and $G_3^{(1)}$-type respectively. We call these {\it
{affine exceptional}} type.\newline\newline
{\bf {2. Affine Weyl type isomorphism}}\newline
{\bf {2.1.}} In this section,  we introduce a family $\{ L_i\}$ of isomorphisms between
Lie superalgebras $\g\Epip$ of $\Epip$ of affine $ABCD$-type (see also [FSS]). 
The isomorphism $L_i$ can be considered
as a super-version of Weyl group action. However our isomorphisms 
change Dynkin diagrams. We shall also introduce other Lie superalgebras 
 $\bg\Epip$ of $\Epip$ of affine ABCD-type. We shall show that $\bg\Epip$'s
are universal superalgebras satisfying: (i) $\g\Epip$ is a quotient of $\bg\Epip$, (ii) $L_i$ can be lifted to isomorphisms of $\bg\Epip$'s. Finally we shall show that $\bg\Epip=\bbg\Epip$. We have already introduced $\bbg\Epip$ in a concrete way in \S 1 for $\Epip$ of affine $ABCD$-type. We have known that $\bbg\Epip\ne\g\Epip$ if and only if $\g\Epip=\hat{sl}(\tcE)^{(1)},\,\hat{sl}(\tcE_{even})^{(2)},\,\hat{sl}(\tcE)^{(4)}$ and  $\displaystyle{\sum_{i=1}^N \bd_i =0}$.\par Our idea using Weyl-group-type isomorphisms relates to [LS].   
\newline
{\bf {2.2.}}  Let
$$ {\breve {E}}_{ij}(k)=
[...[[E_j,\underbrace{
 {E_i],E_i],\ldots,E_i]}}_{k-
 \mbox{ times}}, \quad
 {\breve {F}}_{ij}(k)=
[...[[F_j,\underbrace{
 {F_i],F_i],\ldots,F_i]}}_{k-
 \mbox{ times}}.
$$ 

The calculations of the following lemma are useful. \newline\newline 
{\bf Lemma 2.2.1.} (i) If $i\ne j$, then
$[[E_j,E_i],F_i]=-(\al_i,\al_j)E_j$, 
$[E_i,[F_j,F_i]]=(-1)^{p(\al_i)p(\al_j)}(\al_i,\al_j)F_j$
$[[E_j,E_i],[F_j,F_i]]=$
$(-1)^{p(\al_i)p(\al_j)}(\al_i,\al_j)H_{\al_i+\al_j}$.\newline
(ii) For the
$i$-th simple root $\al_i\in\Pi$  satisfying \,$(\al_i, \al_i) \ne 0$,
let $\displaystyle{a_{ij}={\frac {2(\al_j,\al_i)} {(\al_i,\al_i)}}}$. Assume
 that $a_{ij}$ is an even integer if $p(\al_i)=1$. Put
\[ <k;-a_{ij}> = \left\{
 \begin{array}{cl}
k(-a_{ij}-k+1)  &\quad\mbox{if $p(\al_i)=0$,}\\
k  &\quad\mbox{if $p(\al_i)=1$ and $k$ is even,}\\
-a_{ij}-k+1  &\quad\mbox{if $p(\al_i)=1$ and $k$ is odd.}\\
\end{array}\right.\]
We put $<k;-a_{ij}>!=\prod_{r=1}^k<r;-a_{ij}>$.
Then

$$
[E_i,{\breve {F}}_{ij}(k)]
=-(-1)^{p(\al_i)p(\al_j)}<k;-a_{ij}>{\breve {F}}_{ij}(k-1), 
$$

$$
[{\breve {E}}_{ij}(k),F_i]
=(-1)^{(k-1)p(\al_j)}<k;-a_{ij}>{\breve {E}}_{ij}(k-1), 
$$

$$
[{\breve {E}}_{ij}(k), {\breve {F}}_{ij}(k)]
=(-1)^k(-1)^{p(\al_i)p(\al_j)}<k;-a_{ij}>!\,H_{k\al_i+\al_j}. 
$$

$$
$$
{\bf Lemma 2.2.2.} {\it Let $\g=\g\Epip$ be a Kac-Moody 
superalgebra with the triangular decomposition
$\g={\cal {\cal N}^+}\otimes{\cal H}\otimes{\cal {\cal N}^-}$.
If $X\in{\cal {\cal N}^+}$ (resp. $Y\in{\cal {\cal N}^-}$)
satisfies $[X,F_k]=0$ (resp. $[E_k,Y]=0$) for any $k$,
then $X=0$ (resp. $Y=0$) in $\g$.}\newline
{\bf Proof.} We can assume that $X$ is in an root space. Let $r_+(X)$ be the ideal of ${\cal N}^+$
generated by $X$. Then $r_+(X)$ is an ideal of $\g$
such that $r_+(X) \cap {\cal H} =0$. Hence $X=0$.
\par \hfill Q.E.D. \newline\newline
As an immediate consequence of Lemma 2.2.1 and Lemma 2.2.2, we have:
\newline\newline
{\bf Lemma 2.2.3} {\it (i) For the
$i$-th simple root $\al_i\in\Pi$  satisfying \,$(\al_i, \al_i) \ne 0$,
let $a_{ij}=$ $2(\al_j,\al_i)/(\al_i,\al_i)$.
Assume that $a_{ij}$ is an even integer if $p(\al_i)=1$. Then 
${\breve {E}}_{ij}(-a_{ij})=0$,
${\breve {F}}_{ij}(-a_{ij})=0$ in $\g$. \newline
(ii) If $(\al_i, \al_i) = 0$, then $[E_i,E_i]=[F_i,F_i]=0$
in $\g$.}
\newline
{\bf Proof.} Direct calculations. \newline\newline
{\bf Proposition 2.2.1.}  
{\it Let $\g=\g\Epip$ be a Kac-Moody superalgebra such that an
$i$-th simple root $\al_i\in\Pi$  satisfies \,$(\al_i, \al_i) = 0$.
Let $\Pi'=\{\al_1',\ldots\al_n' \}$ $=\{-\al_i, \al_j + \al_i\,$
$(j\ne i, (\al_i,\al_j)\ne 0),\,$ 
$\al_j \, (j\ne i, (\al_i,\al_j)= 0)\, \}$. 
Put $\g'=\g({\cal E}, \Pi', p)$.
Then there are 
isomorphisms $\phi:\g\rightarrow
 \g'$ such that
$$
\phi(H_\gamma)=H_\gamma. \eqno{(2.2.1)}
$$(In particular,
\[ \quad  \phi(H_\gamma) = \left\{
 \begin{array}{rl}
-H_{\al_i'} &\quad\mbox{if $\gamma=\al_i$,}\\ 
H_{\al_i'+\al_j'} &\quad\mbox{if $\gamma=\al_j$ and $(\al_i,\al_j)\ne 0$,}\\ 
H_{\al_j'} &\quad\mbox{if $\gamma=\al_j$ and $(\al_i,\al_j)= 0$.}\\ 
 \end{array}\right.  \]) 

$$
\phi(E_i)=-(-1)^{p(\al_i)}F_i,\quad\phi(F_i)=-E_i,  \eqno{(2.2.2)}
$$
 
$$
\phi(E_j)=
-{\frac {(-1)^{p(\al_i)p(\al_j)}} {(\al_i,\al_j)}}[E_j,E_i],\quad
\phi(F_j)=
-[F_j,F_i],\,\,(\,\,(\al_i,\al_j)\ne 0\,\,) \eqno{(2.2.3)}
$$
$$
\phi(E_j)=E_j,\quad
\phi(F_j)=F_j\quad(\,\,i\ne j,\,\,(\al_i,\al_j)= 0\,\,) \eqno{(2.2.4)}
$$
}
{\bf Proof.} Let ${\cal H}'$ be the Cartan subalgebra of $\g'$. 
Denote the right hand sides of (2.2.1-4)
by $H_\gamma'$, $E'_j$ and $F'_j$.
We can show that there is an epimorphism 
$y:{\tilde \g}\Epip\rightarrow\g'$.
such that  
$y(H_\gamma)=H_\gamma'$, 
$y(E_j)=E'_j$ and $y(F_j)=F'_j$. 
For example, by Lemma 2.2.3 (ii), we can show that the elements satisfy 
(1.2.1-3). Clearly 
$y_{|{\cal H}}:{\cal H}\rightarrow{\cal H}'$ 
 is isomorphism.
Hence there exists an epimorphism
$\phi_1:\g'\rightarrow\g$
such that $\phi_1(H'_\gamma)=H_\gamma$, 
$\phi_1(E'_j)=E_j$ and $\phi_1(F'_j)=F_j$.
Since ${\phi_1}_{|{\cal H}'}$ is injective, $\phi_1$ is isomorphism.
$\phi_1$ is nothing else but $\phi^{-1}$.
\par\hfill Q.E.D. \newline\newline
Keep the notations in the statement of Proposition 2.2.1.
We still assume $(\al_i,\al_i)=0$.
For the Dynkin diagram $\Gamma$ of $\Epip$, we denote
the Dynkin diagram of $({\cal E}, \Pi', p)$ by $\Gamma^{<i>}$.
Similarly to the proof 
of Proposition 2.2.1, we have following Propositions 2.2.2-3.\newline\newline
{\bf Proposition 2.2.2.}  
{\it Let $\g=\g\Epip$ be a Kac-Moody superalgebra such that an
$i$-th simple root $\al_i\in\Pi$  satisfies \,$(\al_i, \al_i) \ne 0$.
Let $d_i=(\al_i,\al_i)/2$. Let 
\[ (-a_{ij})_s! = \left\{
 \begin{array}{rl}
(-a_{ij})!  &\quad\mbox{if $p(\al_i)=0$,}\\
({\frac {-a_{ij}} {2}})! 2^{{\frac {-a_{ij}} {2}}}
 &\quad\mbox{if $p(\al_i)=1$,}\\ \end{array}\right.\]
Assume that $a_{ij}$ is an even integer if $p(\al_i)=1$.
Then there is an 
isomorphisms $\phi:\g\rightarrow\g$ such that
$$
\phi(H_\gamma)=H_{\gamma-
{\frac {2(\gamma,\al_i)} {(\al_i,\al_i)}}\al_i}\eqno{(2.2.7)}
$$
$$
\phi(E_i)=-(-1)^{p(\al_i)}F_i,\quad\phi(F_i)=-E_j,\eqno{(2.2.8)}
$$
$$
\phi(E_j)= {\frac {1} {(-a_{ij})_s! d_i^{-a_{ij}}}}
{\breve {E}}_{ij}(-a_{ij}),\eqno{(2.2.9)}
$$
$$
\phi(F_j)= (-1)^{-a_{ij}}{\frac {1} {(-a_{ij})_s!}}
{\breve {F}}_{ij}(-a_{ij}).\eqno{(2.2.10)}
$$}
{\bf Proposition 2.2.3.} {\it For 
$\Pi=\{\al_1,\ldots,\al_n\}$, let $a:\{1,\ldots,n\}$
$\rightarrow\{1,\ldots,n\}$ be a bijective map
such that $(\al_{a(i)},\al_{a(j)})=$ $(\al_i,\al_j)$
 $(1 \le i,\, j \le n)$.
Then there is an 
isomorphisms $\phi:\g\rightarrow\g$ such that
$$
H_{\al_{i}}=H_{\al_{a(i)}},\eqno{(2.2.11)} 
$$
$$
\phi(E_i)=E_{a(i)}, \phi(F_i)=F_{a(i)}.\eqno{(2.2.12)}
$$} $$ $$
{\bf {2.3.}} Here we fix a positive integer $N$. Let $\Theta_N$ be the set of affine $ABCD$-type Dynkin diagrams $\Gamma$
satisfying that the number of the dots of $\Gamma$ is $N-1$ 
if $\Gamma$ is of the affine $A$-type, $N$ otherwise. 
Let $D(\Theta_N)$ be the set of the data $({\cal E},
\Pi=\{\al_0, \al_1,\ldots, \al_n \}, p)$
whose Dynkin diagrams belong to $\Theta_N$. 
Let $\Gamma
\in \Theta_N$ and $\Epip \in D(\Theta_N)$.
For $0 \le i \le n$, we define $\Gamma^{\sigma(i)}
=\Gamma' \in
\Theta_N$
and $({\cal E}^{\sigma(i)}=$ $\oplus_{i=1}^N C\be_i'\oplus C\dl
\oplus C\Lambda_0,$
 $\Pi^{\sigma(i)}, p^{\sigma(i)})$ 
$\in D(\Theta_N)$ by following (i)-(iii):
 (We put $\bd_i'=(\be_i',\be_i')$.)
\newline\newline
(i) If $1 \le i \le N-1$ and $(\al_i,\al_i) = 0$,
then $\Epipsi$ satisfies that $\Gamma^{\sigma(i)}=\Gamma^{<i>}$
 (read the sentences before Proposition 2.2.2) and 
$\bd_i'=\bd_{i+1}$, $\bd_{i+1}'=\bd_i$, $\bd_j'=\bd_j$
$(j\ne i, i+1)$. \newline\newline
(ii) If $\Gamma$ is affine $A$ type and $i=0$, then $\Epipsi$ satisfies that 
$\Gamma^{\sigma(i)}=\Gamma^{<i>}$ and 
$\bd_1'=\bd_N$, $\bd_N'=\bd_1$ and $\bd_j'=\bd_j$
($j\ne 1, N$).\newline\newline  
(iii) Otherwise we put $\Gamma^{\sigma(i)}=\Gamma$
and $\Epipsi=\Epip$.\newline\newline
{\bf {2.4.}}  For $\Epip \in D(\Theta_N)$, let  
$({\cal E^\dagger},$ $\Pi^\dagger=\{\al_0^\dagger
,\ldots,\al_n^\dagger\},$ $p^\dagger)$ $\in D(\Theta_N)$
be an another datum satisfying:
\newline\newline
(i) ${\cal E^\dagger}={\cal E},$\newline\newline
(ii) For $\Epip \in D(\Theta_N)$, the
type of the Dynkin diagram $\Gamma^\dagger$ of
 $\Epipdagger \in D(\Theta_N)$ is: \[ \left\{
 \begin{array}{rl}
\mbox{$(AA)$}  &\quad\mbox{if $\Gamma$ is type $(AA)$,}\\
\mbox{$(BB)$}  &\quad\mbox{if $\Gamma$ is type $(BB)$,}\\
\mbox{$(CB)$}  &\quad\mbox{if $\Gamma$ is type $(CB)$ or $(DB)$,}\\
\mbox{$(CC)$}  &\quad\mbox{if $\Gamma$ is type $(CC)$,
 $(CD)$, $(DC)$ or $(DD)$.}\\
 \end{array}\right.\] 
We shall not need $p^\dagger$. So we merely denote $\Epipdagger$ by
$\Epidagger$.
\newline\newline 
{\bf {Remark 2.4.1.}} We can easily see that two $\Gamma^\dagger$'s 
defined for $\Epip$ and $\Epipsi$ are same.  
\newline\newline
{\bf Lemma 2.4.1.} 
$P_+
=Z_+\al_0 \oplus \ldots \oplus Z_+\al_n$
$\subset$
$P_+^\dagger
=Z_+\al_0^\dagger \oplus \ldots \oplus Z_+\al_n^\dagger$.
\newline\newline
{\bf {2.5.}}  On ${\cal E}=\oplus_{i=1}^N C\be_i\oplus C\dl
\oplus C\Lambda_0$, we define another symmetric form 
$((\,,\,)):{\cal E} \times {\cal E} \rightarrow C$ by
$$
((\be_i,\be_i))=\dl_{ij}, ((\be_i,\dl))=0, ((\dl,\dl))=0, 
((\dl,\Lambda_0))=1, ((\Lambda_0,\Lambda_0))=1.
$$
For $\Epip \in D(\Theta_N)$ and $0 \le i \le n$, define 
$I^{\sigma(i)}:{\cal E}=\oplus_{i=1}^N C\be_i\oplus C\dl
\oplus C\Lambda_0 \rightarrow 
{\cal E}^{\sigma(i)}=\oplus_{i=1}^N C\be_i'\oplus C\dl
\oplus C\Lambda_0$ by $I^{\sigma(i)}(\be_i)=\be_i'$,
$I^{\sigma(i)}(\dl)=\dl$ and 
$I^{\sigma(i)}(\Lambda_0)=\Lambda_0$.
We also define a linear map 
$\sigma(i):{\cal E}\rightarrow{\cal E}^{\sigma(i)}$ by
$$
\sigma(i)\,(v)=I^{\sigma(i)}(v-{\frac {2((v,\al_i^\dagger))} 
{((\al_i^\dagger,\al_i^\dagger))}}\al_i^\dagger).
$$
As an easy consequence from Propositions 2.2.1-3, we have:
\newline\newline
{\bf Theorem 2.5.1.} {\it For $\Epip \in D(\Theta_N)$ and $0 \le i \le n$,
Let ${\cal H}\oplus \oplus_{\al\in {\Phi}^{\sigma(i)}}
 \g_{\al}^{\sigma(i)}$
be the root space decomposition of $\g\Epipsi$.\par
(i) There is an isomorphism $L_i:\g\Epip\rightarrow\g\Epipsi$
such that:
$$
 L_i(H_{\gamma})=H_{\sigma(i)(\gamma)} \, (\gamma \in {\cal E}).  \eqno{(2.5.1)}
$$ In particular, $L_i$ satisfies:
$$
L_i(\g_{\al})=\g^{\sigma(i)}_{\sigma(i)(\al)} 
\quad (\al \in \Phi). 
  \eqno{(2.5.2)}
$$ \par
(ii) Let $L'_i$ and $L_i$ be two isomorphisms satisfying (2.5.1)
Then there exist $a_j \in C^*$ such that
$$
L'_i(E_j)=a_jL_i(E_j), \,\, L'_i(F_j)=a_j^{-1}L_i(F_j). \eqno{(2.5.3)} 
$$}  
{\bf Proof.} (i) We can choose one of $\phi$'s in Propositions 2.2.1-3
as $L_i$. It is obvious $L_i$ satisfies (2.5.1-2). \par
(ii) By (2.5.2) and the fact of $\dim\g_{\al_j}=1$,  
$\dim\g^{\sigma(i)}_{\sigma(i)(\al_j)}=1$. By the fact of
$[L_i(E_j),L_i(F_j)]$ $=H_{\sigma(i)(\al_j)}$
 $=[L'_i(E_j),L'_i(F_j)]$, we get (2.5.3).
\newline
\newline
{\bf Remark 2.5.1.} For example, $\sigma(i)$ and $L_i$
move as follows:

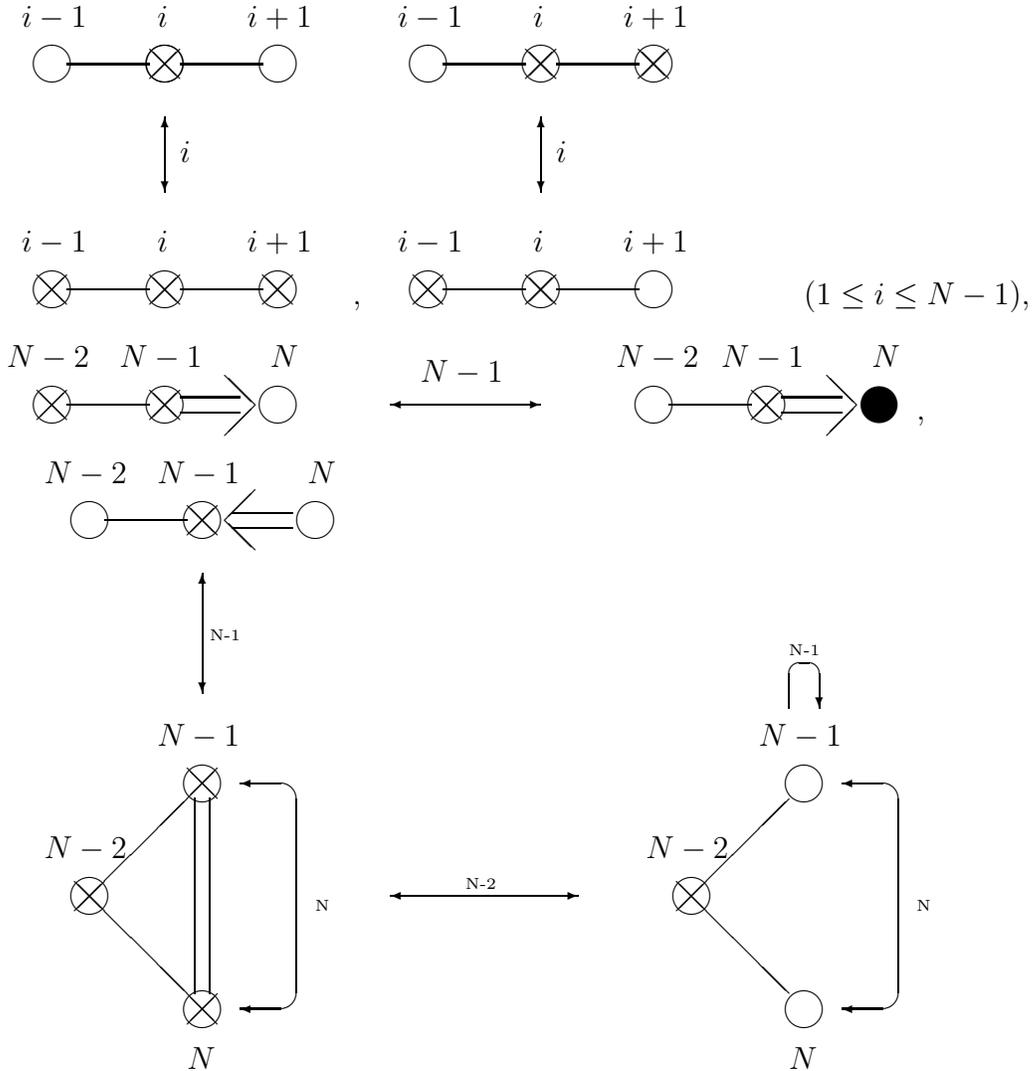
\setlength{\unitlength}{1mm}
\begin{picture}(180,45)(0,0)

                           \put(1, 40){$i-1$}  

\put(5, 35){\circle{5}}
\put(7, 35){\line(1,0){11}}

\put(18, 33){\line(1,1){4}} \put(19, 40){$i$} 
\put(18, 37){\line(1,-1){4}}
\put(20, 35){\circle{5}}
\put(22, 35){\line(1,0){11}}

\put(20, 35){\circle{5}}
\put(22, 35){\line(1,0){11}}

                           \put(31, 40){$i+1$} 

\put(35, 35){\circle{5}}

\put(20, 28){\vector(0,-1){10}} \put(22, 22){$i$}
\put(20, 18){\vector(0,1){10}}

\put(3, 3){\line(1,1){4}} \put(1, 10){$i-1$}  
\put(3, 7){\line(1,-1){4}}
\put(5, 5){\circle{5}}
\put(7, 5){\line(1,0){11}}

\put(18, 3){\line(1,1){4}} \put(19, 10){$i$} 
\put(18, 7){\line(1,-1){4}}
\put(20, 5){\circle{5}}
\put(22, 5){\line(1,0){11}}

\put(33, 3){\line(1,1){4}} \put(31, 10){$i+1$} 
\put(33, 7){\line(1,-1){4}}
\put(35, 5){\circle{5}}

\put(45, 3){$,$}

                           \put(51, 40){$i-1$}  

\put(55, 35){\circle{5}}
\put(57, 35){\line(1,0){11}}

\put(68, 33){\line(1,1){4}} \put(69, 40){$i$} 
\put(68, 37){\line(1,-1){4}}
\put(70, 35){\circle{5}}
\put(72, 35){\line(1,0){11}}

\put(83, 33){\line(1,1){4}} \put(81, 40){$i+1$} 
\put(83, 37){\line(1,-1){4}}
\put(85, 35){\circle{5}}

\put(70, 28){\vector(0,-1){10}} \put(72, 22){$i$}
\put(70, 18){\vector(0,1){10}}

\put(53, 3){\line(1,1){4}} \put(51, 10){$i-1$}  
\put(53, 7){\line(1,-1){4}}
\put(55, 5){\circle{5}}
\put(57, 5){\line(1,0){11}}

\put(68, 3){\line(1,1){4}} \put(69, 10){$i$} 
\put(68, 7){\line(1,-1){4}}
\put(70, 5){\circle{5}}
\put(72, 5){\line(1,0){11}}

                           \put(81, 10){$i+1$} 
\put(85, 5){\circle{5}}

\put(105, 3){$(1 \le i \le N-1)$,}

\end{picture}

\setlength{\unitlength}{1mm}
\begin{picture}(180,15)(0,0)

\put(3, 3){\line(1,1){4}} \put(-1, 10){$N-2$}  
\put(3, 7){\line(1,-1){4}}
\put(5, 5){\circle{5}}
\put(7, 5){\line(1,0){11}}

\put(18, 3){\line(1,1){4}} \put(14, 10){$N-1$} 
\put(18, 7){\line(1,-1){4}}
\put(20, 5){\circle{5}}
\put(22, 4){\line(1,0){8}}
\put(22, 6){\line(1,0){8}}
\put(32, 5){\line(-1,-1){4}}
\put(32, 5){\line(-1,1){4}}

																									 \put(34, 10){$N$} 
\put(35, 5){\circle{5}}

\put(50, 5){\vector(1,0){20}} \put(54, 8){$N-1$}
\put(70, 5){\vector(-1,0){20}}

																													 \put(80, 10){$N-2$}  
\put(85, 5){\circle{5}}
\put(87, 5){\line(1,0){11}}

\put(98, 3){\line(1,1){4}} \put(94, 10){$N-1$} 
\put(98, 7){\line(1,-1){4}}
\put(100, 5){\circle{5}}
\put(102, 4){\line(1,0){8}}
\put(102, 6){\line(1,0){8}}
\put(112, 5){\line(-1,-1){4}}
\put(112, 5){\line(-1,1){4}}

																										 \put(114, 10){$N$} 
\put(115, 5){\circle*{5}}
\put(120, 3){$,$}

\end{picture}

\setlength{\unitlength}{1mm}
\begin{picture}(180,85)(0,0)

\put(8, 23){\line(1,1){4}} \put(4, 30){$N-2$}  
\put(8, 27){\line(1,-1){4}}
\put(10, 25){\circle{5}}
\put(12, 27){\line(1,1){11}}
\put(12, 23){\line(1,-1){11}}

\put(23, 38){\line(1,1){4}} \put(19, 45){$N-1$} 
\put(23, 42){\line(1,-1){4}}
\put(25, 40){\circle{5}}

\put(24, 38){\line(0,-1){26}}
\put(26, 38){\line(0,-1){26}}

\put(23,  8){\line(1,1){4}} \put(23, 2){$N$} 
\put(23, 12){\line(1,-1){4}}
\put(25, 10){\circle{5}}

\put(35, 25){\oval(5,30)[r]}
\put(35, 40){\vector(-1,0){5}}
\put(35, 10){\vector(-1,0){5}}
\put(40, 23){{\tiny N}}


\put(50, 25){\vector(1,0){25}}\put(60, 26){{\tiny {N-2}}}
\put(75, 25){\vector(-1,0){25}}


\put(88, 23){\line(1,1){4}} \put(84, 30){$N-2$}  
\put(88, 27){\line(1,-1){4}}
\put(90, 25){\circle{5}}
\put(92, 27){\line(1,1){11}}
\put(92, 23){\line(1,-1){11}}

                            \put(99, 45){$N-1$} 

\put(105, 40){\circle{5}}

                            \put(103, 2){$N$} 

\put(105, 10){\circle{5}}

\put(115, 25){\oval(5,30)[r]}
\put(115, 40){\vector(-1,0){5}}
\put(115, 10){\vector(-1,0){5}}
\put(120, 23){{\tiny N}}

\put(105, 50){\oval(4,12)[t]}
\put(107, 54){\vector(0,-1){4}}
\put(103, 57){{\tiny {N-1}}}


\put(25, 52){\vector(0,1){16}} \put(26, 59){{\tiny {N-1}}}
\put(25, 68){\vector(0,-1){16}}


                            \put( 4, 80){$N-2$}  
\put(10, 75){\circle{5}}
\put(12, 75){\line(1,0){11}}

\put(23, 73){\line(1,1){4}} \put(19, 80){$N-1$} 
\put(23, 77){\line(1,-1){4}}
\put(25, 75){\circle{5}}
\put(28, 75){\line(1,1){4}}
\put(28, 75){\line(1,-1){4}}
\put(29, 74){\line(1,0){8}}
\put(29, 76){\line(1,0){8}}

																									 \put(39, 80){$N$} 
\put(40, 75){\circle{5}}

\end{picture}

{\bf Proposition 2.5.1.} For $\Epip \in D(\Theta_N)$, let
$\g=\g\Epip$. Then $\dim \g_\al=1$
if $\al \in \Phi\setminus Z\dl$.\newline
{\bf {Proof.}} We consider the case of 
$\al \in \Phi_+\setminus Z_+\dl$. 
We use an induction on
the height with respect to $\Pi^\dagger$. By 
$((\al,\al)) > 0$, $((\al,\al_i))> 0$ for some $i$.
If $\al\notin Z\al_i$ then the height of
 $\sigma(i)(\al)$ is smaller than the height
of $\al$. Thus by $\sigma(i)$'s, we can get a path from
$\al$ to $r\al_i\in\Pi'$ 
of another datum $\Epipdash$
where $r\in\{1,2\}$ if $p(\al_i)=1$ and 
$(\al_i,\al_i)\ne 0$, $r=1$ otherwise. 
On the other hand, 
$\dim \g'_{r\al_j}=1$ by Lemma 1.2.1.
\par\hfill Q.E.D. \newline\newline

{\bf {2.6.}}
Let $D(\Theta_N)_{(XY)}$ denote 
$\{\Epip \in \Theta_N | \Epip\,\,is\,\,(XY)-type\}$.
Then $\{\sigma(i)|1\le i\le n\}$ preserve 
$D(\Theta_N)_1$ $=D(\Theta_N)_{(AA)}$, 
$D(\Theta_N)_2$ $=D(\Theta_N)_{(BB)}$, 
$D(\Theta_N)_3$ $=D(\Theta_N)_{(CB)}\cup D(\Theta_N)_{(DB)}$ or
$D(\Theta_N)_4$ $=D(\Theta_N)_{(CC)}$
$\cup D(\Theta_N)_{(CD)}$
$\cup D(\Theta_N)_{(DC)}$
$\cup D(\Theta_N)_{(DD)}$.
Let $W=W(i)$ denote a group generated by $\{\sigma(i)|1\le i\le n\}$
acting on $D(\Theta_N)_i$. Then
$W(i)$ is isomorphic to the affine Weyl group of the corresponding $\Epidagger$.\par Let $D(\Theta_N)_i= \cup_{j=1}^{a_i}D(\Theta_N)_{ij}$
be the orbit decomposition. If a datum $\Epip$ associated to
a Dynkin diagram
$\Gamma \in \Theta_N$ belongs to $D(\Theta_N)_{ij}$,
then we denote $D(\Theta_N)_{ij}$ by $D(\Theta_N)[\Gamma]$.
\newline\newline
{\bf {2.7.}}
{\bf Proposition 2.7.1.} {\it Fix an orbit $D(\Theta_N)[\Gamma]$. Fix an isomorphism $\L_i:\g\Epip\rightarrow\g\Epipsi$ for each $\Epip$ and $\sigma(i)$.
Then there exists a unique family  $\{\bg\Epip$
 $|\Epip\in D(\Theta_N)[\Gamma]\}$  
of Lie superalgebras
satisfying following (1), (2) and (3).
\newline\newline
(1) For $\Epip\in D(\Theta_N)[\Gamma]$,
there is a sequence of  epimorphisms 
\[{\tilde \g}\Epip \stackrel{{\scriptscriptstyle {\tilde \Psi}\Epip}}
{\longrightarrow\!\!\!\!\rightarrow}
\bg\Epip
\stackrel{{\scriptscriptstyle \Psi\Epip}}
{\longrightarrow\!\!\!\!\rightarrow}
\g\Epip\]
\[(\,H_\gamma, E_i, F_i \stackrel{}\,\,{\longrightarrow}\,\,
H_\gamma, E_i, F_i
\stackrel{}\,\,{\longrightarrow}\,\,H_\gamma, E_i, F_i).\]
$$
$$
(2) For the isomorphism 
$L_i:\g\Epip\rightarrow \g\Epipsi$,
there is an isomorphism 
${\bar L}_i:\bg\Epip\rightarrow \bg\Epipsi$
satisfying the following commuting diagram:

\setlength{\unitlength}{1mm}
\begin{picture}(100,55)(0,0)

\put(15, 45){$\bg\Epip$} 
\put(40, 45){\vector(1,0){20}}\put(50, 47){${\bar L}_i$} 
\put(65, 45){$\bg\Epipsi$}

\put(25, 40){\vector(0,-1){20}}\put(13, 30){${\scriptscriptstyle \Psi\Epip}$}
\put(75, 40){\vector(0,-1){20}}\put(78, 30){${\scriptscriptstyle \Psi\Epipsi}$} 

\put(15, 15){$\g\Epip$} 
\put(40, 15){\vector(1,0){20}}\put(50, 17){$L_i$}  
\put(65, 15){$\g\Epipsi$}

\put(50, 32){\oval(6,6)[t]}
\put(50, 32){\oval(6,6)[br]}
\put(50, 29){\vector(-1,0){2}}

\end{picture}

$$
$$
(3) If  there is a family $\{ 
\bg^{(\lambda)}\Epip$ $|\Epip\in D(\Theta_N)[\Gamma]\}$  
of Lie superalgebras satisfying (1) and (2),
 then there are epimorphisms 
$$
{\bar \Psi}^{(\lambda)}\Epip
: \bg\Epip
\rightarrow\!\!\!\!\rightarrow
{\tilde \g}^{(\lambda)}\Epip \quad\quad 
(\Epip \in D(\Theta_N)[\Gamma])
$$ satisfying
following commutative diagram
\begin{center}
Diagram 2.7.1.
\setlength{\unitlength}{1mm}
\begin{picture}(180,85)(0,0)

\put(15, 75){$\bg\Epip$} 
\put(40, 75){\vector(1,0){40}}\put(60, 77){${\bar L}_i$} 
\put(85, 75){$\bg\Epipsi$}

\put(25, 70){\vector(0,-1){60}}\put(13, 60){${\scriptscriptstyle \Psi\Epip}$} 
\put(25, 70){\vector(0,-1){58}}
\put(95, 70){\vector(0,-1){60}}\put(70, 60){${\scriptscriptstyle \Psi\Epipsi}$}  
\put(95, 70){\vector(0,-1){58}} 

\put(15,  5){$\g\Epip$} 
\put(40,  5){\vector(1,0){40}}\put(60, 7){$L_i$}
\put(85,  5){$\g\Epipsi$}

\put(32, 71){\vector(1,-3){8}}
\put(32, 71){\vector(1,-3){7.5}}
\put(36, 62){${\scriptscriptstyle {\bar \Psi}^{(\lambda)}\Epip}$} 
\put(35, 40){$\bg^{(\lambda)}\Epip$} 
\put(40, 35){\vector(-1,-3){8}}
\put(40, 35){\vector(-1,-3){7.5}}
\put(38, 22){${\scriptscriptstyle \Psi^{(\lambda)}\Epip}$}
\put(62, 40){\vector(1,0){40}} \put(80, 42){${\bar L}^{(\lambda)}_i$}

\put(102, 71){\vector(1,-3){8}}
\put(102, 71){\vector(1,-3){7.5}}
\put(106, 62){${\scriptscriptstyle {\bar \Psi}^{(\lambda)}\Epipsi}$}
\put(105, 40){$\bg^{(\lambda)}\Epipsi$} 
\put(110, 35){\vector(-1,-3){8}}
\put(110, 35){\vector(-1,-3){7.5}}
\put(108, 22){${\scriptscriptstyle \Psi^{(\lambda)}\Epipsi}$}

\end{picture}
\end{center}
(ii) The set $\{\bg\Epip$ $|\Epip\in D(\Theta_N)[\Gamma]\}$ 
does not depend on the choice of $\{L_i\}$. } \newline\newline
{\bf Proof.} (i) Let $\{\,C^{(\lambda)}= 
\{\bg^{(\lambda)}\Epip|\Epip\in D(\Theta_N)[\Gamma]\}\,\}
_{(\lambda)\in(\Lambda)}$ be the family of the families 
$C^{(\lambda)}$ $((\lambda)\in(\Lambda))$ satisfying (1) and (2).
Let 
$
{\tilde \Psi}^{(\lambda)}\Epip:{\tilde \g}\Epip$ 
$\rightarrow$
$\bg^{(\lambda)}\Epip$
be the epimorphism in (1) for $(\lambda)\in(\Lambda)$.
Let ${\bar r}^{(\lambda)}\Epip=\ker{\tilde \Psi}^{(\lambda)}\Epip$ and
${\bar r}\Epip=\cap_{(\lambda)\in(\Lambda)}{\bar r}^{(\lambda)}\Epip$.
Put $\bg\Epip = {\tilde \g}\Epip/{\bar r}\Epip$.
Let $[L_i(E_j)]$, $[L_i(F_j)]$
$\in {\bg}\Epipsi_{\sigma(i)(\alpha_j)}$
be representatives of
$L_i(E_j)$, $L_i(F_j)$.
Similarly to Proposition 2.5.1, we can show 
$\dim\bg^{(\lambda)}\Epip_\al=1$ if 
$\al \in \Phi\setminus Z\dl$. 
Then ${\bar r}^{(\lambda)}\Epipsi_{\sigma(i)(\alpha_j)}$ $=$
$r\Epipsi_{\sigma(i)(\alpha_j)}$.
Hence $[L_i(E_j)] \equiv {\bar L_i}^{(\lambda)}(E_j)$, 
$[L_i(F_j)] \equiv {\bar L_i}^{(\lambda)}(F_j)$
(mod ${\bar r}^{(\lambda)}\Epipsi$). \newline Hence,
in $\bg^{(\lambda)}\Epipsi$, the elements 
$H_{\sigma(i)(\gamma)}$, $[L_i(E_j)]$ and $[L_i(F_j)]$ satisfy
(1.2.1-3) for $\Epip$ whence, even in $\bg\Epip$, the elements satisfy (1.2.1-3). Hence there is a morphism
${\tilde {\bar L_i}}:{\tilde \g}\Epip
\rightarrow \bg\Epipsi$ such that 
$${\bar \Psi}^{(\lambda)}\Epipsi\circ{\tilde {\bar L_i}}
 ={\bar L_i}^{(\lambda)} \circ {\tilde \Psi}^{(\lambda)}\Epip. \eqno{(2.7.1)}$$
Therefore there exists 
${\bar L_i}:\bg\Epip
\rightarrow \bg\Epipsi$ 
such that ${\tilde {\bar L_i}}={\bar L_i}\circ{\tilde \Psi}\Epip$.
Clearly ${\bar L_i}$ satisfying the commutative diagram of (2).
Using Lemma 2.2.1, we can see that 
${\bar L_i}\circ{\bar L_i}:$ 
$\bg\Epip\rightarrow \bg\Epip$
satisfies that ${\bar L_i}\circ{\bar L_i}(H_\gamma)=(H_\gamma)$,
and ${\bar L_i}\circ{\bar L_i}(E_j)$, ${\bar L_i}\circ{\bar L_i}(F_j)$
are nonzero scalar multiples of $E_j$, $F_j$ respectively. Hence
${\bar L_i}$ is an isomorphism. \newline
(ii) By theorem 2.5.1,
if $L'_i$ is another $L_i$, then $L'_i(E_j)$ 
$=a_jL_i(E_j)$, $L'_i(F_j)$ $=a_j^{-1}L_i(F_j)$ 
for some $a_j\in C\setminus\{0\}$ 
$(0\leq j \leq n)$. Let $\phi_a : \g\Epip$
$\rightarrow\g\Epip$ be an isomorphism defined by
$\phi_a(E_j)=a_jE_j$, $\phi_a(F_j)=a_j^{-1}F_j$, 
$\phi_a(H)=H$. Then $L'_i=L_i\circ\phi_a$.
The ideal ${\bar r}\Epip$ in the proof of (i)
satisfies ${\cal H} \cap {\bar r}\Epip =0$. Then 
 ${\bar r}\Epip$ is the homogeneous ideal.
Then we can also define an isomorphism ${\bar \phi}_a : {\bar\g}\Epip$
$\rightarrow{\bar\g}\Epip$ similar to $\phi_a$. 
Denote ${\bar L}_i$ defined for $L'_i$ by ${\bar L}'_i$. 
By the universality of ${\bar L}_i$, it follows that 
${\bar L}'_i={\bar L}_i\circ{\bar \phi}_a$. In particular,
${\bar\g}\Epip$ is determined independently of the choice
of $L_i$.
\par \hfill Q.E.D. \newline\newline
{\bf Lemma 2.7.1.} {\it Let $\Epip \in D(\Theta_N)$
and $\bg=\bg\Epip$.
Let ${\bar {\cal N}}^+$, ${\bar {\cal N}}^-$
be the subalgebras of $\bg$
generated by $E_i$, $F_i$ respectively.
Then we have the triangular decomposition
$$
\bg=
{\bar {\cal N}}^+\oplus {\cal H}\oplus {\bar {\cal N}}^-.
$$ Here ${\cal H}$ can be identified with
${\cal H}$ of $\g$. We have 
${\bar {\cal N}}^+ \cong {\bar {\cal N}}^-$
($E_i \leftrightarrow F_i$).}\newline
{\bf Proof.} The triangular decomposition is clear because ${\bar r}\cap{\cal H}=0$.
Let ${\bar r}$ be the ideal ${\bar r}\Epip$ of 
${\widetilde \g}\Epip$
in the proof of Proposition 2.7.1.
Let ${\bar r}_{\pm}={\bar r}\cap {\widetilde {\cal N}}^{\pm}$.
Let ${\bar r}^1_{-}$ (resp. ${\bar r}^1_{+}$)
be the ideal defined as the image of ${\bar r}_{+}$
of the map 
${\widetilde {\cal N}}^{+}\rightarrow {\widetilde {\cal N}}^{-}$ 
($E_i\rightarrow F_i$)
(resp.
${\widetilde {\cal N}}^{-}\rightarrow {\widetilde {\cal N}}^{+}$ 
($F_i\rightarrow E_i$) \,). Put 
${\bar r}^1={\bar r}^1_{-}\oplus{\bar r}^1_{-}$
and $\bg^1={\widetilde \g}/{\bar r}^1$.
By the universality, we can show 
$\bg^1=\bg$. Then
we have ${\bar r}^1_{\mp}={\bar r}_{\pm}$
\par \hfill Q.E.D. \newline\newline
{\bf Lemma 2.7.2.} {\it For $\al_i\in\Pi$, we have
$\dim\bg_{\al_i}=1$. If $p(\al_i)=1$
and $(\al_i,\al_i)\ne 0$, then 
 $\dim\g_{2\al_i}=1$.} \newline
{\bf Proof.} The proof is obtained similarly to the proof
of Lemma 1.2.1.
\par \hfill Q.E.D. \newline\newline
{\bf Proposition 2.7.2.} {\it Let $\Phi$ be the set of the roots
of $\g\Epip$ associated with $\Epip \in D(\Theta_N)$. 
For $\bg=\bg\Epip$, let 
$\bg_\gamma = \{x\in\bg\,|\,
[H,x]=\gamma(H)x\}$. Then we have \newline\newline
(i) $\dim\bg_\gamma=1$ 
if $\gamma\in\Phi\setminus Z\dl$. \newline
(ii) $\dim\bg_\gamma\ge\dim\g_\gamma$
if $\gamma\in Z\dl$. \newline
(iii)  $\dim\bg_\gamma=0$
if $\gamma\notin\Phi\cup\{0\}$. \newline\newline
In particular, 
$$
\ker \Psi\Epip\subset 
\oplus_{r\ne 0}\bg_{r\dl}.
$$}
{\bf Proof.} 
By Lemma 2.7.1, we may assume $\gamma \in P_+$. (ii) is clear.
Using Lemma 2.7.2, we can proof (i) similarly to the proof
of Proposition 2.5.1. Moreover (iii) can be proved similarly to the proof
of Proposition 2.5.1: If $\gamma\notin\Phi\cup\{0\}$, then,
by $\sigma(i)$'s, we can get a path from
$\gamma$ to $(\oplus_{\al_i\in\Pi'}Z\al_i)\setminus
(P'_+ \cup -P'_+ )$ 
of another datum $\Epipdash$.
By the triangular decomposition of
$\bg\Epipdash$, we have 
$\dim\bg_\gamma=0$.
\par \hfill Q.E.D. \newline\newline
\hspace{1cm}\newline
{\bf {2.8.}} {\bf Proposition 2.8.1.} 
 {\it Let $\Epip\in D(\Theta_N)$. \par
(i) There exists a unique Lie superalgebra 
 $\bg^{\ddag}\Epip$
satisfying following (1), (2) and (3). \par
(1) For $\Epip\in D(\Theta_N)$,
there is a sequence of  epimorphisms 
\[{\tilde \g}\Epip 
\stackrel{{\scriptscriptstyle {\tilde \Psi}^{\ddag}\Epip}}
{\longrightarrow\!\!\!\!\rightarrow}
\bg^{\ddag}\Epip
\stackrel{{\scriptscriptstyle \Psi^{\ddag}\Epip}}
{\longrightarrow\!\!\!\!\rightarrow}
\g\Epip\]
\[(\,H_\gamma, E_i, F_i \stackrel{}\,\,{\longrightarrow}\,\,
H_\gamma, E_i, F_i
\stackrel{}\,\,{\longrightarrow}\,\,H_\gamma, E_i, F_i).\] \par
(2) By (1), there is a root space decomposition
\[
  \bg^{\ddag}\Epip
= {\cal H}\oplus(\oplus_{\al\in \Psi^{\ddag}}\bg^{\ddag}_{\al})  
\] such
that ${\cal E}={\cal H}^{*}\supset\Psi^{\ddag}\supset\Psi$. Here 
$\Psi$ is the set of the roots of $\g\Epip$. Then the assumption 
(2) is that
$\Psi^{\ddag}=\Psi$ and $\dim\bg^{\ddag}_{\al}=1$
if $\al \in \Psi \setminus Z\dl$. \par
(3) If  there is a Lie superalgebra
 $\bg^{\ddag(\lambda)}\Epip$ satisfying (1) and (2),
 then there is an epimorphism: 
$$
{\bar \Psi}^{\ddag(\lambda)}\Epip
: \bg^{\ddag}\Epip
\rightarrow\!\!\!\!\rightarrow
\bg^{\ddag(\lambda)}\Epip
$$ \par
(ii) $\bg^{\ddag}\Epip$
is isomorphic to
$\bg\Epip$.
}\newline\newline
{\bf {Proof.}} (i) 
Let $\{{\bg}^{\ddag(\lambda)}\Epip\}
_{(\lambda)\in(\Lambda)}$ be the family of 
Lie superalgebras
 satisfying (1) and (2).
Let 
${\tilde \Psi}^{\ddag(\lambda)}\Epip:$ ${\tilde \g}\Epip$ 
$\rightarrow$
$\bg^{\ddag(\lambda)}\Epip$
be the epimorphism in (1) for $(\lambda)\in(\Lambda)$.
Let ${\bar r}^{\ddag(\lambda)}\Epip
=\ker{\tilde \Psi}^{\ddag(\lambda)}\Epip$ and
${\bar r}^{\ddag}\Epip
=\cap_{(\lambda)\in(\Lambda)}{\bar r}^{\ddag(\lambda)}\Epip$.
Put $\bg^{\ddag}\Epip 
= {\tilde \g}\Epip/{\bar r}^{\ddag}\Epip$.
Then $\bg^{\ddag}\Epip$ satisfies (1), (2) and (3).\par
(ii) By the universality of $\bg^{\ddag}\Epip$,
there is an epimorphism $L_i^{\ddag}:$ 
$\bg^{\ddag}\Epip$ $\rightarrow$
$\bg^{\ddag}\Epipsi$ such that
$L_i^{\ddag}(H_{\gamma})=H_{\sigma(i)(\gamma)}$,
$L_i^{\ddag}(E_j)={\Psi}^{\ddag}\Epip^{-1}L_i(E_j)$.
$L_i^{\ddag}(E_j)={\Psi}^{\ddag}\Epip^{-1}L_i(E_j)$.
Clearly $\{L_i\}$ satisfy (2) of Proposition 2.7.1.
Then we have an epimorphism $\bg\Epip$ 
$\rightarrow$
$\bg^{\ddag}\Epip$. By the universality of 
$\bg^{\ddag}\Epip$,
this map is isomorphism. 
\par \hfill Q.E.D. \newline\newline\hspace{1cm}\par
As an immediate consequence, we have:\newline\newline
{\bf Lemma 2.8.1.} {\it For $\Epip\in D(\Theta_N)$,
there is an epimorphism
$$
{\Psi}^{\dag}\Epip
: \bg\Epip
\rightarrow\!\!\!\!\rightarrow
\bbg\Epip.
$$ where $\bbg\Epip$ has already introduced in \S 1 for $\Epip$
of affine $ABCD$-type.} 
\newline\par
In $\S$3, we will show that ${\Psi}^{\dag}\Epip$ is isomorphism.
$$
$$\newline
{\bf {3. The estimation of $\dim\bg_{r\dl}$}}\newline 
{\bf {3.1.}} {\bf Proposition 3.1.1.} {\it (i) Let $\Epip\in D(\Theta_N)$. 
Put $D_N=\sum_{i=1}^N\bd_i$. Then
\[ (\dl,2\rho)= \left\{
 \begin{array}{ll}
2D_N  &\quad\mbox{$(AA)$,}\\
2D_N  &\quad\mbox{$(BB)$,}\\
2\bd_1+4D_N  &\quad\mbox{$(CB)$,} \\
-2\bd_1+4D_N  &\quad\mbox{$(DB)$,} \\
2\bd_1+4D_N+2\bd_N  &\quad\mbox{$(CC)$,} \\
-2\bd_1+4D_N+2\bd_N  &\quad\mbox{$(DC)$,} \\
2\bd_1+4D_N-2\bd_N  &\quad\mbox{$(CD)$,} \\
-2\bd_1+4D_N-2\bd_N  &\quad\mbox{$(DD)$.} \\
 \end{array}\right.\]
(ii) For exceptional type, we have: 
\[ (\dl,2\rho)=\,\, \left\{
 \begin{array}{ll}
 0  &\quad\mbox{$D(2,1,;x)^{(1)}$,}\\
 -12  &\quad\mbox{$F_4^{(1)}$,}\\
 12 &\quad\mbox{$G_3^{(1)}$.} \\
 \end{array}\right.\]}
{\bf Proof.} Direct calculations. 
\par \hfill Q.E.D. \newline\newline

Here we again remark, if there is a relation of weight 
$\beta\in P_+$ such that $(\beta,\beta)\ne 2(\beta, \rho)$,
then the relation can be obtained by relations of lower weights
(see Proposition 1.2.1). Therefore have:\newline\newline
{\bf Lemma 3.1.1} If
$\Epip\in D(\Theta_N)$ satisfies
$(\dl,\dl)\ne 2(\dl, \rho)$, then $\bg\Epip
\cong\g\Epip$. 
\newline\newline
{\bf {3.2.}} {\bf Proposition 3.2.1.} {\it 
 Let $\Epip\in D(\Theta_N)[\Gamma]$. 
Assume that $x$-th simple root $\al_x\in\Pi$ and
$y$-th simple root $\al_y\in\Pi$ satisfy 
$\al_x+\al_y\in\Phi_+$ (i.e. $(\al_x,\al_y)\ne 0$),
$\al_x+\al_y=\sigma(x)(\al_y^\dagger)$ and
that $\Pi_{(x,y)}$ $=(\Pi\setminus\{\al_x,\al_y\})$ 
$\cup\{\al_x+\al_y\}$ is 
affine $ABCD$ type. 
\par (i) Let $\Pi_{(x,y)}^\dagger$ be $\Pi^\dagger$ of 
$\Pi_{(x,y)}$. Then $\Pi_{(x,y)}^\dagger$ 
$=(\Pi^\dagger\setminus\{\al^\dagger_x,\al^\dagger_y\})$ 
$\cup\sigma(x)(\al_y^\dagger)$
Let $W_{(x,y)}$ denote a subgroup of $W$
generated by 
$(\{\sigma(0),\ldots,\sigma(n)\}\setminus
\{\sigma(x),\sigma(y)\})$ $\cup\{\sigma(x)\sigma(y)\sigma(x)\}$.
Then $W_{(x,y)}$ is $W$ defined for $\Epipxy$. \par
(ii) Let $\Epipxy$ be a datum such that the set of
simple roots is $\Pi_{(x,y)}$.  
Then there is a homomorphism 
$i:\bg\Epipxy\rightarrow\bg\Epip$
such that 
\[ i(H_{\al_j})=\,\, \left\{
 \begin{array}{ll}
H_{\al_x+\al_y} &\quad\mbox{$\al_j=\al_x+\al_y$,}\\
H_{\al_j} &\quad\mbox{$\al_j\ne\al_x+\al_y$,}\\
 \end{array}\right.\]
\[ i(E_j)=\,\, \left\{
 \begin{array}{ll}
(-1)^{p(\al_x)p(\al_y)}
(\al_x,\al_y)^{-1}[E_x,E_y] &\quad\mbox{$\al_j=\al_x+\al_y$,}\\
E_j &\quad\mbox{$\al_j\ne\al_x+\al_y$,}\\
 \end{array}\right.\]
\[ i(F_j)=\,\, \left\{
 \begin{array}{ll}
[F_x,F_y] &\quad\mbox{$\al_j=\al_x+\al_y$,}\\
F_j &\quad\mbox{$\al_j\ne\al_x+\al_y$.}\\
 \end{array}\right.\]}
\newline\newline
{\bf Proof.} (i) We can check the fact for each affine type. \par 
(ii) Let $\{{\bar L}_0,\ldots, {\bar L}_n\}$ be the isomorphisms defined in 
Proposition 2.7.1.
We consider an orbit 
$Orbit(\bg\Epip)_{(x,y)}$ through $\bg\Epip$ under the 
action of a subgroup
generated by 
$(\{{\bar L}_0,\ldots, {\bar L}_n\}\setminus
\{{\bar L}_x,{\bar L}_y\})\cup\{{\bar L}_x{\bar L}_y{\bar L}_x\}$.
Then $Orbit(\bg\Epip)_{(x,y)}$ satisfies the conditions of 
$\{\bg^{(\lambda)}\Epipxy\}$ in Proposition 2.7.1
where $\{{\bar L}_i^{(\lambda)}\}$ is 
$(\{{\bar L}_0,\ldots, {\bar L}_n\}\setminus
\{{\bar L}_x,{\bar L}_y\})\cup\{{\bar L}_x{\bar L}_y{\bar L}_x\}$.
By the universality of $\bg\Epipxy$,
we get the homomorphism 
$i:\bg\Epipxy\rightarrow\bg\Epip$.
\par \hfill Q.E.D. \newline\newline

\hspace{1cm}\newline
{\bf {3.3.}} Let $\Epip\in D(\Theta_N)[\Gamma]$.
Here we shall show that $\bg\Epip$ is isomorphic to
$\bbg\Epip$. By Proposition 2.8.1, we have to show
$\dim\bbg\Epip_{n\dl}=$
$\dim\bg\Epip_{n\dl}$.
{\it Here we only prove the fact in the case of
$\Epip$ of Diagram 1.11.6 with $\sum p(\al_i)\equiv 1$.}
Because we can prove the fact in another case by the 
similar way. In this case, 
$\bbg={\widehat {osp}}({\tilde {\cal E}}_0, e)^{(I)}$
(see 1.11). Then we have to show:
$$
\dim\bg\Epip_{n\dl}\leq  N \quad (n \ne 0).\eqno{(3.3.1)} 
$$ We
start with $N=3$. In this case, its Dynkin diagram is:
\begin{center}
Diagram 3.3.1 
\newline
\DiagramCDd
\end{center}
The isomorphism $\sigma(1)$ divert Diagram 3.3.1 into: 
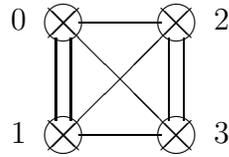
\begin{center}
Diagram 3.3.2 
\newline\newline
\DiagramDDd
\end{center} Therefore
it is sufficient to show (3.3.1) for $\Epip$ of Diagram 3.3.2.
By Proposition 3.2.1, we have the homomorphism
$$
  i:\bg\Epipnisan \rightarrow \bg\Epip.
$$ Here
the Dynkin diagram of $\Epipnisan$ is:
\begin{center}
Diagram 3.3.3 
\newline\newline
\DiagramDDdd
\end{center}
This is equivalent to Diagram 1.6.2
as the Dynkin diagram of the Kac-Moody Lie superalgebra. 
By Lemma 3.1.1,
$\bg\Epipnisan$
$=\g\Epipnisan$. Hence 
$$
\dim\bg\Epipnisan_{n\dl}=2\quad (n \ne 0). \eqno{(3.3.2)}
$$ We
assume $n > 0$. For a root $\gamma \notin Z\dl$ of
 $\g\Epip$ (resp. $\g\Epipnisan$), let $E_\gamma$
(resp. $E_\gamma^{(2,3)}$) denote a non-zero element of
 $\bg\Epip_{\gamma}$ (resp.
 $\bg\Epipnisan_{\gamma}$). By (3.3.2), 
$\{ [E_{n\dl-\al_0}^{(2,3)},E_{\al_0}^{(2,3)}]$,
$[E_{n\dl-\al_1}^{(2,3)},E_{\al_1}^{(2,3)}]$,
$[E_{n\dl-(\al_2+\al_3)}^{(2,3)},E_{\al_2+\al_3}^{(2,3)}]\}$
are linearly dependent. Transposing by $i$, we see that
$\{[E_{n\dl-\al_0},E_0]$,
$[E_{n\dl-\al_1},E_1]$,
$[E_{n\dl-(\al_2+\al_3)},[E_2,E_3]]\}$
are linearly dependent. Since
$[E_{n\dl-(\al_2+\al_3)},[E_2,E_3]]$
$=[[E_{n\dl-(\al_2+\al_3)},E_2],E_3]$
$-(-1)^{p(\al_2)p(\al_3)}[[E_{n\dl-(\al_2+\al_3)},E_3],E_2]$,
$\{[E_{n\dl-\al_0},E_0]$,
$[E_{n\dl-\al_1},E_1]$,
$[E_{n\dl-\al_2},E_2]$,
$[E_{n\dl-\al_3},E_3]\}$ are linearly dependent.
Hence $\dim\bg\Epip\leq 3$. Then we could show (3.3.1)
for $N=3$. \par
Next we show (3.3.1) for $N\geq 4$ by induction. By Proposition 3.2.1,
we can use the homomorphism
$$
  i:\bg\EpipNniNsan \rightarrow \bg\Epip.
$$ Using a similar argument to that in the case of $N=3$, we can show
(3.3.1). \par
Using a similar argument to that in the above case, we can get:
\newline\newline
{\bf Theorem 3.3.1.} {\it Let $\Epip \in D(\Theta_N)$. Then
$\bg\Epip$ is isomorphic to 
$\bbg\Epip$.}
\newline\newline\newline
{\bf 4. Relations of Affine $ABCD$-types}
\newline
{\bf 4.1.} Let $\Epip \in D(\Theta_N)$. Using the definition of 
$\bg\Epip$ given in 
Proposition 2.7.1, we can directly calculate defining relations
of $\bg\Epip$.
\newline\newline
{\bf {Theorem 4.1.1.}} {\it Let $\Epip \in D(\Theta_N)$ 
(i.e., $\Epip$ is affine $ABCD$ type).
The Lie superalgebra $\bg\Epip$ is defined 
by generators $H \in {\cal H}$, $E_i$, $F_i$
$(0 \leq i \leq n)$
with parities $p(H)=0$, $p(E_i)=p(E_i)=p(\al_i)$ and relations:
$$
\leqno{(1)}\quad [H, H']=0, \quad (H, H' \in {\cal H}) 
$$
$$
\leqno{(2)}\quad [H, E_i]=\al_i(H)E_i,\quad  [H, F_i]=-\al_i(H)F_i,  
$$
$$
\leqno{(3)}\quad [E_i, F_j]=\delta_{ij}H_{\al_i},  
$$\newline
(4)\quad Relations of $E_i$'s.
\newline\newline
(i)\quad $[E_i, E_j]=0$
\begin{flushright}
if $(\al_i,\al_j)=0$, $(i\ne j)$, \end{flushright}
\vspace{1cm}
(ii)\quad $[E_i, E_i]=0$
 
\begin{flushright}
\setlength{\unitlength}{1mm}
\begin{picture}(40,15)(0,0)
\put(23, 5){if}
\put(33, 3){\line(1,1){4}} \put(34, 10){$i$} 
\put(33, 7){\line(1,-1){4}}
\put(35, 5){\circle{5}}
\end{picture}
\end{flushright}

\vspace{1cm}
(iii)\quad $\cwhl E_i, \cwhl E_i, \ldots , \cwhl E_i, E_j
 \cwhr\!\ldots\!\cwhr\cwhr =0 $\,\, 
\newline ($E_i$ appears $1-2(\al_i,\al_j)/(\al_i,\al_i)$ times)\newline
\begin{flushright}
if $(\al_i,\al_i)\ne 0$ and 
$(-1)^{\{p(\al_i){{2(\al_j,\al_i)} \over  {(\al_i,\al_i)}}\}}=1$,
\end{flushright}\vspace{1cm}
(iv)\quad $[\cwhl \cwhl E_i, E_j \cwhr ,E_k \cwhr , E_j]=0$
 
\begin{flushright}
\setlength{\unitlength}{1mm}
\begin{picture}(90,15)(0,0)
\put(23, 5){if}
\put(33, 3){\line(1,1){4}} \put(34, 10){$i$} 
\put(33, 7){\line(1,-1){4}}
\put(37, 5){\line(1,0){16}} \put(42, 8){$-x$} 
\put(53, 3){\line(1,1){4}} \put(54, 10){$j$} 
\put(53, 7){\line(1,-1){4}}
\put(55, 5){\circle{5}}
\put(57, 5){\line(1,0){16}} \put(64, 8){$x$}
\put(73, 3){\line(1,1){4}} \put(74, 10){$k$} 
\put(73, 7){\line(1,-1){4}}
\put(83, 5){$(x\ne0)$,}
\end{picture}
\end{flushright}
\vspace{1cm}
(v)\quad $[\cwhl \cwhl E_i, E_j \cwhr ,\cwhl 
\cwhl E_i, E_j \cwhr ,E_k \cwhr \cwhr , E_j]=0$
 
\begin{flushright}
\setlength{\unitlength}{1mm}
\begin{picture}(90,15)(0,0)
\put(23, 5){if}
\put(33, 3){\line(1,1){4}} \put(34, 10){$i$} 
\put(33, 7){\line(1,-1){4}}
\put(35, 5){\circle{5}}
\put(37, 5){\line(1,0){16}} 
\put(53, 3){\line(1,1){4}} \put(54, 10){$j$} 
\put(53, 7){\line(1,-1){4}}
\put(55, 5){\circle{5}}
\put(58, 5){\line(1,1){4}}
\put(58, 5){\line(1,-1){4}}
\put(60, 4){\line(1,0){13}}
\put(60, 6){\line(1,0){13}}
\put(75, 5){\circle{5}} \put(74, 10){$k$} 
\put(80, 5){,}
\end{picture}
\end{flushright}
\vspace{1cm}
(vi)\quad $[\cwhl \cwhl \cwhl \cwhl \cwhl 
E_i, E_j \cwhr , E_k \cwhr , E_l \cwhr ,
E_k \cwhr , E_j \cwhr , E_k]=0$

\begin{flushright}
\setlength{\unitlength}{1mm}
\begin{picture}(105,15)(0,0)
\put(23, 5){if}
\put(33, 3){\line(1,1){4}} \put(34, 10){$i$} 
\put(33, 7){\line(1,-1){4}}
\put(37, 5){\line(1,0){16}} 
                           \put(54, 10){$j$} 
\put(55, 5){\circle{5}}
\put(57, 5){\line(1,0){16}}
\put(73, 3){\line(1,1){4}} \put(74, 10){$k$}
\put(73, 7){\line(1,-1){4}}
\put(75, 5){\circle{5}}  
\put(78, 5){\line(1,1){4}}
\put(78, 5){\line(1,-1){4}}
\put(80, 4){\line(1,0){13}}
\put(80, 6){\line(1,0){13}}
\put(95, 5){\circle{5}} \put(94, 10){$l$} 
\put(100, 5){,}
\end{picture}
\end{flushright}
\vspace{1cm}
(vii)\quad $(-1)^{p(\al_i)p(\al_k)}
(\al_i,\al_k)
\cwhl \cwhl E_i, E_j \cwhr , E_k \cwhr
=
(-1)^{p(\al_i)p(\al_j)}
(\al_i,\al_j)
\cwhl \cwhl E_i, E_k \cwhr , E_j \cwhr $

\begin{flushright}
\setlength{\unitlength}{1mm}
\begin{picture}(45,35)(0,0)

\put(3, 15){if}
\put(20, 15){.}		 \put(19, 20){$i$} 

\put(22, 17){\line(1,1){11}} \put(24, 25){$a$}
\put(22, 13){\line(1,-1){11}} \put(24, 3){$b$}
\put(35, 30){.}  	  	  \put(40, 29){$j$} 
  
\put(35, 28){\line(0,-1){26}} \put(38, 15){$c$}
\put(35, 0){.}   		  \put(40, 0){$k$}

\end{picture} \\
  $(a=(\al_i,\al_j), b=(\al_i,\al_k), c=(\al_j,\al_k),$ \\
            $abc\ne0, a+b+c=0,$ \\ 
            $p(\al_i)p(\al_j)+p(\al_i)p(\al_j)
+p(\al_i)p(\al_j)\equiv1),$ 
\end{flushright}
\vspace{1cm}
(viii)
$\cwhl \cwhl \cwhl  E_i, E_j \cwhr ,
 \cwhl  E_j, E_k \cwhr \cwhr ,
 \cwhl  E_j, E_l \cwhr \cwhr 
= \cwhl \cwhl \cwhl  E_i, E_j \cwhr ,
 \cwhl  E_j, E_l \cwhr \cwhr ,
 \cwhl  E_j, E_k \cwhr \cwhr$

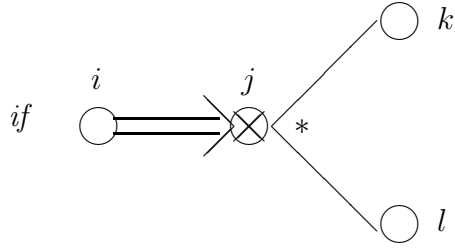
\begin{flushright}
\setlength{\unitlength}{1mm}
\begin{picture}(85,35)(0,0)

\put(23, 15){if}
\put(35, 15){\circle{5}} \put(34, 20){$i$}
\put(37, 14){\line(1,0){14}} 
\put(37, 16){\line(1,0){14}} 
\put(53, 15){\line(-1,1){4}}  
\put(53, 15){\line(-1,-1){4}}
\put(53, 13){\line(1,1){4}} \put(54, 20){$j$} 
\put(53, 17){\line(1,-1){4}}
\put(55, 15){\circle{5}} \put(61, 14){$*$}  
\put(58, 15){\line(1,1){14}} 
\put(58, 15){\line(1,-1){14}} 
                             \put(80, 28){$k$} 
\put(75, 29){\circle{5}}  
                             \put(80, 1){$l$}
\put(75, 2){\circle{5}}

\end{picture}
\end{flushright}
\vspace{1cm}
(ix)
\[ \begin{array}{l}
  [\cwhl \cwhl E_k, \cwhl E_l, \cwhl E_k, E_j, 
 \cwhr \cwhr \cwhr,
\cwhl E_k, \cwhl E_l, \cwhl E_k, \cwhl E_j, E_i, 
 \cwhr \cwhr \cwhr \cwhr \cwhr, 
E_j] \\
 = 2\cwhl \cwhl E_k, E_j \cwhr,
\cwhl \cwhl E_k, \cwhl E_j, E_i \cwhr \cwhr, 
\cwhl E_k, \cwhl E_l, \cwhl E_k, \cwhl E_j, E_i, 
 \cwhr \cwhr \cwhr \cwhr \cwhr \cwhr \\ 
\end{array}  \]

\begin{flushright}
\setlength{\unitlength}{1mm}
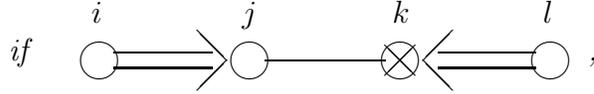
\begin{picture}(105,15)(0,0)
\put(23, 5){if}
\put(35, 5){\circle{5}} \put(34, 10){$i$} 
\put(50, 4){\line(-1,0){13}} 
\put(50, 6){\line(-1,0){13}} 
\put(52, 5){\line(-1,1){4}}
\put(52, 5){\line(-1,-1){4}}
                           \put(54, 10){$j$} 
\put(55, 5){\circle{5}}
\put(57, 5){\line(1,0){16}}
\put(73, 3){\line(1,1){4}} \put(74, 10){$k$}
\put(73, 7){\line(1,-1){4}}
\put(75, 5){\circle{5}}  
\put(78, 5){\line(1,1){4}}
\put(78, 5){\line(1,-1){4}}
\put(80, 4){\line(1,0){13}}
\put(80, 6){\line(1,0){13}}
\put(95, 5){\circle{5}} \put(94, 10){$l$} 
\put(100, 5){,}
\end{picture}
\end{flushright}
\vspace{1cm}
(x)
$$
\cwhl E_j, \cwhl E_k, \cwhl E_j, \cwhl E_k, E_i \cwhr \cwhr \cwhr \cwhr =
\cwhl E_k, \cwhl E_j, \cwhl E_k, \cwhl E_j, E_i \cwhr \cwhr \cwhr \cwhr
$$

\begin{flushright}
\setlength{\unitlength}{1mm}
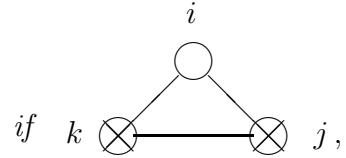
\begin{picture}(35,20)(0,0)
\put(1, 5){if}

\put(13, 3){\line(1,1){4}} \put(8, 4){$k$} 
\put(13, 7){\line(1,-1){4}}
\put(15, 5){\circle{5}}
\put(17, 5){\line(1,0){16}} 
                            
\put(33, 3){\line(1,1){4}}  
\put(33, 7){\line(1,-1){4}}
\put(35, 5){\circle{5}} \put(41, 4){$j\,,$}

\put(23,13){\line(-1,-1){6}} 
\put(25,15){\circle{5}} \put(24, 20){$i$}
\put(27,13){\line(1,-1){6}}

\end{picture}
\end{flushright}
\vspace{1cm}
(5)\quad Relations of $F_i$'s
defined as the same relations as (4).}\newline\newline
{\bf Proof.} Direct calculations.
\par \hfill Q.E.D. \newline\newline\vspace{1cm}
{\bf 5. Relations of 
$D(2;1,x)^{(1)}$, $F_4^{(1)}$ and $G_3^{(1)}$}\newline 
{\bf 5.1.} Let $\Eopip$ $=({\cal E}_0=
\oplus_{i=1}^{3} C\al_i,
 \Pi_0=\{\al_1, \al_2, \al_3\}, p_0)$
be the datum whose Dynkin diagram is:
\begin{center}
Diagram 5.1.1.
\end{center}\begin{center}
\twistdo .
\end{center}
Let $\Epip=\Eopip^{(I)}$ (see 1.5). Then its Dynkin diagram is:
\begin{center}
Diagram 5.1.2.
\end{center}\begin{center}
\twistda .
\end{center} Using
the same argument
as the one for affine $ABCD$-type, we can calculate the defining
relation of $\g\Epip$. In the argument, $\Epidagger$
$=({\cal E^\dagger}=\oplus_{i=0}^{3}C\al^\dagger_i,
 \Pi^\dagger=\{\al^\dagger_i, (0\leq i \leq 3)\})$ is defined by a Dynkin diagram:
\begin{center} 
Diagram 5.1.3.
\end{center}
\begin{center} 
\setlength{\unitlength}{1mm}
\begin{picture}(75,15)(0,0)
\put(5, 5){\circle{5}} \put(4, 10){$0$} 
\put(20, 4){\line(-1,0){13}} 
\put(20, 6){\line(-1,0){13}} 
\put(22, 5){\line(-1,1){4}}
\put(22, 5){\line(-1,-1){4}}
                           \put(24, 10){$1$} 
\put(25, 5){\circle{5}}
\put(27, 5){\line(1,0){16}}
\put(44, 10){$2$}
\put(45, 5){\circle{5}}  
\put(48, 5){\line(1,1){4}}
\put(48, 5){\line(1,-1){4}}
\put(50, 4){\line(1,0){13}}
\put(50, 6){\line(1,0){13}}
\put(65, 5){\circle{5}} \put(64, 10){$3$} 
\put(70, 5){,}
\end{picture}
\end{center} Then the Weyl-group-type 
isomorphism $\sigma(i)$ and $L_i$ $(0 \leq i \leq 3)$
move as follows;
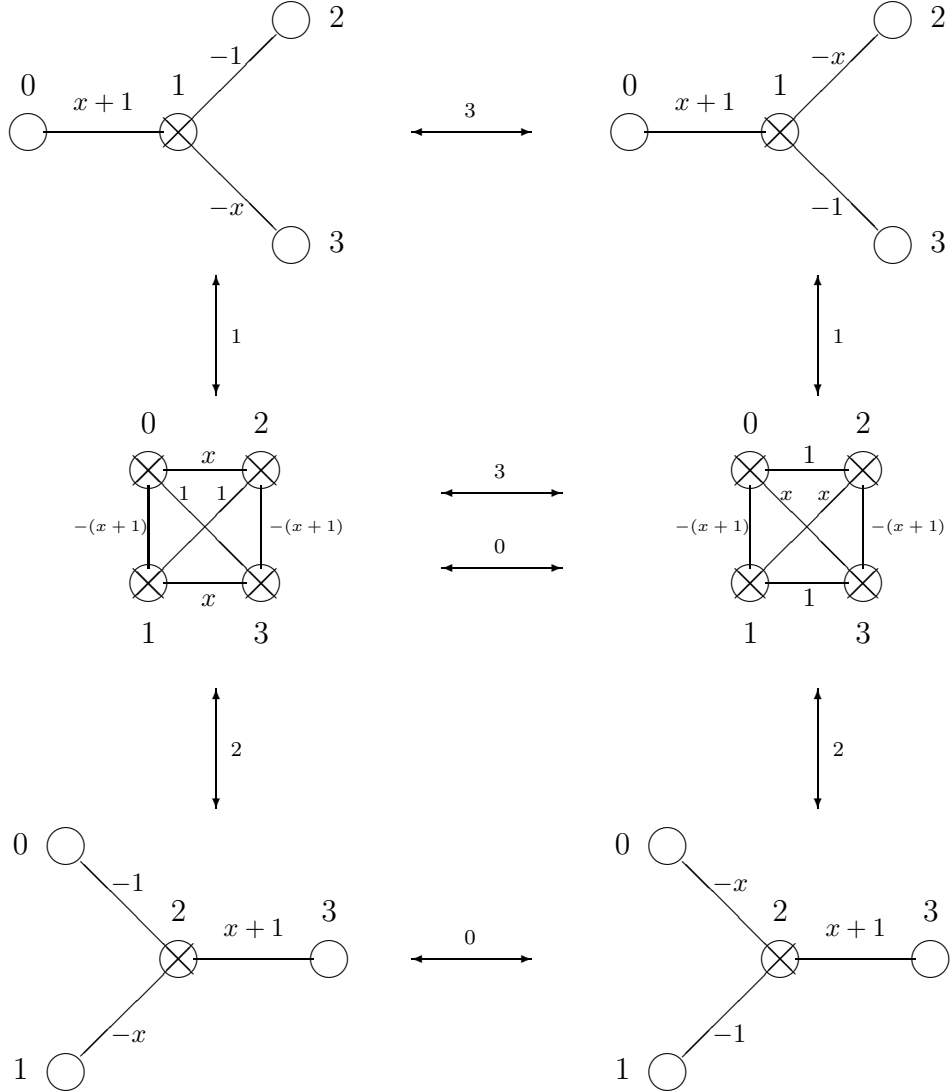
\begin{center} 
Diagram 5.1.4.
\comtwistd 
\end{center} {\bf
Theorem 5.1.1.} {\it Let $\Epip$ be the data of which Dynkin
diagrams in Diagram 5.1.4. Then 
$\bg\Epip$ $=\g\Epip$ and its defining 
relations are ones defined by replacing (vii) of
Theorem 4.1.1 with:
\newline\newline
(vii)
$\cwhl \cwhl \cwhl  E_i, E_j \cwhr ,
 \cwhl  E_j, E_k \cwhr \cwhr ,
 \cwhl  E_j, E_l \cwhr \cwhr 
= x\cwhl \cwhl \cwhl  E_i, E_j \cwhr ,
 \cwhl  E_j, E_l \cwhr \cwhr ,
 \cwhl  E_j, E_k \cwhr \cwhr$

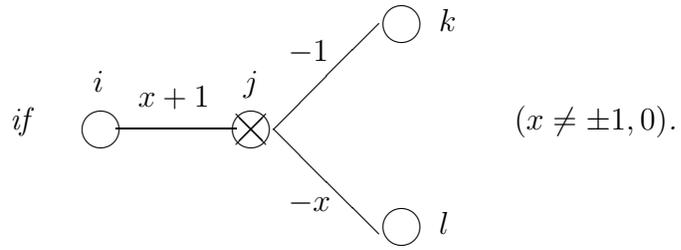
\begin{flushright}
\setlength{\unitlength}{1mm}
\begin{picture}(110,35)(0,0)

\put(23, 15){if}
\put(35, 15){\circle{5}} \put(34, 20){$i$}
\put(37, 15){\line(1,0){16}} \put(40,18){$x+1$}
\put(53, 13){\line(1,1){4}} \put(54, 20){$j$} 
\put(53, 17){\line(1,-1){4}}
\put(55, 15){\circle{5}}  
\put(58, 15){\line(1,1){14}} \put(60, 24){$-1$}
\put(58, 15){\line(1,-1){14}} \put(60, 4){$-x$}
                             \put(80, 28){$k$} 
\put(75, 29){\circle{5}}  
                             \put(80, 1){$l$}
\put(75, 2){\circle{5}}
\put(90, 15){$(x\ne{\pm 1}, 0)$.}

\end{picture}
\end{flushright}
}\vspace{1cm}

{\bf 5.2.} The same argument can still apply to affine $F$-type. Let
 $\Eopip$ $=({\cal E}_0=
\oplus_{i=1}^{4} C\al_i,
 \Pi_0=\{\al_1, \al_2, \al_3, \al_4\}, p_0)$
be the datum whose Dynkin diagram is:
\begin{center}
Diagram 5.2.1.
\end{center}\begin{center}
\sfoo .
\end{center}
 Using
the same argument
as the one for affine $ABCD$-type, we can calculate the defining
relation of $\g\Epip$. In the argument, $\Epidagger$
$=({\cal E^\dagger}=\oplus_{i=0}^{3}, C\al^\dagger_i,
 \Pi^\dagger=\{\al^\dagger_i, (0\leq i \leq 4)\})$ is defined by a Dynkin diagram:
\begin{center}
Diagram 5.2.2.
\end{center}\begin{center}
\foa  .
\end{center} Then the affine $F_4$ Weyl group type 
isomorphism $\sigma(i)$ and $L_i$ $(0 \leq i \leq 4)$
move as follows: (In the diagrams, $i+j+\cdots$
denote $\al^\dagger_i+\al^\dagger_j+\cdots$.)

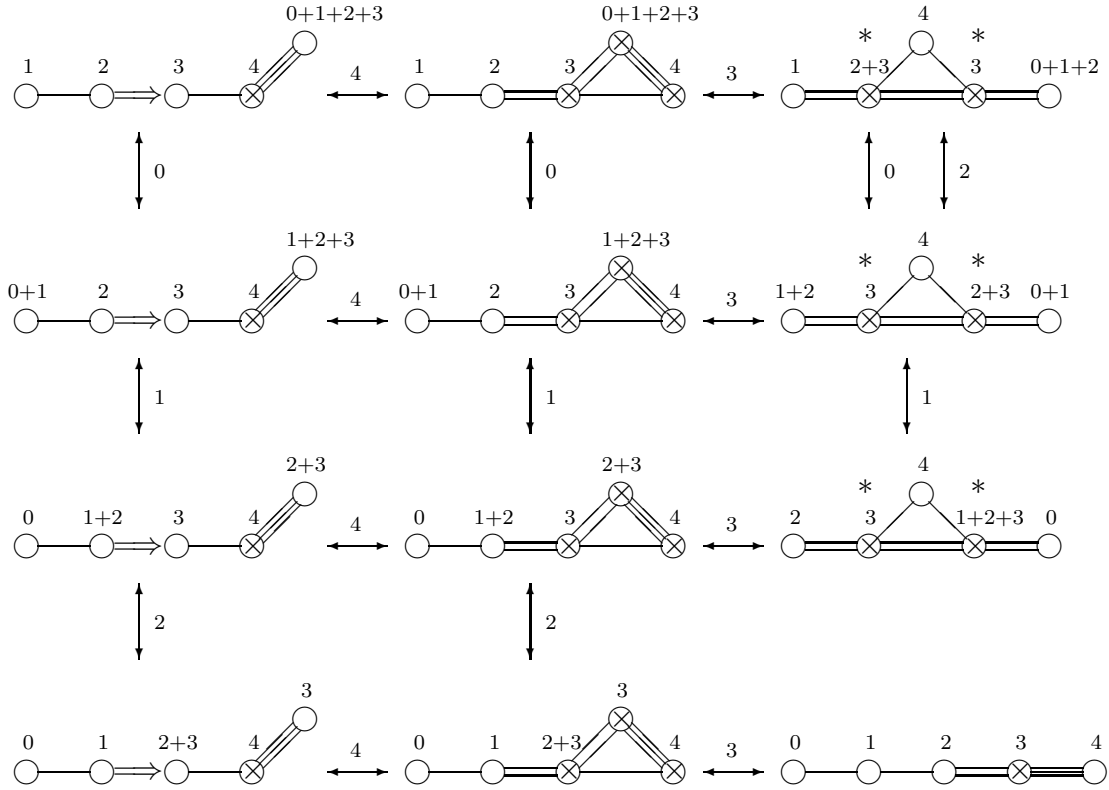
\begin{center}
Diagram 5.2.3.
\end{center}\begin{center}
\fyoncom
\end{center} {\bf
Theorem 5.2.1} {\it Let $\Epip$ be the data of which Dynkin
diagrams existing in Diagram 5.2.3. Then 
$\bg\Epip$ $=\g\Epip$ and its 
relations are defined by adding the following relations
to the ones of
Theorem 4.1.1 with: \newline\newline
(4) (xi) $[\cwhl \cwhl \cwhl \cwhl \cwhl \cwhl \cwhl \cwhl \cwhl 
E_i, E_j \cwhr ,
 E_k \cwhr , E_l \cwhr ,E_k \cwhr , E_j \cwhr ,
 E_k \cwhr , E_l \cwhr ,E_k \cwhr , E_j \cwhr ,
 E_k]=0$
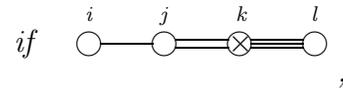
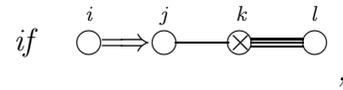
\begin{flushright}
\fxi ,
\end{flushright} (xii)
 $\cwhl \cwhl \cwhl \cwhl \cwhl 
E_l, E_k \cwhr , E_j \cwhr , E_i \cwhr 
, E_k \cwhr , E_j \cwhr 
= 2\cwhl \cwhl \cwhl \cwhl \cwhl 
E_l, E_k \cwhr , E_j \cwhr , E_i \cwhr 
, E_j \cwhr , E_k \cwhr$

\begin{flushright}
\fxii ,
\end{flushright} (5)
\quad Relations of $F_i$'s
defined as the same relations as (4).}
\newline\newline\newline
{\bf 5.3.} 
An argument for affine $G$-type shall be different from the one
for another affine type. Let
 $\Eopip$ $=({\cal E}_0=
\oplus_{i=1}^{4} C\al_i,
 \Pi_0=\{\al_1, \al_2, \al_3\}, p_0)$
be the datum whose Dynkin diagram is:
\begin{center}
Diagram 5.3.1.
\end{center}\begin{center}
\sgoo .
\end{center} Then the positive roots of 
$\g\Eopip$ are:\newline\newline
$\Phi_{0,+}$ $=\{a\al_1+b\al_2+c\al_3\,|$
$(a,b,c)=$ $(1,0,0),$ $(1,1,0),$ $(1,1,1),$ $(1,2,1),$ $(1,3,1),$
$(1,3,2),$ $(1,4,2),$ $(2,4,2),$ $(0,0,1),$ $(0,1,1),$
$(0,3,2),$ $(0,2,1),$ $(0,3,1),$ $(0,1,0)\}$. 
\newline\newline
Let $\Epip=\Eopip^{(I)}$
(see 1.5). Then its Dynkin diagram is: 
\begin{center}
Diagram 5.3.2.
\end{center}\begin{center}
\sgoa .
\end{center} Here the null root $\delta$ of $\g=\g\Epip$
 is given by
$\delta = \al_0+2\al_1+4\al_2+2\al_3$. For $i=0, 2, 3$, we define
 $\sigma (i):{\cal E}\rightarrow{\cal E}$ by 
$\sigma (i)(\gamma)=\gamma - {\frac {2(\gamma, \al_i)}
{(\al_i,\al_i)}}\al_i$. For the set $\Phi_{+}$ of the positive roots of
$\g$, we put
\[ \begin{array}{l}
\Phi_{+}^{\flat} = \Phi_{+} \cup \bigcup_{i=0}^3(\Phi_{+}+\al_i),
\\
 \Phi_{+}^{\sharp}
 =\{\gamma\in\Phi_{+}^{\flat}|(\gamma,\gamma)=2(\rho,\gamma)\}.
\end{array}\] Then as a more precise fact than $(1.2.4)$, it follows:
$$
r_{+}=\bigcup_{\gamma\in\Phi_{+}^{\sharp}}r_{+,\leq \gamma} \eqno{(5.3.1)}
$$ Since $|2(\rho,\delta)|=12$, it is clear that sufficiently large 
element of $\Phi_{+}^{\flat}$ doesn't belong to $\Phi_{+}^{\sharp}$.
By direct calculation, we can get: \newline\newline 
$\Phi_{+}^{\sharp}=\{\al_1,$ $2\al_1,$ $\al_0,$ $\al_3,$ $\al_2,$
$\al_0+3\al_1+4\al_2+\al_3,$ $\al_1+2\al_2,$ $\al_0+\al_1+\al_2+\al_3,$
$2\al_0+3\al_1+5\al_2+2\al_3,$  $\al_0+2\al_1+4\al_2+4\al_3,$
$\al_2+2\al_3,$ $\al_0+2\al_1+2\al_2+\al_3,$ $4\al_2+\al_3\}$.
\newline\newline
Let $r_{+}^{\sharp}$ be the ideal of ${\widetilde {\cal N}}^{+}$
 generated by the relations (i), (ii), (iii) of Theorem 4.1.1 and
$$
2\cwhl \cwhl \cwhl \cwhl \cwhl 
E_0, E_1 \cwhr , E_2 \cwhr , E_3 \cwhr 
, E_1 \cwhr , E_2 \cwhr 
= 3\cwhl \cwhl \cwhl \cwhl \cwhl 
E_0, E_1 \cwhr , E_2 \cwhr , E_3 \cwhr 
, E_2 \cwhr , E_1 \cwhr\,. \eqno{(5.3.2)}
$$ We also
define the ideal $r_{-}^{\sharp}$ of ${\widetilde {\cal N}}^{-}$ in the same way.
By the criterion of Lemma 2.2.2, we can see
 $r_{\pm}^{\sharp} \subset r_{\pm}$. We can also see that
$r^{\sharp}=r_{-}^{\sharp} \oplus r_{+}^{\sharp}$ is an ideal of
${\tilde \g}$. Let
 $\g^{\sharp}={\tilde \g}/r^{\sharp}$.
Then we have a triangular decomposition 
$\g^{\sharp}={\cal N}^{\sharp,+}\oplus{\cal H}\oplus {\cal N}^{\sharp,-}$
where ${\cal N}^{\sharp,\pm}={\widetilde {\cal N}}^{\pm}/r^{\sharp}_\pm$. Since the relations (iii) of 
Theorem 4.1.1 and their $F_i$'s version hold in $\g^{\sharp}$, the 
automorphisms 
$L^{\sharp}_i:\g^{\sharp}\rightarrow\g^{\sharp}$
 $(i = 0, 2, 3)$ given by
$$
 L_i^{\sharp}=\exp(\mbox{ad}E_i)
\exp(\mbox{ad}(-{\frac {2} {(\al_i,\al_i)}}F_i))
\exp(\mbox{ad}E_i)
$$ are well-defined (see also [K1]). Let
 $\g^{\sharp}={\cal H}\oplus(\oplus_{\al\in P_+\cup P_{-}}
\g^{\sharp}_\al)$ be the root space decomposition.
Then we have $L_i^{\sharp}(\g^{\sharp}_\al)=
\g^{\sharp}_{\sigma(i)(\al)}$. By Proposition 2.2.2, we have already 
had the automorphism $L_i:\g\rightarrow\g$ such that
$L_i(\g_\al)=\g_{\sigma(i)(\al)}$. Clearly
$\dim\g^{\sharp}_\al=\dim\g_\al$ if $\al\in
\{\al_0, \al_1, 2\al_1, \al_2, \al_3\}$ Therefore, if $\beta\in P_{+}$
is a minimal element under the order $\leq$ (see 1.2) among 
$\dim\g^{\sharp}_\beta > \dim\g_\beta$, then
$$
\beta\in\Phi_{+}^{\flat}\,\,\mbox{and}\,\,
(\beta,\al_0)\geq 0,\, (\beta,\al_2)\leq 0,\,(\beta,\al_3)\leq 0
\eqno{(5.3.3)}
$$ because $\sigma(i)(\beta)\leq\beta$ $(i=0,2,3)$. 
The unique element satisfying (5.3.3) is
$\al_0+2\al_1+2\al_2+\al_3$. However, using the relation (5.3.2),
we see $\dim \g_{\al_0+2\al_1+2\al_2+\al_3}^{\sharp}\leq 1$.
Hence $\dim \g_{\al_0+2\al_1+2\al_2+\al_3}^{\sharp}$
$=\dim \g_{\al_0+2\al_1+2\al_2+\al_3}$. Hence such $\beta$
doesn't exist. Hence $\g^{\sharp}$ is isomorphic
 to $\g$. \par
By Proposition 2.2.1,we can get other Dynkin diagrams of $\g$.
Those are:

\begin{center}
Diagram 5.3.3.
\end{center}\begin{center}
\gsancom .

\end{center} {\bf Theorem 5.3.1.} {\it Let $\Epip$ be the data of which Dynkin
diagrams existing in Diagram 5.2.3. Then defining relations of
$\g\Epip$ are defined by adding the following relations
to the ones of
Theorem 4.1.1 with: \newline\newline
(4) (xiii)
\quad $[\cwhl \cwhl E_i, E_j \cwhr ,\cwhl \cwhl E_i, E_j \cwhr ,\cwhl 
\cwhl E_i, E_j \cwhr ,E_k \cwhr \cwhr \cwhr , E_j]=0$

\begin{flushright}
\gxiii ,
\end{flushright} (xix) $\cwhl E_j, \cwhl E_k, \cwhl E_k, \cwhl E_j,
E_i \cwhr\cwhr\cwhr\cwhr
=\cwhl E_k, \cwhl E_j, \cwhl E_k, \cwhl E_j,
E_i \cwhr\cwhr\cwhr\cwhr$

\begin{flushright}
\gxix ,
\end{flushright} (xx) $2\cwhl \cwhl \cwhl \cwhl \cwhl 
E_l, E_k \cwhr , E_j \cwhr , E_i \cwhr 
, E_k \cwhr , E_j \cwhr 
= 3\cwhl \cwhl \cwhl \cwhl \cwhl 
E_l, E_k \cwhr , E_j \cwhr , E_i \cwhr 
, E_j \cwhr , E_k \cwhr$

\begin{flushright}
\gxx ,
\end{flushright} (xxi) \newline
$[\cwhl \cwhl \cwhl \cwhl \cwhl \cwhl \cwhl \cwhl
 \cwhl \cwhl \cwhl \cwhl \cwhl 
E_i, E_j \cwhr ,
 E_k \cwhr , E_l \cwhr ,E_k \cwhr , E_j \cwhr ,
 E_k \cwhr , E_l \cwhr ,E_k \cwhr , E_j \cwhr ,
 E_k \cwhr , E_l \cwhr ,E_k \cwhr , E_j \cwhr ,
 E_k]=0$
\begin{flushright}
\gxxi ,
\end{flushright} (5)
\quad Relations of $F_i$'s
defined as the same relations as (4).}\newline\newline 
{\bf 6. Quantization of relations}\newline
{\bf 6.1.} Let $\Epip$ be a datum. Let $C[[h]]$ denote
the $C$-algebra of formal power series in $h$. 
In [Y1], we defined an $h$-adic topological Hopf 
superalgebra $\uhg$ $=U_h({\cal G}\Epip)$ in an abstract manner. 
For the terminologies of topological algebras, Hopf superalgebras etc.,
see [Y1].
Let $\tufhbs$ be a non-topological $C[[h]]$-algebra 
defined with generators $K_\lambda$ $(\lambda \in Z\Pi)$,
$E_\al$ $(\al \in \Pi)$, $\sigma$ and relations:
\[\begin{array}{l}
\sigma^2=1,\, \sigma K_\lambda \sigma =K_\lambda,\,
\sigma E_\al \sigma = (-1)^{p(\al)}E_\al, \\
K_0=1,\,K_\lambda K_\mu= K_{\lambda+\mu},\,
K_\lambda E_\al {K_\lambda}^{-1} =\exp ((\lambda,\al)h)E_\al .
\end{array}
\] 

It is easy to see that $\tufhbs$ is a Hopf algebra with coproduct $\Delta$,
antipode $S$ and counit $\varepsilon$ such that
\[\begin{array}{l}
\Delta (\sigma) = \sigma \otimes \sigma,\, 
\Delta (K_\lambda) = K_\lambda \otimes K_\lambda,\,
\Delta (E_\al) = E_\al \otimes 1 + K_\al \sigma^{p(\al)} \otimes E_\al \\
S(\sigma)=\sigma,\, S(K_\lambda)=K_\lambda^{-1},\, 
S(E_\al)=-K_\al^{-1}\sigma^{p(\al)}E_\al \\
\varepsilon(\sigma)=1,\,\varepsilon(K_\lambda)=1,\,\varepsilon(E_\al)=0.
\end{array}\]  
We note that $\tufhbs$ is not a Hopf superalgebra but
a Hopf algebra.\par
By [Y1], we have:\newline\newline
{\bf Lemma 6.1.1.}
 {\it (i)
$$
E_{\al(1)}\cdots E_{\al(r)}
K_\lambda
\sigma^c\quad\quad (\al(j) \in \Pi,\, 
\lambda \in Z\Pi,\, c \in \{0, 1\}) 
$$ form a $C[[h]]$-basis of $\tufhbs$. In particular, as topological
modules,
$$
\tufhbs \cong {\tilde N}^+ \otimes \chkl \otimes \chs. \eqno{(6.1.1)}
$$ Here ${\tilde N}^+$, $\chkl$ and $\chs$
 denote 
the free algebra generated by $E_\al$ $(\al \in \Pi)$, 
the Laurent polynomial algebra in $K_\al^{\pm 1}$ $(\al \in \Pi)$
and the group ring of $\{1, \sigma\}$ respectively.
\par
(ii) There is a symmetric 
Hopf pairing 
$\langle\,,\,\rangle:\tufhbs \times \tufhbs \rightarrow C[[h]]$ such that;
\newline (a) ${\tilde N}^+$ and $\chkl \otimes \chs$ are orthogonal,
\newline (b) $\langle E_\al, E_\beta \rangle = \delta_{\al,\beta}$
$(\al, \beta \in \Pi)$, 
\newline (c) $\langle K_\lambda\sigma^{c}, K_\mu\sigma^{d} \rangle 
= \exp((\lambda,\mu)h)(-1)^{cd}$.}
\newline\newline
{\bf Remark 6.1.1} In [Y1], we introduced an another topological 
Hopf algebra
${\tilde U}'_{\sqrt{h}}(b_+)^\sigma$. Then $\tufhbs$ 
is given as the non-topological subalgebra of
${\tilde U}'_{\sqrt{h}}(b_+)^\sigma$ generated by $E_\al$,
$K_\al^\pm = \exp (\pm{\sqrt{h}}H'_\al)$ and $\sigma$. 
\newline\newline
{\bf Lemma 6.1.2.} {\it (See [Y1])} {\it
Let $I^+=\{X \in {\tilde N}^+ |
 \langle X,Y \rangle = 0\, (Y \in {\tilde N}^+)\}$.
Then $\ker \langle\,,\rangle \cong I^+ \otimes \chkl \otimes \chs$ under (6.1.1).
In particular, letting $\ufhbs=\tufhbs/\ker \langle\,,\rangle$ and 
$N^+={\tilde N}^+/I^+$, it follows that $\ufhbs \cong
N^+ \otimes \chkl \otimes \chs$.}
\newline\newline
By [Y1] and Lemma 6.1.1, we have:\newline\newline
{\bf Theorem 6.1.1.} {\it For the datum $\Epip$, there is an h-adic 
topological $C[[h]]$-Hopf superalgebra $\uhg$ $=U_h({\cal G}\Epip)$
with generators $H_\lambda$ $(\lambda \in {\cal E})$, $E_\al$, $F_\al$
$(\al \in \Pi)$ with parities $p(H_\lambda)=0$,
$p(E_\al)=p(F_\al)=p(\al)$
 satisfying following (a) and (b):
\newline\newline
(a) In $\uhg$, we have:
\par\vspace{0.5cm}
\quad $\left[H_\lambda, H_\mu\right]=0,\,
\left[H_\lambda,E_\al\right]=(\lambda,\al)E_\al,\,
\left[H_\lambda,F_\al\right]=-(\lambda,\al)F_\al,$ \par\quad
$\left[E_\al,F_\beta\right]=\delta_{\al,\beta}
{\frac {\sinh(hH_\al)} {\sinh(h)}}$. \hfill} (6.1.3) \newline\newline
{\it Put $K_\lambda=\exp(hH_\lambda)$. Then $(\uhg, \Delta, S,
\varepsilon)$ is a Hopf superalgebra such that\newline\newline
\quad
$\Delta (H_\lambda) = H_\lambda \otimes 1 +
1 \otimes H_\lambda,\,
\Delta (E_\al) = E_\al \otimes 1 + K_\al \otimes E_\al$,\par\quad
$\Delta (F_\al) = F_\al \otimes K_\al^{-1} + 1 \otimes F_\al$, \par\quad
$S(H_\lambda)=-H_\lambda,\, 
S(E_\al)=-K_\al^{-1}E_\al,\, S(F_\al)=-F_\al K_\al$, \par\quad
$\varepsilon(H_\lambda)=\varepsilon(E_\al)=\varepsilon(F_\al)=0$.}
\hfill (6.1.4) \newline\newline {\it
(b) Let $\chh$ (resp. $N^+$ or $N^-$) be the (non-topological)
 subalgebra of $\uhg$ generated by $H_\lambda$
(resp. $E_\al$ or $F_\al$). Then $\chh$ is the polynomial ring in 
$H_\lambda \in $ ${\cal H} = {\cal E}^*$.
There are algebra isomorphisms 
${\tilde N}^+/I^+\rightarrow N^+$ $(E_\al\rightarrow E_\al)$ and
$N^+\rightarrow N^-$ $(E_\al \rightarrow F_\al)$. 
There is a topological module isomorphism;
$\uhg \leftarrow (N^- \otimes \chh \otimes N^+)^{\wedge}$ \quad
$(YQX \leftarrow Y\otimes Q\otimes X)$. (Here $(\,)^\wedge$
 denotes completion.)
In particular, $\uhg$ is topologically free as an h-adic module.}
\newline\newline
{\bf 6.2.}  Let $\uhg =\oplus_{\gamma\in Z\Pi} \uhg_\gamma$ be the weight
space decomposition. (Here $\uhg_\gamma
=\{X\in \uhg | [H_\lambda, X]=(\lambda,\gamma)X\,(\lambda\in{\cal E})\}$.
Putting $N^+_\gamma = N^+ \cap \uhg_\gamma$, we have
$N^+=\oplus_{\gamma\in P_+} N^+_\gamma$.\newline\newline
{\bf Lemma 6.2.1.} {\it There is an anti-homomorphism
$t:N^+\rightarrow N^+$ $(E_\al\rightarrow E_\al)$ $(\al \in \Pi)$.}\newline
{\bf Proof.} Let $S$ be the antipode of $\ufhbs$. 
Define $t$ by putting
 $t(X)$ $=(-1)^{\sum_{i<j}p(\al(i))p(\al(j))}$
$\cdot\exp(-\sum_{i<j}(\al(i),\al(j))h)K_\gamma S(X)$ for $X \in {\cal N}^\gamma$ with
$\gamma=\sum \al(i)$ $(\al(i) \in \Pi)$. By Lemma 6.1.2, we can easily check that $t$
is the anti-homomorphism because $S$ is so.
\begin{flushright}  
Q.E.D.
\end{flushright}\par
Let $C((h))$ be the quotient field of $C[[h]]$. 
Put $\fchh = (\chh)^\wedge \otimes_{C[[h]]} C((h))$.  
Let $\fchkl$ be the Laurent polynomial $C((h))$-algebra in $K_\al^{\pm1}$
$(\al \in \Pi)$. Then $\fchkl=\chkl \otimes_{C[[h]]} C((h))$ 
and there is an epimorphism 
$\fchkl \hookrightarrow \fchh$ $(K_\lambda \rightarrow \exp(hH_\lambda))$
$(\lambda\in Z\Pi)$. We define $e:\uhg\rightarrow\fchh$ as the composition:
$$
\uhg\cong(N^- \otimes \chh \otimes N^+)^{\wedge}
\stackrel{\varepsilon\otimes id \otimes\varepsilon}
{\longrightarrow}\chh\hookrightarrow\fchh .
$$ For $\gamma \in Z\Pi$ and $T \in \fchkl$, denote the coefficient 
of $K_\gamma$ of $T$ by $Q_\gamma(T)\in C((h))$. denote the isomorphism
$N^-\rightarrow N^+$ $(F_\al \rightarrow E_\al)$ by $j$. 
Put $q=\exp(h)$ and $q^a=\exp(ah)$. Let $\langle, \rangle:$
$N^+\times N^+ \rightarrow C[[h]]$ be the non-degenerate symmetric 
pairing induced from $\langle, \rangle$ of Lemma 6.1.1 (ii).  
\newline\newline
{\bf Lemma 6.2.2.} {\it  Let $\gamma =\sum_{\al\in\Pi} l_\al \al \in
P_+$ $(l_\al\in Z_+)$. Let $X\in N^+_\gamma$ and $Y \in N^-_\gamma$.
Then $e(XY)\in\fchkl$ and we have:
$$
Q_\gamma(e(XY))={\frac {\langle t(X),j(Y)\rangle} {(q-q^{-1})^{\sum {l_\al}}}}
\in C((h)). \eqno{(6.2.1)}
$$} {\bf Proof.} The first statement is clear.
Describe $\gamma\in P_+$ as two sums 
$\gamma=\sum_{i=1}^r \al(i)=\sum_{i=1}^r \al'(i)$
 $(\al(i), \al'(i)\in\Pi)$. We compare two
calculations of $\langle E_{\al(r)}E_{\al(r-1)}\cdots E_{\al(1)},
E_{\al'(1)}\cdots E_{\al'(r-1)}E_{\al'(r)}\rangle$ and
$E_{\al(r)}\cdots E_{\al(r-1)}E_{\al(r)}
F_{\al'(1)}\cdots F_{\al'(r-1)}F_{\al'(r)}$.
By $\Delta (E_{\al(i)}) = E_{\al(i)} \otimes 1
 + K_{\al(i)} \sigma^{p({\al(i)})} \otimes E_{\al(i)}$, we can calculate:
\[\begin{array}{l}
\langle E_{\al(r)}E_{\al(r-1)}\cdots E_{\al(1)},
E_{\al'(1)}\cdots E_{\al'(r-1)}E_{\al'(r)}\rangle \\
=\quad\langle E_{\al(r)}\otimes E_{\al(r-1)}\cdots E_{\al(1)},
\Delta(E_{\al'(1)}\cdots E_{\al'(r-1)}E_{\al'(r)})\rangle \\
=\quad\displaystyle{\sum_{{1\leq x \leq r}\atop{\al(r)=\al'(x)}}}
\langle E_{\al(r)}\otimes E_{\al(r-1)}\cdots E_{\al(1)}, \\
\quad K_{\sum_{i=1}^{x-1}\al'(i)}
\sigma^{p(\sum_{i=1}^{x-1}\al'(i))}E_{\al'(x)}
K_{\sum_{i=x+1}^r\al'(i)}
\sigma^{p(\sum_{i=x+1}^r\al'(i))}
\otimes E_{\al'(1)} \stackrel{{x}\atop{\vee}}{\cdots} E_{\al'(r)}\rangle \\
=\quad\displaystyle{\sum_{{1\leq x \leq r}\atop{\al(r)=\al'(x)}}}
(-1)^{p(\al(r))p(\sum_{i=x+1}^r\al'(i))}
q^{-(\al(r),\sum_{i=x+1}^r\al'(i))}
\langle E_{\al(r-1)}\cdots E_{\al(1)},
E_{\al'(1)} \stackrel{{x}\atop{\vee}}{\cdots} E_{\al'(r)}\rangle
\end{array}\] 
and
\[\begin{array}{l}
E_{\al(1)}\cdots E_{\al(r-1)}E_{\al(r)}
F_{\al'(1)}\cdots F_{\al'(r)} \\
=(-1)^{p(\al(r))p(\sum_{i=1}^r\al'(i))} 
E_{\al(1)}\cdots E_{\al(r-1)}
F_{\al'(1)}\cdots F_{\al'(r)}E_{\al(r)} \\
+{\frac {1} {(q-q^{-1})}}
\displaystyle{\sum_{{1\leq x \leq r}\atop{\al(r)=\al'(x)}}}
(-1)^{p(\al(r))p(\sum_{i=1}^{x-1}\al'(i))} \\
\Bigl\{q^{-(\al(r),\sum_{i=x+1}^r\al'(i))}
E_{\al(1)}\cdots E_{\al(r-1)}
F_{\al'(1)} \stackrel{{x}\atop{\vee}}{\cdots} F_{\al'(r)}
K_{\al(r)} \\
-q^{(\al(r),\sum_{i=x+1}^r\al'(i))}
E_{\al(1)}\cdots E_{\al(r-1)}
F_{\al'(1)} \stackrel{{x}\atop{\vee}}{\cdots} F_{\al'(r)}
K_{\al(r)}^{-1} \Bigr\}.
\end{array}\]
Since $(-1)^{p(\al(r))p(\sum_{i=1}^{x-1}\al'(i))}$
 $=(-1)^{p(\al(r))+p(\gamma)p(\al(r))}
(-1)^{p(\sum_{x+1}^{r}\al'(i))}$, we can reach (6.2.1) inductively.
\begin{flushright}  
Q.E.D.
\end{flushright}
\vspace{0.5cm}
{\bf Proposition 6.2.1.} {\it If $X\in N^+$ satisfies
$$
 [X, F_\al]=0\quad \mbox{for all}\,\, \al\in\Pi, \eqno{(6.2.2)} 
$$ then $X=0$.} \newline
{\bf Proof.} Since 
$[N^+_\gamma,F_\al]\subset N^+_{\gamma-\al}$,
 we may assume $X\in N^+_\gamma$ for some 
 $\gamma\in P_+$. By (6.2.2), we see that
 $[X, F_{\al(1)}F_{\al(2)}\cdots F_{\al(r)}]=0$ for any
$\{\al(i)\}$ with $\sum_{i=1}^r \al(i)=\gamma$. Hence, by Lemma 6.2.2,
it follow that $\langle X, X_1
\rangle = 0$ for any $X_1\in N^+_\gamma$. It is clear that
 the decomposition $N^+=\oplus_{\gamma\in P_+}N^+_\gamma$ is orthogonal.
Hence $X\in \ker \langle\,,\rangle$ whence $X=0$. 
\begin{flushright}  
Q.E.D.
\end{flushright}
\vspace{0.5cm}
{\bf 6.3.} For $X_\al\in \uhg_\al$, $X_\beta\in \uhg_\beta$, we put:
$$
\whl X_\al, X_\beta \whr 
= X_\al X_\beta -(-1)^{p(\al)p(\beta)}q^{-(\al,\beta)}X_\beta X_\al.
$$ Let $[x]={\frac {\sinh(xh)} {\sinh(h)}} \in C[[h]]$. 
By Proposition 6.2.1 and direct calculation,
we have:\newline\newline
{\bf Proposition 6.3.1.} {\it In $N^+$, we have:\newline
\newline
(i) $[E_i, E_j]=0$ \hfill $\mbox{if}\, (\al_i,\al_j)=0, (i\ne j)$ 
\newline
(ii) $[E_i, E_i]=0$ \hfill  \aii
\newline\newline
(iii) $\whl E_i, \whl E_i, \ldots , \whl E_i, E_j
 \whr\!\ldots\!\whr\whr =0 $\, 
($E_i\, \mbox{appear}\, 1-{\frac {2(\al_i,\al_j)} {(\al_i,\al_i)}}
\mbox{times}$) \par
\hfill $\mbox{if}\, (\al_i,\al_i)\ne 0\, \mbox{and}\, 
p(\al_i){\frac {2(\al_j,\al_i)} {(\al_i,\al_i)}}\,
\mbox{is even,}$
\newline
(iv) $[\whl \whl E_i, E_j \whr ,E_k \whr , E_j]=0$ \hfill  \aiv 
\newline
(v)  $[\whl \whl E_i, E_j \whr ,\whl 
\whl E_i, E_j \whr ,E_k \whr \whr , E_j]=0$ \hfill \av 
\newline
(vi) $[\whl \whl \whl \whl \whl 
E_i, E_j \whr , E_k \whr , E_l \whr ,
E_k \whr , E_j \whr , E_k]=0$ \hfill  \avi
\newline\newline
(vii) $(-1)^{p(\al_i)p(\al_k)}
[(\al_i,\al_k)]
\whl \whl E_i, E_j \whr , E_k \whr
=
(-1)^{p(\al_i)p(\al_j)}
[(\al_i,\al_j)]
\whl \whl E_i, E_k \whr , E_j \whr$
\par \hfill \avii
\newline\newline
(viii)
$\whl \whl \whl  E_i, E_j \whr ,
 \whl  E_j, E_k \whr \whr ,
 \whl  E_j, E_l \whr \whr 
= [x]\whl \whl \whl  E_i, E_j \whr ,
 \whl  E_j, E_l \whr \whr ,
 \whl  E_j, E_k \whr \whr$\par
\hfill \aviii
\newline\newline
(ix) $\whl \whl E_k, \whl E_l, \whl E_k, E_j, 
 \whr \whr \whr,
\whl E_k, \whl E_l, \whl E_k, \whl E_j, E_i, 
 \whr \whr \whr \whr \whr, 
E_j]$\par 
$= [2]\whl \whl E_k, E_j \whr,
\whl \whl E_k, \whl E_j, E_i \whr \whr, 
\whl E_k, \whl E_l, \whl E_k, \whl E_j, E_i, 
 \whr \whr \whr \whr \whr \whr$ \par 
\hfill \aix
\newline\newline
(x) $\whl E_j, \whl E_k, \whl E_j, \whl E_k, E_i \whr \whr \whr \whr =
\whl E_k, \whl E_j, \whl E_k, \whl E_j, E_i \whr \whr \whr \whr$ \par
\hfill \ax 
\newline\newline
(xi) $[\whl \whl \whl \whl \whl \whl \whl \whl \whl 
E_i, E_j \whr ,
 E_k \whr , E_l \whr ,E_k \whr , E_j \whr ,
 E_k \whr , E_l \whr ,E_k \whr , E_j \whr ,
 E_k]=0$
\begin{flushright}
\fxi ,
\end{flushright} (xii)
 $\whl \whl \whl \whl \whl 
E_l, E_k \whr , E_j \whr , E_i \whr 
, E_k \whr , E_j \whr 
= [2]\whl \whl \whl \whl \whl 
E_l, E_k \whr , E_j \whr , E_i \whr 
, E_j \whr , E_k \whr$

\begin{flushright}
\fxii ,
\end{flushright} (xiii)
\quad $[\whl \whl E_i, E_j \whr ,\whl \whl E_i, E_j \whr ,\whl 
\whl E_i, E_j \whr ,E_k \whr \whr \whr , E_j]=0$

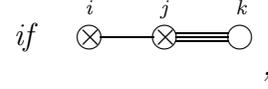
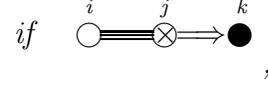
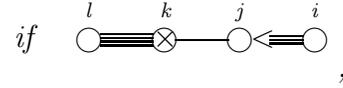
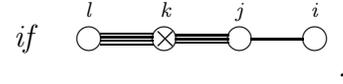
\begin{flushright}
\gxiii ,
\end{flushright} (xix) $\whl E_j, \whl E_k, \whl E_k, \whl E_j,
E_i \whr\whr\whr\whr
=\whl E_k, \whl E_j, \whl E_k, \whl E_j,
E_i \whr\whr\whr\whr$

\begin{flushright}
\gxix ,
\end{flushright} (xx) $[2]\whl \whl \whl \whl \whl 
E_l, E_k \whr , E_j \whr , E_i \whr 
, E_k \whr , E_j \whr 
= [3]\whl \whl \whl \whl \whl 
E_l, E_k \whr , E_j \whr , E_i \whr 
, E_j \whr , E_k \whr$

\begin{flushright}
\gxx ,
\end{flushright} (xxi) \newline
$[\whl \whl \whl \whl \whl \whl \whl \whl
 \whl \whl \whl \whl \whl 
E_i, E_j \whr ,
 E_k \whr , E_l \whr ,E_k \whr , E_j \whr ,
 E_k \whr , E_l \whr ,E_k \whr , E_j \whr ,
 E_k \whr , E_l \whr ,E_k \whr , E_j \whr ,
 E_k]=0$
\begin{flushright}
\gxxi .
\end{flushright}

} In 6.4, we shall describe how we calculate the relations (i)-(xix).
\newline\newline
{\bf 6.4.} Let ${\cal S}$ be a $C[[h]]$ or $C((h))$-superalgebra.
For $a \in C[[h]]^\times$, we put:
$$
[X, Y]_a = XY-(-1)^{p(X)p(Y)}aYX \quad (X,\, Y \in {\cal S}).
$$
Then we have
$$
[[X,Y]_a, Z]_b = [X, [Y, Z]_c]_{abc^{-1}}
+(-1)^{p(Y)p(Z)}c[[X,Z]_{bc^{-1}},Y]_{ac^{-1}}, \eqno{(6.4.1)}
$$ 
and
$$
[X,[Y, Z]_a]_b = [[X, Y]_c, Z]_{abc^{-1}}
+(-1)^{p(X)p(Y)}c[Y,[X,Z]_{bc^{-1}}]_{ac^{-1}}\,. \eqno{(6.4.2)}
$$
Hence, for $\uhg$ and $X_\nu \in \uhg_\nu$, $X_\mu \in \uhg_\mu$,
$X_\eta \in \uhg_\eta$, we have:
$$
\whl\whl X_\nu,X_\mu\whr, X_\eta\whr
 = \whl X_\nu, \whl X_\mu, X_\eta \whr\whr
+(-1)^{p(\mu)p(\eta)}q^{-(\mu,\eta)}
[\whl X_\nu,X_\eta\whr,X_\mu]_{q^{(\mu,\eta-\nu)}}, \eqno{(6.4.3)}
$$ 
and
$$
\whl X_\nu, \whl X_\mu, X_\eta\whr\whr
 = \whl\whl X_\nu, X_\mu \whr, X_\eta \whr
+(-1)^{p(\mu)p(\nu)}q^{-(\mu,\nu)}
[X_\mu,\whl X_\nu,X_\eta \whr]_{q^{(\mu,\nu-\eta)}}. \eqno{(6.4.4)}
$$ \par
We can get the relations (i)-(xix) of Proposition 6.3.1 by
Proposition 6.2.1 and direct calculation
using (6.4.3-4). Here we only show how to get (ix) because the
other relations can be also gotten  similarly.
 We replace the letters
$i$, $j$, $k$, $l$ with 0, 1, 2, 3. We assume $(\al_1,\al_1)=-2$. Then the 
diagram can be rewritten as:
\begin{center}
\aixa
\end{center}\vspace{0.5cm}
Put $E_{...dcba}$ $=\whl\ldots\whl E_d, \whl E_c, \whl E_b , E_a 
\whr\whr\whr\mbox{\scriptsize{...}}\whr$. Then (ix) is rewritten as:
$$
-\whl\whl E_{2321}, E_{23210} \whr, E_1 \whr
+(q+q^{-1})\whl E_{21} \whl E_{321}, E_{23210} \whr\whr = 0
\eqno{(6.4.5)}
$$ We 
denote the LHS of (6.4.5) by ${\cal X}$. Showing (6.2.2) is equivalent to
showing $\whl {\cal X}, F_aK_a^{-1}\whr=0$
for all $0\leq a\leq 3$
 where we put $K_a=K_{\al_a}$. We note:
$$
\whl E_\al, F_\beta K_\beta^{-1} \whr
=\delta_{\al,\beta}{\frac {1-K_\beta^{-2}} {q-q^{-1}}}.
\eqno{(6.4.6)}
$$ First we show $\whl {\cal X}, F_3K_3^{-1}\whr=0$.
In following equations, 
 $LHS\sim  RHS$ mean  
$LHS = a\cdot RHS$ for some $a\in C[[h]]^\times$. By (6.4.3) and (6.4.6),
$\whl E_{321}, F_3K_3^{-1}\whr$ 
$\sim [{\frac {1-K_3^{-2}} {q-q^{-1}}},E_{21}]_{q^{(2+2)}}$ 
$\sim E_{21}$. Hence
\begin{eqnarray*}
\lefteqn{\whl E_{2321}, F_3K_3^{-1} \whr}\\
& &\sim \whl E_2, \whl E_{321}, F_3K_3^{-1}\whr\whr \quad \,
\mbox{by (6.4.3)}\\ 
& &\sim \whl E_2, E_{21}\whr \\
& &=0\quad \,
\mbox{by (ii) of Proposition 6.3.1}.
\end{eqnarray*} Hence we also have
$\whl E_{23210}, F_3K_3 \whr =0$.
Hence we have
$\whl \whl \whl E_{2321}, E_{23210} \whr, E_1 \whr, F_3K_3^{-1} \whr =0.$
On the other hand,
\begin{eqnarray*}
\lefteqn{\whl \whl E_{21} \whl E_{321}, 
E_{2321}\whr\whr, F_3K_3^{-1} \whr}\\
& &\sim \whl E_{21} \whl E_{21}, E_{23210} \whr \whr
\quad \mbox{by (6.4.3) and}\,\, E_{321}^2=0 \\
& & = 0\quad \mbox{since}\,\, E_{21}^2=0.
\end{eqnarray*} Then we have $\whl {\cal X}, F_3K_3 \whr =0$.\par
Next we show $\whl {\cal X}, F_2K_2 \whr =0$.
First we calculate: \newline $\whl E_{2321}, F_2K_2^{-1} \whr$
$=\whl \whl E_2, E_{321}\whr, F_2K_2^{-1} \whr$
$=-q^{-1}[{\frac {1-K_2^{-2}} {q-q^{-1}}}, E_{321}]_{q^{(1+1)}}$ $=E_{321}$,
\newline\newline
 $\whl E_{23210}, F_2K_2^{-1} \whr$ $=E_{3210}$,
\newline\newline
$\whl\whl E_{2321}, E_{23210} \whr, F_2K_2^{-1} \whr$
$=\whl E_{2321}, E_{3210} \whr$ 
$+q^{-1}[E_{321},E_{23210}]_{q^{(1-0)}}$ \par
$=-q^{-2} [E_{23210}, E_{321}]_{q^{(2+1)}}$ 
$+q^{-1}[E_{321},E_{23210}]_{q^{(1-0)}}$ 
(since $E_{321}^2=0$) \par
$=(q+q^{-1})\whl E_{321},E_{23210} \whr$,
\newline\newline
$\whl\whl E_{321},E_{23210} \whr, F_2K_2^{-1} \whr =0$
and $\whl E_{21}, F_2K_2^{-1} \whr = E_1$.
\newline\newline
Using these, we have:
\begin{eqnarray*}
\lefteqn{\whl {\cal X}, F_2K_2^{-1} \whr} \\
& &= -q[(q+q^{-1})\whl E_{321}, E_{23210} \whr, E_1 ]_{q^{(-1-2)}} 
- (q+q^{-1})q^{-2}
[E_1, \whl E_{321}, E_{23210} \whr]_{q^{(2+1)}} \\
& & =0 \\
\end{eqnarray*} 
Similarly we can show $\whl {\cal X}, F_1K_1^{-1}\whr$
$=\whl {\cal X},F_0K_0^{-1}\whr$ $=0$. By Proposition 6.2.1, it follows that
${\cal X}=0 \in N^+$. \par
By similar calculation, we can get the relations (i)-(xix) of
Proposition 6.3.1. \newline\newline
{\bf 6.5.} Let $({\cal A}, \Delta)$  be a cocommutative Hopf
C-superalgebra. Let
$$ 
P({\cal A})= \{x\in{\cal A}| \Delta(X)
=X\otimes 1 + 1\otimes X\}.
$$ Then $P({\cal A})$ is a Lie $C$-superalgebra with a bracket
$[\,,\,]$ given by $[X,Y]=XY-(-1)^{p(X)p(Y)}YX$. \par
Let ${\cal G}$ be a Lie $C$-superalgebra and $U({\cal G})$ its universal
enveloping superalgebra. Then $U({\cal G})$ is a cocommutative Hopf
$C$-superalgebra with coproduct $\Delta$ such that
$\Delta(X)
=X\otimes 1 + 1\otimes X$ for $X\in U({\cal G})$. It is known that:
\newline\newline
{\bf Theorem 6.5.1}(Milnor-Moor [MM]) {\it  
Let ${\cal CSH}$ be the category of cocommutative Hopf $C$-superalgebras. 
and ${\cal SL}$ the category of Lie $C$-superalgebras.
Define morphisms ${\cal P}$ and ${\cal U}$ by
${\cal P}: {\cal CSH}\rightarrow{\cal SL}$ (${\cal A}\rightarrow P({\cal A})$),
${\cal U}: {\cal SL}\rightarrow{\cal CSH}$ (${\cal G}\rightarrow U({\cal G})$).
Then ${\cal PU}=id_{{\cal CSH}}$ and ${\cal UP}=id_{{\cal SL}}$.}
\vspace{1cm}\par
As an immediate consequence of Theorem 6.5.1, we have:
\newline\newline
{\bf Lemma 6.5.1} {\it For a datum $\Epip$, let ${\cal G}={\cal G}\Epip$ and
$U_0({\cal G})$ a cocommutative Hopf $C$-superalgebra defined by
$U_0({\cal G})=U_h({\cal G})/hU_h({\cal G})$. Then
there is an epimorphism
$$
\phi: U_0({\cal G})\rightarrow U({\cal G})\quad
(H,\, E_\al,\, F_\al\rightarrow H,\, E_\al,\, F_\al).
$$}\newline
{\bf Proof.} Let ${\cal G}_0=P(U_0({\cal G}))$. By theorem 6.5.1,
we have $U_0({\cal G})=U({\cal G}_0)$. Hence ${\cal G}_0$ should be a Lie 
$C$-superalgebra generated by $H,\,E_\al,\,F_\al$. By Theorem 6.1.1 (a),  
$H,\,E_\al,\,F_\al$ satisfy (1.2.1-3). Since $\uhg$
has the triangular decomposition by Theorem 6.1.1 (b), 
$U_0({\cal G})$ also has a triangular decomposition. In particular,
${\cal H}$ can be embedded into ${\cal G}_0 (\subset U_0({\cal G}))$.
By definition of the Kac-Moody Lie superalgebra ${\cal G}$ (see [K1]), we have
epimorphism $\phi_{|{\cal G}_0}: {\cal G}_0 \rightarrow {\cal G}$ 
$(H,\, E_\al,\, F_\al\rightarrow H,\, E_\al,\, F_\al)$. Hence we have $\phi$.
\par
\hfill Q.E.D.  
\newline
{\bf 6.6.} {\bf Theorem 6.6.1} {\it Let $\Epip$ be a datum of affine type.
Assume ${\bar {\cal G}}\Epip$ $={\cal G}\Epip$,i.e.,
${\cal G}\Epip$ is not of ${\hat sl}(m|m)^{(i)}$ $(i=1,2,4)$.
Put ${\cal G}={\cal G}\Epip$. Then the defining relations
of $\uhg$ are given by \newline
(i) The relations of Theorem 6.1.1 (a),\newline
(ii) The relations of Proposition 6.3.1 (i)-(xix),\newline
(iii) The same relations of $F_\al$'s as (ii).}\newline\newline    
{\bf Proof.} By Theorem 6.1.1,
 Proposition 6.3.1, $\uhg$ satisfies the relations (i)-(iii).
Therefore, by Serre type theorems of Chapters 4 and 5,
we have an epimorphism
$\psi: U({\cal G})\rightarrow U_0({\cal G})$ 
$(H,\, E_\al,\, F_\al\rightarrow H,\, E_\al,\, F_\al)$.
By Lemma 6.5.1, $\psi$ should be the inverse map of $\phi$. 
Then $U_0({\cal G})=U({\cal G})$. The relations given by putting $h=0$ on
(i)-(iii) are the defining relations of $U({\cal G})$. Hence, by the topologically 
freedom of $\uhg$, the relations (i)-(iii) should be the defining relations of
$\uhg$.\par
\hfill Q.E.D.
\newline\newline
{\bf 6.7} {\bf Lemma 6.7.1.} {\it Let $\uhgd$ be the topologically free Hopf 
C[[h]]-superalgebra with generators $\{H, E_\al, F_\al\}$ ($\al\in\Pi$)
such that \newline\newline
(1) $\{H, E_\al, F_\al\}$ satisfy (6.1.3-4). \newline
(2) The map ${\cal H}\rightarrow\uhgd$ $(H\rightarrow H)$
is injective. \newline\newline
Then there is a Hopf superalgebra epimorphism 
$$
{\dot J}: \uhgd\rightarrow\uhg \quad
 (H, E_\al, F_\al\rightarrow H, E_\al, F_\al).
$$}
{\bf Proof.} By (2), the topological freedom of $\uhgd$ and Theorem 6.5.1,
$\chh$ is embedded into $\uhgd$. Let $U_h({\dot {\cal N}}^+)$, $
U_h({\dot {\cal N}}^-)$ be 
non-topological subalgebras generated by $E_\al$, $F_\al$ respectively.
Put $U_0({\dot {\cal G}})=\uhgd/h\uhgd$ and 
$U_0({\dot {\cal N}}^\pm)=
U_h({\dot {\cal N}}^\pm)/hU_h({\dot {\cal N}}^\pm)$.
By Minor-Moor's Theorem 6.5.1, $U_0({\dot {\cal G}})$ (resp. 
$U_0({\dot {\cal N}}^\pm)$) 
is the universal enveloping algebra
$U({\dot {\cal G}})$ (resp. $U({\dot {\cal N}}^\pm)$) 
of a Lie $C$-superalgebra ${\dot {\cal G}}
={\cal P}(U_0({\dot {\cal G}}))$.
(resp. ${\dot {\cal N}}^\pm =
{\cal P}(U_0({\dot {\cal N}}^\pm))$
By (2), ${\cal H}$ is embedded into
${\dot {\cal G}}$. Hence we have the triangular decompositions
${\dot {\cal G}}={\dot {\cal N}}^- \oplus {\cal H} \oplus {\dot {\cal N}}^+$
and $U({\dot {\cal G}})=U({\dot {\cal N}}^-) \otimes
U({\cal H}) \oplus U({\dot {\cal N}}^+)$. By the topological freedom of $\uhgd$,
we have the triangular decomposition
$\uhgd = U_h({\dot {\cal N}}^-) {\hat \otimes} \chh {\hat \otimes}
U_h({\dot {\cal N}}^+)$. \par
Let ${\dot I}^+=\ker ({\tilde N}^+\rightarrow {\dot {\cal N}}^+)$.
For $\gamma\in P_+$,
put ${\tilde N}^+_\gamma$ (resp. ${\dot I}^+_\gamma$)
$=\{X\in {\tilde N}^+
\,\mbox{(resp ${\dot I}$)}\, | [H_\lambda, X]=(\lambda,\gamma)X\}$.
Then we have 
${\tilde N}^+=\oplus_{\gamma\in P_+} {\tilde N}^+_\gamma$ and
${\dot I}^+=\oplus_{\gamma\in P_+} {\dot I}^+_\gamma$.
{\it Keep notations in 6.2.} Since ${\dot I}^+$ is an ideal of ${\tilde N}^+$,
by the triangular decomposition of $\uhgd$, 
$$
e([{\dot I}^+_\gamma, F_{\al(1)}\cdots F_{\al(r)}])=0\quad
\mbox{for}\, \sum \al (i) = \gamma\,\,(\al (i) \in \Pi).
$$
Hence, by a $\uhgd$-version of Lemma 6.2.2, since
${\tilde N}^+=\oplus_{\gamma\in P_+} {\tilde N}^+_\gamma$
is orthogonal with respect to $\langle\,,\,\rangle$, we have
${\dot I}^+_\gamma\subset I^+=\ker\langle\,,\,\rangle$. 
Then we have the algebra epimorphism
$U_h({\dot {\cal N}}^+)\rightarrow N^+$. $(E_\al\rightarrow E_\al$). \par
Next we should show the existence of the epimorphism
${\dot N}^-\rightarrow N^-$. However a Hopf superalgebra 
$(\uhgd, s\tau\circ\Delta, S^{-1}, \varepsilon)$ with generators
$\{-H_\lambda, F_\al, (-1)^{p(\al)}E_\al \}$ also satisfies (1) and (2).
(Here $s\tau(X\otimes Y)= (-1)^{p(X)p(Y)}Y \otimes X)$. Then the same argument
can be applied for the subalgebra generated by $F_\al$. Hence we can show 
the existence. \par
Eventually we get a $C[[h]]$-module surjective map:
$$
{\dot J}: \uhgd= ({\dot N}^-\otimes \chh \otimes {\dot N}^+)^\wedge
\uhg= (N^-\otimes \chh \otimes N^+)^\wedge.
$$
Considering under the two triangular decompositions,
it is clear that ${\dot J}$ preserve product. Hence ${\dot J}$
is the algebra epimorphism. Clearly ${\dot J}$
is the Hopf superalgebra epimorphism.
\par\hfill Q.E.D.
\newline\newline
{\bf 7. Quantization of Weyl-group-type isomorphisms}\newline
{\bf 7.1}
Let $({\cal C}, \Delta, S, \varepsilon)$ be a topological Hopf
 $C[[h]]$-algebra.
Define $\tau : {\cal C} {\hat \otimes} {\cal C} \rightarrow 
{\cal C} {\hat \otimes} {\cal C}$
by $\tau(x\otimes y)=y \otimes x$. Let $\Delta'=\tau\circ\Delta$. Let
${\cal C}_0$ be a Hopf subalgebra of ${\cal C}$. Let
 $R=\sum a_i \otimes b_i$ be an invertible element of
${\cal C}_0\otimes{\cal C}_0$ satisfying:
\newline\par
\hspace{2.5cm} $R\Delta(x)R^{-1}=\Delta'(x)$, \hfill (7.1.1) 
\par
$(\Delta\otimes I)(R)=R_{13}R_{23}$,\,\,
$(I\otimes \Delta)(R)=R_{13}R_{12}$. \hfill (7.1.2)
\newline\newline
where
$R_{12} =R \otimes I$, $R_{23}=I\otimes R$ and
$R_{13} =\sum a_i \otimes I \otimes b_i$. \par
By Drinfeld[D2], we have known:
\newline\newline
{\bf Proposition 7.1.1.}(Drinfeld[D2]){\it (i) $R$ satisfies:
\newline\par
\hspace{2.5cm}$R_{12}R_{13}R_{23}=R_{23}R_{13}R_{13}$, \hfill (7.1.3)
\par
\hspace{1.5cm}$(S\otimes I)(R)=R^{-1}=(I\otimes S^{-1})(R)$,\hfill (7.1.4)
\par
\hspace{2cm}$(\varepsilon\otimes I)(R)=1
=(I\otimes \varepsilon)(R)$.\hfill (7.1.5)
\newline\newline
(ii) For $R=\sum a_i \otimes b_i$, following equations hold in 
${\cal C}$:\newline\par
\hspace{1.5cm}$\sum a_iS^{-2}(b_i)=
\sum S(a_i)S^{-1}(b_i)=\sum S^2(a_i)b_i$, \hfill (7.1.6)
\par
\hspace{2.5cm}$
\sum a_iS(b_i)=\sum S^{-1}(a_i)b_i$. \hfill (7.1.7)
\newline\newline
Let $u_4$, $v_4\in{\cal C}$ be the elements of (7.1.6), (7.1.7)
respectively. Then $u_4v_4=1=v_4u_4$.}
\newline\newline
{\bf 7.2.} {\bf Proposition 7.2.1.} {\it  Keep notations in 7.1.
For the Hopf algebra ${\cal C}$ and the element $R=\sum a_i \otimes b_i$ satisfying (7.1.1-2),
there is an another Hopf algebra structure
$({\cal C}^{(R)})=({\cal C},\Delta^{(R)}, S^{(R)}, \varepsilon)$ given by:
$$
\Delta^{(R)}(x)=R\Delta(x)R^{-1},\,S^{(R)}(x)=u_4^{-1}S(x)u_4.
$$}
{\bf Proof.} First we show $(I\otimes \Delta^{(R)})\circ\Delta^{(R)}
=(\Delta^{(R)}\otimes I)\circ\Delta^{(R)}$. By (7.1.1-2), We have:

\begin{eqnarray*}
\lefteqn{(I\otimes \Delta^{(R)})\circ\Delta^{(R)}(x)}\\
& &= R_{23}(I\otimes \Delta)(R\Delta(x)R^{-1})R_{23}^{-1}\\ 
& &= R_{23}R_{13}R_{12}(I\otimes \Delta)(\Delta(x))
R_{12}^{-1}R_{13}^{-1}R_{23}^{-1}\\ 
& &= R_{12}R_{13}R_{23}(\Delta\otimes I)(\Delta(x))
R_{23}^{-1}R_{13}^{-1}R_{12}^{-1}\\ 
& &=(\Delta^{(R)}\otimes I)\circ\Delta^{(R)}(x).
\end{eqnarray*}  
Let $m: {\cal C} \otimes {\cal C}\rightarrow {\cal C}$ be the multiplication,
which is defined by $m(x\otimes y)=xy$. Next we show 
$m\circ (I\otimes S^{(R)})\circ\Delta^{(R)}=\varepsilon$
$=m\circ (S^{(R)}\otimes I)\circ\Delta^{(R)}$. By (7.1.4) and
(7.1.6-7), for $x\in {\cal C}$ with $\Delta(x)=\sum x_i^{(1)}\otimes
x_i^{(2)}$, we have:

\begin{eqnarray*}
\lefteqn{(I\otimes S^{(R)})\circ\Delta^{(R)}(x)}\\
& &=\sum m((1\otimes u_4^{-1})(I\otimes S)
(a_ix_j^{(1)}S(a_l)\otimes b_ix_j^{(2)}b_l)(1\otimes u_4)) \\
& &=\sum a_ix_j^{(1)}S(a_l) \cdot a_yS(b_l)S(b_i)S(x_j^{(2)})S(b_i)u_4 \\
& &=\sum a_ix_j^{(1)}S(x_j^{(2)})S(b_i)u_4\quad \mbox{since}\,\,
(S\otimes I)(R)R=1 \\
& &=\varepsilon(x)\sum a_iS(b_i)u_4=\varepsilon(x)v_4u_4=
\varepsilon(x)
\end{eqnarray*}
and

\begin{eqnarray*}
\lefteqn{(S^{(R)}\otimes I)\circ\Delta^{(R)}(x)}\\
& &=\sum m((u_4^{-1}\otimes 1)(S\otimes I)
(a_ix_j^{(1)}S(a_l)\otimes b_ix_j^{(2)}b_l)(u_4\otimes 1)) \\
& &=\sum u_4^{-1}S^2(a_l)S(x_j^{(1)})S(a_i)S(a_y)S^{-1}(b_y)\cdot
b_ix_j^{(2)}b_l \\
& &=\sum u_4^{-1}S^2(a_l)S(x_j^{(1)})x_j^{(2)}b_l \quad \mbox{since}\,\,
(I\otimes S^{-1})(R)R=1 \\
& &=\varepsilon(x)u_4^{-1}\sum S^2(a_l)b_l=\varepsilon(x)u_4^{-1}u_4=
\varepsilon(x).
\end{eqnarray*}
To show other formulae of the axiom of the Hopf algebra are easy.\par
\hfill Q.E.D.
\newline\newline
{\bf 7.3.} Let $\Epip$ be a datum and ${\cal G}$ $={\cal G}\Epip$.
Let $\uhg$ be a topological Hopf $C[[h]]$-superalgebra
 introduced in 6.1. Put $\uhgs = \uhg \otimes_{C[[h]]}\chs$.
Then $\uhgs$ is an algebra with a formula $\sigma X \sigma $ 
$= (-1)^{p(X)}X$ $(X\in \uhg)$. By [Y1], 
$(\uhgs, \Delta , S, \varepsilon)$ is a Hopf algebra such that 

\[\begin{array}{l}
\Delta (H) = H \otimes 1 + 1 \otimes H,\,
\Delta (E_\al) = E_\al \otimes 1 + K_\al\sigma^{p(\al)} \otimes E_\al,\, \\
\Delta (F_\al) = F_\al \otimes K_\al^{-1} + \sigma^{p(\al)} \otimes F_\al, \\
S(H)=-H,\, 
S(E_\al)=-K_\al^{-1}\sigma^{p(\al)}E_\al,\, 
S(F_\al)=-\sigma^{p(\al)}F_\al K_\al \\
\varepsilon(H)=\varepsilon(E_\al)=\varepsilon(F_\al)=0.
\end{array}\]
Put $\uhhs =C[[h]][{\cal H}]\otimes \chs$. Then $\uhhs$ is a 
Hopf subalgebra of $\uhgs$. Put $t_0=\sum H_{\delta_i}
\otimes H_{\delta_i}$ $\in {\cal H}\otimes{\cal H}$ where $\{\delta_i\}$
is a $C$-basis of ${\cal H}$ such that $(\delta_i,\delta_j)=\delta_{ij}$. Then,
by the quantum double construction (see [D] (also [Y1])),
$$
R_T={\frac {1} {2}}(\sum_{c,d=0,1} (-1)^{cd}\sigma^c\otimes \sigma^d)
\cdot\exp(-ht_0)\in\uhhs\otimes\uhhs
$$
satisfies (7.1.1-2). Clearly $R_T^{-1}$ also satisfies (7.1.1-2). 
\par
 For $t\in C[[h]]$ and $n>0$, we put
$\{n\}_t={\frac {t^n-1} {t-1}}$, 
$\{n\}_t!=\{n\}_t\{n-1\}_t\cdots\{1\}_t$ and
\[{n\brace m}_t=
\left\{
\begin{array}{ll}
{\frac {\{n\}_t!} {\{m\}_t!\{n-m\}_t!}}
&\quad\mbox{if $n\leq m\leq 0$,} \\
0  & \quad\mbox{otherwise.}
\end{array}\right. \]
For $u \in h\uhgs$, put $\displaystyle{e(u,t)=
\sum_{n=0}^{\infty} {\frac {u^n} {\{n\}_t!}}}$. 
It is easy to show that\newline\par 
$e(-u,t^{-1})=e(u,t)^{-1}$,\hfill (7.3.1)
\par
$e(u,t)Xe(u,t)^{-1}=
\sum_{n=0}^{\infty}{\frac {1} {\{n\}_t!}} ad_{t^{n-1}}(u)ad_{t^{n-2}}(u)
\cdots ad_1(u)(X)$ \hfill (7.3.2)\newline\newline
where 
$ad_x(u)(X)$ $=[u,X]_{-,x}$ $=uX-xXu$.\par
For $\al\in\Pi$, let $\uhgas$ be a topological subalgebra of $\uhgs$
generated by $\uhhs$ and $E_\al, F_\al$. 
By the quantum double construction, we see that
$$
R_{\al}=e(-(q-q^{-1})E_\al\otimes F_\al\sigma^{p(\al)},
(-1)^{p(\al)}q^{(\al,\al)})\cdot R_T
\in \uhgas\otimes\uhgas 
$$
satisfies (7.1.1-2). Let $(\uhgs)^{(\al)}$
$=(\uhgs, \Da, S^{(\al)}, \varepsilon)$ be an another Hopf algebra 
defined as $((\uhgs)^{(R_\al)})^{(R_T^{-1})}$. Put 
$$
\hRa=e(-(q-q^{-1})\sigma^{p(\al)}K_\al^{-1}E_\al\otimes F_\al K_\al,
(-1)^{p(\al)}q^{(\al,\al)}).
$$
Then we get $\hRa=R_T^{-1}R_\al$. Hence
$$
\Da(X)=\hRa\Delta(X)\hRa^{-1} \quad (X\in \uhgs).
$$
{\bf Proposition 7.3.1.} {\it For $\al, \beta\in \Pi$.
Put $E^\vee_{\beta + s\al}=\whl...\whl\whl E_\beta, E_\al \whr ,
E_\al \whr \ldots E_\al \whr$, 
$F^\vee_{\beta + s\al}=\whl...\whl\whl F_\beta, F_\al \whr ,
F_\al \whr \ldots F_\al \whr$ 
($E_\al$, $F_\al$ appears $s$-times).
\newline\newline 
(i)\quad\,\, $\Da(E_\al K_\al^{-1})=
E_\al K_\al^{-1}\otimes K_\al + \sigma^{p(\al)}\otimes E_\al K_\al^{-1}$,
\par
\quad\,\, $\Da(K_\al F_\al)= K_\al F_\al \otimes 1 + 
\sigma^{p(\al)} K_\al^{-1} \otimes K_\al F_\al$.
\newline\newline
(ii) Assume $(\al,\al)\ne 0$. For $\beta\in\Pi$, assume
 $r=r_{(\al,\beta)}=-{\frac {2(\al,\beta)} {(\al,\al)}}
\in Z_+$ and $p(\al)\cdot r$ is even. Then
\[\begin{array}{c}
\Da(E^\vee_{\beta+r\al})=
E^\vee_{\beta+r\al}\otimes 1 + K_{\beta+r\al}\sigma^{p(\beta+r\al)}
\otimes E^\vee_{\beta+r\al}, \\
\Da(F^\vee_{\beta+r\al})=
F^\vee_{\beta+r\al}\otimes K_{\beta+r\al}^{-1} + \sigma^{p(\beta+r\al)}
\otimes F^\vee_{\beta+r\al}.
\end{array}\]
(iii) Assume $(\al,\al)=0$ and $(\al,\beta)\ne 0$. Then
\[\begin{array}{c}
\Da(E^\vee_{\beta+\al})=
E^\vee_{\beta+\al}\otimes 1 + K_{\beta+\al}\sigma^{p(\beta+\al)}
\otimes E^\vee_{\beta+\al}, \\
\Da(F^\vee_{\beta+\al})=
F^\vee_{\beta+\al}\otimes K_{\beta+\al}^{-1} + \sigma^{p(\beta+\al)}
\otimes F^\vee_{\beta+\al}.
\end{array}\]}
\newline
(iv) $\Da(H)=H\otimes 1 + 1\otimes H$, 
$\Da(\sigma)=\sigma\otimes\sigma$. \newline
{\bf Proof.} Here we calculate $\Da(E^\vee_{\beta+r\al})$ of (ii). 
Put $t_\al=(-1)^{p(\al)}q^{(\al,\al)}$
and $t_{\al,\beta}=(-1)^{p(\al)p(\beta)}q^{(\al,\beta)}$.
By direct calculation, we have:
\begin{eqnarray*}
\lefteqn{\Delta(E^\vee_{\beta+u\al})=E^\vee_{\beta+u\al}\otimes 1}\\
& & +\sum_{s=0}^u {u\brace s}_{t_\al} \prod_{k=1}^{u-s}
(t_{\al,\beta}^{-1}-t_{\al,\beta}^{-1}t_{\al}^{k-u})
E_\al^{u-s}K_{\beta+s\al}\sigma^{p(\beta+s\al)}
\otimes E^\vee_{\beta+s\al}.
\end{eqnarray*}
Hence, for $r$ of (ii),
\begin{eqnarray*}
\lefteqn{\Delta(E^\vee_{\beta+r\al})=E^\vee_{\beta+r\al}\otimes 1}\\
& & +\sum_{s=0}^r {r\brace s}_{t_\al}t_{\al,\beta}^{r-s}
\prod_{k=1}^{r-s}
(1-t_{\al}^k)
E_\al^{r-s}K_{\beta+s\al}\sigma^{p(\beta+s\al)}
\otimes E^\vee_{\beta+s\al}.
\end{eqnarray*} \par
Put $X=-\sigma^{p(\al)}K_\al^{-1}\otimes F_\al K_\al$. Then
\begin{eqnarray*}
\lefteqn {\bigl[-X,E_\al^s K_{\beta+(r-s)\al}\sigma^{p(\beta+(r-s)\al)}
\otimes E^\vee_{\beta+(r-s)\al}\bigr]_{-,t^{-s}}} \\
& & =
\left\{\begin{array}{l}
-{\frac {t_{\al,\beta}} {q-q^{-1}}} t_\al^{-s}
{\frac {(t_\al^{r-s}-1)(t_\al^{s+1}-1)} {t_\al-1}}
E_\al^{s+1} K_{\beta+(r-s-1)\al}\sigma^{p(\beta+(r-s-1)\al)}
\otimes E^\vee_{\beta+(r-s-1)\al} \\
\quad\quad \mbox{if $r>s$,}
\\
0\quad\quad \mbox{if $r=s$}.
\end{array}\right. 
\end{eqnarray*}
Hence
\begin{eqnarray*}
\lefteqn{ad_{t_\al^{-(s-1)}}(X)ad_{t_\al^{-(s-2)}}(X)\cdots ad_{1}(X)
(K_{\beta+r\al}\sigma^{p(\beta+r\al)}
\otimes E^\vee_{\beta+r\al})}\\
& & =
\left\{\begin{array}{l}
=\Bigl({\frac {t_{\al,\beta}} {q-q^{-1}}}\Bigr)^s 
t_\al^{-{{\frac {s(s-1)} {2}}}} (t_\al-1)^s
{r\brace s}_{t_\al}
E_\al^s K_{\beta+(r-s)\al}\sigma^{p(\beta+(r-s)\al)}
\otimes E^\vee_{\beta+(r-s)\al}\\
\quad \mbox{if $r\geq s$,} \\
0\quad \mbox{if $r <  s$.} 
\end{array}\right. 
\end{eqnarray*}
Hence, by (7.3.1-2), 
\begin{eqnarray*}
\lefteqn{\hRa^{-1}(K_{\beta+r\al}\sigma^{p(\beta+r\al)}
\otimes E^\vee_{\beta+r\al})\hRa}\\
& &
 =\sum_{s=0}^{r}(-t_\al)^s(t-1)^s
{\frac {\{r\}_{t_\al}!} {\{r-s\}_{t_\al}!}}
E_\al^s K_{\beta+(r-s)\al}\sigma^{p(\beta+(r-s)\al)}
\otimes E^\vee_{\beta+(r-s)\al}\,.
\end{eqnarray*}
On the other hand, we can easily show $\hRa^{-1}(E^\vee_{\beta+r\al}
\otimes 1)\hRa=E^\vee_{\beta+r\al}
\otimes 1$. Then we get:
$$
\hRa^{-1}(E^\vee_{\beta+r\al}\otimes 1+K_{\beta+r\al}\sigma^{p(\beta+r\al)}
\otimes E^\vee_{\beta+r\al})\hRa
=\Delta(E^\vee_{\beta+r\al})\,.
$$
We can show other formulae similarly or easily.\par
\hfill Q.E.D. 
\newline\newline
By [Y1], we see:\newline\newline
{\bf Lemma 7.3.1.} {\it For $\al\in\Pi$, $\uhg$ has an another Hopf
superalgebra structure $\uhg^{(\al)}=(\uhg, \Da_s)$
with coproduct $\Da_s$ satisfies formulae given by eliminating
$\sigma^{p(\cdot)}$ in the formulae of $\Da$ of Proposition 7.3.1
(i)-(iv).}
\newline\newline
{\bf 7.4.} {\bf Lemma 7.4.1.} {\it Keep notation in 7.3.
Let $\al\in\Pi$.\newline
(i) Assume $(\al,\al)\ne 0$. Assume $r=r_{\al, \beta}\in Z_+$  
and $p(\al)r \in 2Z$ for any $\beta\in \Pi\setminus \{\al\}$.
Define $\sigma_\al:{\cal H}\rightarrow {\cal H}$ by
$\sigma_\al(H_\lambda)=H_{\lambda
-{\frac {2(\al,\lambda)} {(\al,\al)}}\al}$.
Let $x_\beta, y_\beta \in C[[h]]^\times$ 
$(\beta\in\Pi)$
be such that
\[x_\beta y_\beta = \left\{
   \begin{array}{ll}
    (-1)^{p(\beta)} & (\beta=\al), \\
    (-1)^{r_{\al,\beta}}q^{-(r_{\al,\beta}+1)(\al,\beta)}
    \Bigl(
    {\frac {1-t_\al^{-1}} {q-q^{-1}}}
    \Bigr)^{r_{\al,\beta}} 
    (\{r_{\al,\beta}\}_{t_\al^{-1}}!)^2 
     & (\beta\ne\al,\,(\al,\beta)\ne 0), \\
    1 & (\beta\ne\al,\,(\al,\beta)= 0). 
  \end{array} \right. \]
Put 
$H'_\lambda =H_{\sigma_\al(\lambda)}$, 
$E'_{\al}=x_\al^{-1}F_\al K_\al$,
$F'_{\al}=y_\al^{-1}K_\al^{-1} E_\al$,
$E'_\beta=x_\beta^{-1}E^\vee_{\beta+r_{\al,\beta}\al}$,
$F'_\beta=y_\beta^{-1}F^\vee_{\beta+r_{\al,\beta}\al}$
$(\beta \in \Pi\setminus \{\al\})$. \par
Then $H'_\lambda$, $E'_\beta$, $F'_\beta$ satisfy (6.1.3).
\newline\newline
(ii) Assume $(\al,\al)=0$.
For $\beta\in\Pi\setminus \{\al\}$, put 
\[ r_{\al,\beta}=\left\{
 \begin{array}{ll}
  1 & (\al,\beta)\ne 0 \\
  0 & (\al,\beta)= 0 
  \end{array} \right. \]
Let 
$\Pi'=\{\al'=-\al,$ $\beta'=\al+\beta\, (\beta\in\Pi,\,
(\al,\beta)\ne 0),$ $\gamma'=\gamma \, (\gamma\in 
\Pi\setminus \{\al\},\,(\gamma,\al)=0)\}$.
Let $\sigma_\al:\Epip\rightarrow\Epipc$
by $\sigma_\al(H)=H$. (In particular,
\[ \sigma_\al(H_\beta)= \left\{
   \begin{array}{ll}
   H_{-\al'} & (\beta=\al), \\
   H_{\beta'+r_{\al,\beta}\al'} 
& (\beta\ne\al,\,(\al,\beta)\ne 0).) \\
  \end{array} \right. \]
Let $x_\beta, y_\beta \in C[[h]]^\times$ 
$(\beta\in\Pi)$
be such that
\[x_\beta y_\beta = \left\{
   \begin{array}{ll}
    -1       & (\beta=\al), \\
 t_{\al,\beta}^{-1}
{\frac {q^{(\al,\beta)}-q^{-(\al,\beta)}} {q-q^{-1}}}
     & (\beta\ne\al,\,(\al,\beta)\ne 0), \\
    1 & (\beta\ne\al,\,(\al,\beta)= 0). 
  \end{array} \right. \]
Put 
$E'_{\al}=x_\al^{-1}F_\al K_\al$,
$F'_{\al}=y_\al^{-1}K_\al^{-1} E_\al$,
$E'_\beta=x_\beta^{-1}E^\vee_{\beta+r_{\al,\beta}\al}$,
$F'_\beta=y_\beta^{-1}F^\vee_{\beta+r_{\al,\beta}\al}$
$(\beta \in \Pi\setminus \{\al\})$. \par
Then $H$, $E'_\beta$, $F'_\beta$ satisfy (6.1.3)
for $\Epipc$.}\newline
{\bf Proof.} Here we show how to calculate 
$$
[E^\vee_{\beta+r_{\al,\beta}\al},F^\vee_{\beta+r_{\al,\beta}\al}]
=x_\beta y_\beta {\frac {\sinh(hH_{\beta+r_{\al,\beta}\al})} {\sinh(h)}}.
\eqno{(7.4.1)}
$$ \par
Put $\{k;\beta\}_\al=
{\frac {q^{-(\al,\beta)}(1-t_\al^{-k})(1-t_{\al,\beta}^2 t_\al^{k-1})}
{(q-q^{-1})(1-t_\al^{-1})}}$
and
$\{k;\beta\}_\al !=\displaystyle{\prod_{v=1}^k \{v;\beta\}_\al}$.
First we show:
\[ \begin{array}{c}
\bigl[E_\al,F^\vee_{\beta+k\al}\bigr]
=-t_{\al,\beta}^{-1}q^{-(k-1)(\al,\al)}\{k;\beta\}_\al
F^\vee_{\beta+(k-1)\al}K_\al\,, \\
\bigl[E^\vee_{\beta+k\al},F_\al\bigr]=
(-1)^{(k-1)p(\al)}\{k;\beta\}_\al K_\al^{-1}E^\vee_{\beta+(k-1)\al}\,.
\end{array} \]
Then, by induction on $k$, we can show:
\begin{eqnarray*}
\lefteqn{\bigl[E^\vee_{\beta+k\al},F^\vee_{\beta+(k-1)\al}\bigr]} \\
& & =(-1)^{k(1+p(\al)p(\beta))}
q^{-{\frac {(k-1)(k-2)} {2}}(\al,\al)}
q^{(-k+1)(\al,\beta)} \{k;\beta\}_\al !
E_\al K_{\beta+(k-1)\al},
\end{eqnarray*}
\begin{eqnarray*}
\lefteqn{\bigl[E^\vee_{\beta+(k-1)\al},F^\vee_{\beta+k\al}\bigr]} \\
& & = (-1)^{(k-1)(1+p(\al)+p(\al)p(\beta))}
q^{-{\frac {k(k-1)} {2}}(\al,\al)}
q^{-k(\al,\beta)}
\{k;\beta\}_\al !
K_{\beta+(k-1)\al}^{-1} F_\al
\end{eqnarray*}
and
\begin{eqnarray*}
\lefteqn{\bigl[E^\vee_{\beta+k\al},F^\vee_{\beta+k\al}\bigr]} \\
& & = (-1)^k t_{\al,\beta}^{-k} 
q^{-{\frac {k(k-1)} {2}}(\al,\al)}
\{k;\beta\}_\al !
{\frac {\sinh(hH_{\beta+k\al})} {\sinh(h)}}.
\end{eqnarray*}
Substituting $r_{\al,\beta}$ for $k$, we get (7.3.1). \par
We can show other formulae similarly or easily.\par
\hfill Q.E.D.
\newline\newline
{\bf 7.5.} {\bf Proposition 7.5.1.} {\it Keep notations in 7.4. 
\newline (i)
Let
$\Pi^{\sigma_\al}=\Pi$ if $(\al,\al)\ne 0$ and let 
$\Pi^{\sigma_\al}=\Pi'$ if $(\al,\al)= 0$. 
Put $\uhgsa=U_h({\cal G}({\cal E}, \Pi^{\sigma_\al}, p))$. 
Then there is an isomorphism
$L_\al: \uhg\rightarrow\uhgsa$ such that 
$$
L_\al(H)=\sigma_\al(H),\,
L_\al(E_\beta)=E'_\beta,\,
L_\al(F_\beta)=F'_\beta. \eqno{(7.5.1)}
$$
(ii) 
$$
\Delta(L_\al(X))=\hRa^{-1}(L_\al\otimes L_\al \Delta(X))\hRa
\quad (X\in\uhgs). \eqno{(7.5.2)}
$$}
{\bf Proof.} (i) By Lemma 6.7.1, Lemma 7.3.1 and Lemma 7.4.1, there is an
 epimorphism $L'_\al: \uhg\rightarrow\uhgsa$ satisfying (7.5.1). 
Let $L_\al: \uhg\rightarrow\uhgsa$ denote $L'_\al$ defined by changing
 $\uhg$ and $\uhgsa$. ({\it Keep notations in the proof of
 Lemma 6.5.1.}) Since ${L'_\al L_\al}_{|{\cal H}}=id_{{\cal H}}$
 and (resp. ${L_\al L'_\al}_{|{\cal H}}=id_{{\cal H}}$), $L'_\al L_\al$
 (resp. $L_\al L'_\al$) induce an automorphism of ${\cal G}_0$
 (resp. ${\cal G}^{\sigma_\al}_0$) as well as an automorphism of
 $U_0({\cal G}_0)$
(resp. $U_0({\cal G}^{\sigma_\al}_0)$). Hence $L_\al$ induce an isomorphism
$U_0({\cal G}_0)\rightarrow U_0({\cal G}_0^{\sigma_\al})$.
Hence by topological
freedom of $\uhg$ and $\uhgsa$, $L_\al$ is an isomorphism. \newline
(ii) (7.5.2) is clear from the formulae in Proposition 7.3.1 (i)-(iii).
\par \hfill Q.E.D.
\newline\newline
By direct calculation, we have:\newline\newline
{\bf Lemma 7.5.1.} {\it $L_\al^{-1}:\uhg\rightarrow\uhgsa$ satisfies
(Here we let the region of definition (resp. values) of $L_\al^{-1}$
be $\Epip$ (resp. $({\cal E}, \Pi^{\sigma_\al}, p)$)).
For $\al, \beta\in \Pi$.
Put $E_{\beta + s\al}=\whl E_\al\ldots\whl
E_\al,\whl E_\al, E_\beta \whr\whr...\whr$, 
$F_{\beta + s\al}=\whl F_\al\ldots\whl
F_\al,\whl F_\al, F_\beta \whr\whr...\whr$, 
($E_\al$, $F_\al$ appears $s$-times).\par
 Put 
$H''_\lambda =H_{\sigma_\al(\lambda)}$, 
$E''_{\al}={\dot x}_\al^{-1}K_\al^{-1} F_\al$,
$F''_{\al}={\dot y}_\al^{-1}E_\al K_\al$,
$E''_\beta={\dot x}_\beta^{-1}E_{\beta+r_{\al,\beta}\al}$,
$F''_\beta={\dot y}_\beta^{-1}F_{\beta+r_{\al,\beta}\al}$
$(\beta \in \Pi\setminus \{\al\})$. Here we define
${\dot x}_\beta, {\dot y}_\beta\in C[[h]]^\times$ by:
\[\begin{array}{l}
x_\al {\dot y}_\al = y_\al {\dot x}_\al=1, \\
y_\al^{r_{\al,\beta}} y_\beta{\dot y}_\beta
= (-1)^{r_{\al,\beta}}(-1)^{(p(\al)+p(\al)p(\beta))r_{\al,\beta}}
\{r_{\al,\beta};\beta\}_\al ! , \\
x_\al^{r_{\al,\beta}} x_\beta{\dot x}_\beta
= (-1)^{r_{\al,\beta}}(-1)^{p(\al)p(\beta)r_{\al,\beta}}
q^{r_{\al,\beta}(\al,\al)}
\{r_{\al,\beta};\beta\}_\al ! \quad
(\beta\ne\al,\,(\al,\beta)\ne 0), \\
x_\beta {\dot x}_\beta = y_\beta {\dot x}_\beta=1
\quad (\beta\ne\al,\,(\al,\beta)= 0).
\end{array} \]
(Here $x_\beta, y_\beta\in C[[h]]^\times$
have been defined in Lemma 7.4.1 for 
$L_\al:\uhgsa\rightarrow\uhg$.)
}\newline\newline
{\bf 7.6.} As an immediate consequence of Proposition 7.5.1, we have:
\newline\newline
{\bf Proposition 7.6.1.} {\it (See also [KT].)} {\it Let $\Epip$'s be the data of affine type.
 For the isomorphisms $L_i$ defined for 
${\cal G}\Epip$'s in \S 2, there are isomorphisms $T_i$'s of
$U_h(\g\Epip)$'s such that 
$T_i\rightarrow L_i:U_0(\g\Epip)\rightarrow U_0(\g\Epipsi)$ ($h\rightarrow 0$).}
\newline\newline
{\bf 8. On $U_h({\hat {sl}}(m|m))^{(i)}$ $(i=1,2,\,4)$.}\newline 
In this chapter, we use Beck's method [B]. \newline
{\bf 8.1.} In 8.1, let $\Epip$ be of Diagram 1.6.2 and assume $N\geq 4$.
Let $W$ be the Weyl group defined in 2.6 associated to $\Epidagger$. 
Let $W_0$ be a subgroup of W generated by $\{\sigma(i),\, (1\leq 
i\leq n)\}$. Let $\om^\vee_j \in \oplus_{i=1}^n C\al_i^\dagger$
 $(1\leq j \leq n)$ be such that
${\frac {2((\al_i^\dagger,\om^\vee_j))} {((\al_i^\dagger,\al_i^\dagger))}}
=\delta_{ij}$ $(1\leq i \leq n)$. Put $P^\vee=\oplus Z\om^\vee_i$. Define
${\bar W}=W_0 |\!\!\!\times P^\vee$ by $(s,x)(s',x')=
(ss',s'^{-1}(x)+x')$. We know that there is a certain subgroup ${\cal T}$ of
Dynkin diagram automorphism of $\Epidagger$ such that 
${\bar W} \cong {\cal T} |\!\!\!\times W_0$ $(\tau \sigma(i) \tau^{-1}
=\sigma(\tau(i))\,\, (\tau\in{\cal T}))$. If $W$ is of type $A_{N-1}^{(1)}$, then
${\cal T}\cong Z/NZ$. For the datum $(
{\cal E}=(\oplus_{i=1}^N C\be_i)\oplus C\delta\oplus C\lambda_0,
\Pi=\{\al_i\},p)$ and $\tau\in {\cal T}$, define the datum 
$({\cal E}^\tau=(\oplus_{i=1}^N C\be_i^\tau)\oplus C\delta
\oplus C\lambda_0,\Pi^\tau =\{\al_i^\tau\},p^\tau)$ by\newline\newline
(i) The Dynkin diagram of $\Epipt$ is the same type as the one of
$\Epip$. \newline
(ii) $p^\tau(\al_i)=p(\al_{\tau^{-1}(i)})$. \newline
(iii) $(\be_i^\tau,\be_j^\tau)=(\be_{\tau^{-1}(i)},
\be_{\tau^{-1}(j)})$ \,\, (Here we consider $\tau^{-1}(i)$ under $\mbox
{mod} N$).\newline\newline
For $w\in {\bar W}$ and 
an reduced expression  $w=\tau\sigma(i_1)\cdots\sigma(i_r)$,
put $\Epipw$
$=(\,((\,({\cal E}^{\sigma(i_r)})\cdots)^{\sigma(i_1)})^\tau,$
$((\,(\Pi^{\sigma(i_r)})\cdots)^{\sigma(i_1)})^\tau,$
$((\,(p^{\sigma(i_r)})\cdots)^{\sigma(i_1)})^\tau)$.
Clearly $\Epipw$ doesn't depend on reduced expressions.
\par Let $\uhg'$ be the subalgebra of $\uhg$ generated by 
$\{ H_{\al_i},\,E_i,\, F_i\,\, (0\leq i \leq n)\}$.
By Proposition 7.5.1, Lemma 7.5.1 and direct calculation, we have:\newline
\newline
{\bf Lemma 8.1.1.} {\it Assume that $\Epip$ is type $A_{N-1}^{(1)}$. 
Keep notations in 2.3-5. 
Let $i\in \{0,1,...,N-1\}(=Z/NZ)$. Put $K_i=K_{\al_i}$.
\newline (i) There are isomorphisms
$T_i:U_h({\cal G}\Epip)'\rightarrow U_h({\cal G}\Epipsi)'$ such that
(We put $p'=p^{\sigma(i)}$.) 
\[ \begin{array}{l}
T_iE_i=-\bd'_{i+1}F_iK_i,\,\, T_iF_i=-\bd'_iK_i^{-1}E_i, \\
T_iE_{i-1}=q^{-\bd'_i}\bd'_i\whl E_{i-1}, E_i \whr,\,
T_iE_{i+1}=q^{-\bd'_{i+1}}
(-1)^{p'(\al_i)p'(\al_{i+1})}\bd'_{i+1}\whl E_{i+1}, E_i \whr, \\
T_iF_{i-1}=-(-1)^{p'(\al_i)p'(\al_{i-1})}\whl F_{i-1}, F_i \whr,\,
T_iF_{i+1}=-\whl F_{i+1}, F_i \whr.
\end{array}\]
$T_i^{-1}:U_h({\cal G}\Epip)'\rightarrow U_h({\cal G}\Epipsi)'$ is given by:
\[ \begin{array}{l}

T_i^{-1}E_i=-\bd'_{i+1}K_i^{-1}F_i,\,\, T_i^{-1}F_i=-\bd'_iE_iK_i, \\
T_i^{-1}E_{i-1}=q^{-\bd'_i}(-1)^{p'(\al_i)p'(\al_{i-1})}\bd'_i\whl E_i, E_{i-1} \whr,\,
T_i^{-1}E_{i+1}=q^{-\bd'_{i+1}}
\bd'_{i+1}\whl E_i, E_{i+1} \whr, \\
T_i^{-1}F_{i-1}=-\whl F_i, F_{i-1} \whr,\,
T_i^{-1}F_{i+1}=-(-1)^{p'(\al_i)p'(\al_{i+1})}\whl F_i, F_{i+1} \whr.
\end{array}\]
For $\tau\in{\cal T}$, there is an isomorphism
$T_\tau:U_h({\cal G}\Epip)'\rightarrow U_h({\cal G}\Epipt)'$ such that
$T_\tau(H_{\al_i})=H_{\al_{\tau(i)}}$, 
$T_\tau(E_i)=E_{\tau(i)}$, $T_\tau(F_i)=F_{\tau(i)}$.
\newline\newline 
(ii) $T_i$'s satisfy Braid relation:
$$
T_iT_j=T_jT_i\,((\al_i,\al_j)=0),\quad 
T_iT_jT_i=T_jT_iT_j\,(|(\al_i,\al_j)|=1).
$$ It also hold that $T_\tau T_i T_\tau^{-1} = T_{\tau(i)}$. 
\newline\newline
(iii) By (ii), putting $T_w=T_\tau T_{i_1}\cdots T_{i_r}$ for
$w\in {\bar W}$ whose reduced expression is
  $w=\tau\sigma(i_1)\cdots\sigma(i_r)$,
$T_w$ is well-defined. Moreover we have: 
$$
T_w(E_i)=E_j,\, T_w(F_i)=F_j\quad
\mbox{if $w(\al_i)=\al_j$}.
$$ There
 is an $C$-anti-automorphism $\Omega$ such that 
$$
\Omega(E_i)=\bd_{i+1}F_i,\, \Omega(E_i)=\bd_{i+1}F_i,\, 
\Omega(H)=H,\, \Omega(h)=-h\,.
$$ Moreover
$\Omega T_w = T_w\Omega$ $(w\in {\bar W})$.}\newline\newline
{\bf 8.2.} Put $T_{\om_i}=T_{\om^\vee_i}$. For
$1\leq i \leq n$, $k > 0$ and $s\in Z$, let
$$
\bps^{(s)}_{ik}=K_\delta^{-{\frac {k} {2}}}q^{-(\al_i,\al_i)}
\whl T_{\om_i}^s(E_i), T_{\om_i}^{k+s}(K_i^{-1}F_i) \whr\,.
\eqno{(8.1.1)}
$$
Put $Q_{ij,k}={\frac {q^{k(\al_i,\al_j)}-q^{-k(\al_i,\al_j)}}
{q-q^{-1}}}$ and ${\dot C}_{ij}=q^{(\al_i,\al_j)}K_\delta^{-{\frac {1} {2}}}$.
By [B], we have: \newline\newline
{\bf Lemma 8.2.1.} {\it (i) $K_\delta^{\frac {1} {k}}\bps^{(s)}_{ik}\in{\cal N}_+$
if $s\leq 0$ and $k+s>0$. \newline
(ii) Assume $p(\al_i)=0$. Let $r>0$ and $m\in Z$. Then $\bps^{(s)}_{ir}
=\bps^{(s')}_{ir}$ $(s, s' \in Z)$ and 
\begin{eqnarray*}
\bigl[\bps^{(s)}_{ir},\Tomi^m(F_i)\bigr] =
-K_\delta^{{\frac {1} {2}}}Q_{ii,1}\Bigl\{((q-q^{-1})
\sum_{k=1}^{r-1} {\dot C}_{ii}^{1-k}
\Tomi^{m+k}(F_i)\bps^{(s)}_{i,r-k}) +
{\dot C}_{ii}^{1-r}
\Tomi^{m+r}(F_i)\Bigr\}, \\
\bigl[\bps^{(s)}_{ir},\Tomi^m(E_i)\bigr]=
K_\delta^{-{\frac {1} {2}}}Q_{ii,1}\Bigl\{((q-q^{-1})
\sum_{k=1}^{r-1} {\dot C}_{ii}^{k-1}
\Tomi^{m-k}(E_i)\bps^{(s)}_{i,r-k}) +
{\dot C}_{ii}^{r-1}\Tomi^{m-r}(E_i)\Bigr\} .
\end{eqnarray*}
(ii) Assume $1\leq i\ne j \leq n$. Let $r>0$ and $m\in Z$.
Then:
\[\begin{array}{l}
\bigl[\bps^{(s)}_{ir},\Tomj^m(F_j)\bigr] \\
\quad = \displaystyle{
K_\delta^{{\frac {1} {2}}}Q_{ij,1}\Bigl\{((q-q^{-1})
\sum_{k=1}^{r-1} (-{\dot C}_{ij})^{1-k}
\Tomj^{m+k}(F_j)\bps^{(s)}_{i,r-k}) +
(-{\dot C}_{ij})^{1-r}
\Tomj^{m+r}(F_j)\Bigr\}}, \\
\bigl[\bps^{(s)}_{ir},\Tomj^m(E_j)\bigr] \\
\quad = \displaystyle{
-K_\delta^{-{\frac {1} {2}}}Q_{ij,1}\Bigl\{((q-q^{-1})
\sum_{k=1}^{r-1} (-{\dot C}_{ij})^{k-1}
\Tomj^{m-k}(E_j)\bps^{(s)}_{i,r-k}) +
(-{\dot C}_{ij})^{r-1}\Tomj^{m-r}(E_j)\Bigr\}} .
\end{array}\]} 

Let $o(i)\in\{\pm 1\}$ satisfy that $o(i)\ne o(j)$ if $(\al_i,\al_j)\ne 0$
$(i\ne j)$. Put $\hTomi^mE_i=o(i)^m\Tomi^mE_i$
and $\hTomi^mF_i=o(i)^m\Tomi^mF_i$. Define $\hbps^{(s)}_{ir}$ by replacing
$\Tomi$ of (8.1.1) with $\hTomi$. By Lemma 8.2.1, we have:
\newline\newline
{\bf Lemma 8.2.2.} {\it Assume $j\ne i$ or $p(\al_i)=0$. Then:

\[\begin{array}{l}
\bigl[\hbps^{(s)}_{ir},\hTomj^m(F_j)\bigr] \\
\quad = \displaystyle{-
K_\delta^{{\frac {1} {2}}}Q_{ij,1}\Bigl\{((q-q^{-1})
\sum_{k=1}^{r-1} {\dot C}_{ij}^{1-k}
\hTomj^{m+k}(F_j)\hbps^{(s)}_{i,r-k}) +
{\dot C}_{ij}^{1-r}
\hTomj^{m+r}(F_j)\Bigr\}}, \\
\bigl[\hbps^{(s)}_{ir},\hTomj^m(E_j)\bigr] \\
\quad = \displaystyle{
K_\delta^{-{\frac {1} {2}}}Q_{ij,1}\Bigl\{((q-q^{-1})
\sum_{k=1}^{r-1} {\dot C}_{ij}^{k-1}
\hTomj^{m-k}(E_j)\hbps^{(s)}_{i,r-k}) +
{\dot C}_{ij}^{r-1}\hTomj^{m-r}(E_j)\Bigr\}} .
\end{array}\]}\newline\newline
Define $h_{ik}^{(s)}\in\uhg$ $(k>0)$ by the following generating function in $z$.
$$
\exp((q-q^{-1})\sum_{k=1}^{\infty}h_{ik}^{(s)}z^k)
=1+(q-q^{-1})\sum_{k=1}^{\infty}\hbps_{ik}^{(s)}z^k.
$$
{\it Remark.} For $(\al_i,\al_i)=0$, we have not shown
 $[\hbps^{(s)}_{ik},\hbps^{(s)}_{ir}]=0$ yet.
Hence an uncertainly of the definition of $h_{ik}^{(s)}$ has still remined.
It depends on an order of $\{\hbps^{(s)}_{ik}\}$.
 \newline\newline
By Lemma 8.2.2, we have:\newline\newline
{\bf Lemma 8.2.3.} {\it Assume $j\ne i$ or $(\al_i,\al_i)\ne 0$. Then:

\[\begin{array}{l}
\bigl[h^{(s)}_{ik},\hTomj^m(F_j)\bigr] 
\quad =-{\frac {1} {k}}Q_{ij,k}K_\delta^{\frac {k} {2}}
\hTomj^{m+k}(F_j), \\
\bigl[h^{(s)}_{ik},\hTomj^m(E_j)\bigr] 
\quad ={\frac {1} {k}}Q_{ij,k}K_\delta^{-{\frac {k} {2}}}
\hTomj^{m-k}(E_j).
\end{array}\]}\newline\newline
{\bf 8.3.} {\bf Lemma 8.3.1.} {\it Let $1\leq i\leq n$ and $r\in Z$. Then:
\[\begin{array}{l}
\whl \Tomi^{m+r}(F_i), \Tomi^{m}(F_i) \whr
= -\whl \Tomi^{m+1}(F_i), \Tomi^{m+r-1}(F_i) \whr , \\
\whl \Tomi^{m}(E_i), \Tomi^{m+r}(E_i) \whr
= -\whl \Tomi^{m+r-1}(E_i), \Tomi^{m+1}(E_i) \whr .
\end{array}\]} \newline
{\bf Proof.} For $(\al_i,\al_i)\ne 0$, we have already known these by [B].
For $(\al_i,\al_i)=0$, by Lemma 8.2.3 and $\Tomi^m(F_i)^2=\Tomi^m(E_i)^2=0$,
$$
[\Tomi^{m+r}(F_i), \Tomi^{m}(F_i)]=[\Tomi^{m+r}(E_i), \Tomi^{m}(E_i)]=0,
\eqno{(8.3.1)}
$$ which are nothing else but the formulae we want. \par
\hfill Q.E.D.\newline\newline
{\bf Lemma 8.3.2.} {\it Let $(\al_i,\al_i)=0$ and $r>0$. 
Then $\bps_{ir}^{(s)}=\bps_{ir}^{(s')}$ and 
$$
[h_{ir}^{(s)},\Tomi^m(F_i)]=[h_{ir}^{(s)},\Tomi^m(E_i)]=0.
\eqno{(8.3.2)} $$}\newline
{\bf Proof.} By (8.3.1), we have:
$$
[E_i, \bps_{ir}^{(0)}]=[\bps_{ir}^{(0)},F_i]=0.
\eqno{(8.3.3)}
$$ We use an induction on $r$. Let $1\leq j\leq n$ be such that 
$(\al_i,\al_j)\ne 0$. First we assume $r=1$. Then $h_{i1}^{(s)}=
o(i)\bps_{i1}^{(s)}$. 
\begin{eqnarray*}
 \bps_{i1}^{(-1)}&=&K_\delta^{-{\frac {1} {2}}}[\Tomi^{-1}(E_i),K_i^{-1}F_i] \\
&=& o(i)K_\delta^{-{\frac {1} {2}}} Q_{ji,1}^{-1}K_\delta^{{\frac {1} {2}}}\
[[h_{j1}^{(0)} ,E_i],K_i^{-1}F_i]\quad \mbox{(by Lemma 8.2.3)} \\
&=& o(i)Q_{ji,1}^{-1}K_i^{-1}\cdot o(i)Q_{ji,1}[K_\delta^{-{\frac {1} {2}}}
\Tomi(F_i), E_i] \\
&=& K_\delta^{-{\frac {1} {2}}}[\Tomi(K_i^{-1}F_i), E_i] \\ 
&=&\bps_{i1}^{(0)}.
\end{eqnarray*} Hence, by (8.3.3), we get our formulae for $r=1$. \par
We assume that we have shown the lemma for 1, 2,..., $r-1$.
Firstly we show $[h_{j1}^{(0)},h_{jr-1}^{(0)}]=0$. By Lemma
8.2.3, we have:
$$
[[h_{j1}^{(0)},h_{jr-1}^{(0)}],\hTomk^m(F_k)]=
[[h_{j1}^{(0)},h_{jr-1}^{(0)}],\hTomk^m(E_k)]=0
$$ for $1\leq k\leq n$ and $m\in Z$. By Lemma 8.2.1 (i), 
$K_\delta^{\frac {r} {2}}[h_{j1}^{(0)},h_{jr-1}^{(0)}]\in N^+$. We know the
 fact that $\hTomk^m(F_k)$, $\hTomk^m(E_k)$ and ${\cal H}$ generate $\uhg$.
Hence, by Proposition 6.2.1, we get 
$$
[h_{j1}^{(0)},h_{jr-1}^{(0)}]=0 \eqno{(8.3.4)}
$$  as well as $[h_{j1}^{(0)},\bps_{jr-1}^{(0)}]=0$. Hence:
\begin{eqnarray*}
 \bps_{ir}^{(-1)}&=&K_\delta^{-{\frac {1} {2}}}[\Tomi^{-1}(E_i),
\Tomi^{r-1}(K_i^{-1}F_i)] \\
&=& o(i)K_\delta^{-{\frac {r} {2}}} Q_{ji,1}^{-1}K_\delta^{{\frac {1} {2}}}
[[h_{j1}^{(0)} ,E_i],\Tomi^{r-1}(K_i^{-1}F_i)]\quad \mbox{(by Lemma 8.2.3)} \\
&=& o(i)Q_{ji,1}^{-1}K_\delta^{{\frac {1-r} {2}}}
\Bigl\{[h_{j1}^{(0)},K_\delta^{{\frac {r-1} {2}}}\bps_{ir-1}^{(0)}]+
o(i)Q_{ji,1}[K_\delta^{-{\frac {1} {2}}}\Tomi^r(K_i^{-1}F_i),E_i]\Bigr\} \\
&=&\bps_{ir}^{(0)}.
\end{eqnarray*} Hence, by (8.3.3), we get our formulae. \par
\hfill Q.E.D. \newline\newline
We put $h_{ir}=h_{ir}^{(s)}$ and $\hbps_{ik}\hbps^{(s)}_{ik}$ 
$(r>0,\,1\leq i\leq n)$ which is well defined by
Lemma 8.3.2. Similarly to show (8.3.4), we have:
\newline\newline
{\bf Lemma 8.3.3.} {\it $[h_{ir},h_{ir'}]=0$.}
\newline\newline By Lemma 8.2.3 and Lemma 8.3.2, we have:
\newline\newline
{\bf Lemma 8.3.4.} {\it    
\[\begin{array}{l}
\bigl[h_{ik},\hTomj^m(F_j)\bigr] 
\quad =-{\frac {1} {k}}Q_{ij,k}K_\delta^{\frac {k} {2}}
\hTomj^{m+k}(F_j), \\
\bigl[h_{ik},\hTomj^m(E_j)\bigr] 
\quad ={\frac {1} {k}}Q_{ij,k}K_\delta^{-{\frac {k} {2}}}
\hTomj^{m-k}(E_j).
\end{array}\]}
\newline\newline
We know that $K_\delta^{\frac {k} {2}}h_{ik}\in N^+$.
By Lemma 8.3.4 and Proposition 6.2.1, since $\hTomj^m(E_j)$,
 $\hTomj^m(F_j)$ and ${\cal H}$ generate $\uhg$, we have:
\newline\newline
{\bf Lemma 8.3.5.} {\it Keep notations in 1.6. If $\displaystyle{
\sum_{i=1}^{N}\bd_i=0}$, then
$$
\sum_{i=1}^{n}\Bigl[\sum_{j=1}^{i}k\bd_j\Bigr]K_\delta^{\frac {k} {2}}h_{ik}=0\quad (k\geq 1).
\eqno{(8.3.5)}
$$ in ${\cal N}^+\subset \uhg\Epip$ of $\Epip$ of  Diagram 1.6.2 $(N\geq 4)$.}
\newline\newline After all we have:
\newline\newline
{\bf Theorem 8.3.6.} {\it Let $\Epip$ be the datum of Diagram 1.6.2 
$(N\geq 4)$ with $\sum \bd_i=0$. Then the defining relations of $U_h(\g\Epip)$ 
are defined by adding (8.3.5) 
to the ones of Theorem 6.6.1.}\newline\newline
{\bf 8.4.} For $1\leq i\leq N-1$ and $r\geq 0$, put
\[ \psi_{ir}=\,\, \left\{
 \begin{array}{ll}
 (q-q^{-1})K_i\hbps_{ik}  &\quad(r>0),\\
 K_i  &\quad(r=0).\\
 \end{array}\right.\] and $\varphi_{ir}=\Omega(\psi_{ir})$.
Put $h_{i,-r}=h_{ir}$ $(r>0)$.
For $1\leq i\leq N-1$ and $k\in Z$, put $x^{-}_{ik}=\hTomi^k(F_i)$ and
$x^{+}_{ik}=\hTomi^{-k}(E_i)$.
Similar to [B], we have:\newline\newline
{\bf Theorem 8.4.1.} {\it Let $\Epip$ be the datum of Diagram 1.6.2 
$(N\geq 3)$. Then $U_h(\g\Epip)$ is defined with the generators $\{
H\in{\cal H},\,x^\pm_{ij},\,h_{ik}\}$ and the relations:
\[\begin{array}{c}
[H,x^\pm_{jk}]=(\pm\al_j+k\delta)(H)x^\pm_{jk}, \\
\big[h_{ik},h_{jl}\big]=\delta_{k,-l}{\frac {1} {k}}Q_{ij,k}
{\frac {K_\delta^k-K_\delta^{-k}} {q-q^{-1}}},\\
x^\pm_{ik+1}x^\pm_{jl}-(-1)^{p(\al_i)p(\al_j)}q^{\pm(\al_i,\al_j)} 
x^\pm_{jl}x^\pm_{ik+1}=
(-1)^{p(\al_i)p(\al_j)}q^{\pm(\al_i,\al_j)}x^\pm_{ik}x^\pm_{jl+1}-
x^\pm_{jl+1}x^\pm_{ik}, \\
\big[x^{+}_{ik},x^{-}_{jl}\big]=\delta_{ij}{\frac {
K_\delta^{\frac {k-l} {2}}\psi_{ik+l}-K_\delta^{\frac {l-k} {2}}\phi_{ik+l}}
{q-q^{-1}}}, \\
\end{array}\]  $[x^\pm_{ik}, x^\pm_{il}]=0$
\hfill $\mbox{if}\, (\al_i,\al_i)=0,$ 
\newline\newline 
(In the following equations, $\mbox{Sym}_{k_1,k_2,\ldots,k_s}$
means symmetrization with respect to $\{k_1,k_2,\ldots,k_s\}$.)
\newline\newline 
$\mbox{Sym}_{k_1,k_2}\whl x^\pm_{ik_1},
 \whl x^\pm_{ik_2}, x^\pm_{jl}\whr\whr=0$
\hfill $\mbox{if}\, (\al_i,\al_i)\ne 0\,\,\mbox{and}\,\,(\al_i,\al_j)\ne 0,$ 
\newline\newline
$\mbox{Sym}_{k_1,k_2}
[\whl \whl x^\pm_{il}, x^\pm_{jk_1} \whr ,x^\pm_{um} \whr , x^\pm_{jk_2}]=0$ \hfill  \aivu 
\newline\newline
(Each of the following equations means an equation as a generate function
in an indeterminate $z$.)\newline
\[\begin{array}{c}
\sum_{k\geq 0}\psi_{ik}z^k=K_i\exp((q-q^{-1})\sum_{r\geq 1}h_{ir}z^r), \\ 
\sum_{k\geq 0}\phi_{ik}z^k=K^{-1}_i\exp((q^{-1}-q)\sum_{r\geq 1}h_{i,-r}z^r).
 \\ \end{array} \]

}


\begin{thebibliography} {}
 \bibitem[D1]{}V.G.Drinfeld.
{\it Quantum groups},
Proc. Int. Congr. Math., Berkeley, 1 (1986);
Amer. Math. Soc. (1988), 798-820.
 \bibitem[D2]{}V.G.Drinfeld.
{\it On almost cocommutative Hopf algebras}, 
Leningrad Math. J., 1(1990), 321-342
 \bibitem[FSS]L.Frappet-A.Seiarino-P.Sorba. {\it
Structure of basic Lie superalgebras and of their affine extensions},
Commun. Math. Phys. 121(1989) 457-500.
 \bibitem[GK]{}O.Gabber-V.G.Kac.
{\it On defining relations of certain infinite-dimensional
Lie algebras}, Bull. Amer. Math. Soc. 5(1981), 185-189.
 \bibitem[J]{}M.Jimbo.
{\it A q-defference analogue of U(G) and the Yang-Baxter equation},
Lett. Math. Phys., 10(1985), 798-820.
 \bibitem[K1]{}V.G.Kac.
{\it Infinite dimensional Lie algebras},
3rd ed., Cambridge Univ. Press, Cambridge, 1990.
 \bibitem[K2]{}V.G.Kac,
{\it Lie superalgebras}, 
Adv. Math., 26 (1977), 8-96.
 \bibitem[KT]{}S.M.Khoroshkin-V.N.Tolstoy.
{\it Twisting of quantum (super)algebras. Connection of Drinfeld's and 
Cartan-Weyl realization for quantum affine algebras}, preprint (1993)
 \bibitem[MM]{}J.W.Milnor-J.C.Moore.
 {\it  On the structure of Hopf algebras}, 
Ann. Math. 81(1965) 211-264.
 \bibitem[LS]{}D.Leites-V.Serganava. {\it
Defining relations for classical Lie superalgebras I},
preprint (1991).
 \bibitem[VdL]{}J.W. Van De Leur.
{\it Contragrredient Lie superalgebras of finite growth},
preprint (1986)
 \bibitem[Y1]{}H.Yamane.
 {\it Quantized enveloping algebras associated to simple Lie superalgebras and
Their R-matrices},
Publ. RIMS Kyoto Univ., 30, no. 1 (1994) 15-87.
 \bibitem[Y2]{}H.Yamane.
 {\it A Serre type theorem for affine Lie superalgebras and their 
quantized enveloping superalgebras},
Proc. Japan Acad.,70 Ser.A (1994) 31-36.
\end{thebibliography}
\end{document}